\documentclass[nonblindrev]{informs3_noinforms}              

\usepackage[numbers]{natbib}

\usepackage{enumerate}
\newtheorem{observation}{Observation}

\newtheorem{Assumption}{Assumption}

\newcommand{\E}{\mathbb{E}}

\newcommand{\ignore}[1]{\relax}

\newcommand{\SKOR}{\text{\textsc{SKOR}}}
\newcommand{\CLOCK}{\text{\textsc{CLOCK}}}
\newcommand{\RES}{\text{\textsc{RES}}}
\newcommand{\REN}{\text{\textsc{REN}}}
\newcommand{\ARR}{\text{\textsc{ARR}}}
\newcommand{\DMY}{\text{\textsc{DMY}}}
\newcommand{\ABA}{\text{\textsc{ABA}}}
\newcommand{\ERR}{\text{\textsc{ERR}}}
\newcommand{\OU}{\text{\textsc{OU}}}
\newcommand{\IOU}{\text{\textsc{IOU}}}
\newcommand{\BUF}{\text{\textsc{BUF}}}
\newcommand{\TRIP}{\text{\textsc{TRIP}}}
\newcommand{\PAT}{\text{\textsc{PAT}}}
\newcommand{\VRW}{\text{\textsc{VRW}}}

\renewcommand{\theequation}{\thesection.\arabic{equation}}

\usepackage{natbib}
 \NatBibNumeric
 \def\bibfont{\small}%
 \bibpunct[, ]{[}{]}{,}{n}{}{,}%

\usepackage[colorlinks=true,breaklinks=true,bookmarks=true,urlcolor=blue,
     citecolor=blue,linkcolor=blue,bookmarksopen=false,draft=false]{hyperref}

\def\EMAIL#1{\href{mailto:#1}{#1}}

\TheoremsNumberedThrough     

\EquationsNumberedThrough    

\MANUSCRIPTNO{} 

\begin{document}



\RUNTITLE{Large deviations for the $M/H_2/n + M$ queue in the Halfin-Whitt regime}

\TITLE{Large deviations analysis for the $M/H_2/n + M$ queue in the Halfin-Whitt regime}

\ARTICLEAUTHORS{%
\AUTHOR{\bf Debankur Mukherjee}
\AFF{Eindhoven University of Technology, \EMAIL{d.mukherjee@tue.nl}}
\AUTHOR{\bf Yuan Li}
\AFF{Amazon, \EMAIL{yuan.leo.li@gmail.com}}
\AUTHOR{\bf David A. Goldberg}
\AFF{Cornell University, \EMAIL{dag369@cornell.edu}}
}
\RUNAUTHOR{D. Mukherjee, Y. Li, and D.A. Goldberg}
\ABSTRACT{%
We consider the FCFS $M/H_2/n + M$ queue in the Halfin-Whitt heavy traffic regime, i.e. the first-come-first-served multi-server queue with abandonments, in which inter-arrival and patience times are Markovian, and service times are hyper-exponentially distributed, i.e. the mixture of two exponential distributions.  It is known that if one considers a sequence of such queues under the Halfin-Whitt scaling, the appropriately normalized sequence of steady-state queue-length distributions is tight and converges weakly to a limiting random variable $W^{\infty}$.  However, those works only describe $W^{\infty}$ implicitly as the unique invariant measure of a certain complicated diffusion process with no closed form.  Although it was proven by Gamarnik and Stolyar (2012) that $P\big( W^{\infty} > x \big)$ has a sub-Gaussian decay, and related Lyapunov-type arguments and bounds also appear in the work of Dieker and Gao (2013), the actual value of $\lim_{x \rightarrow \infty} x^{-2} \log\bigg( P\big( W^{\infty} > x \big) \bigg)$ was left open.  In subsequent work, Dai and He (2013) conjectured an explicit form for this exponent, which depends only on the first two moments of the service distribution, and used this conjectured exponent to devise a numerical method for computing the tail probabilities of $W^{\infty}$.  
\\\indent In this work, we explicitly compute the true large deviations exponent for $W^{\infty}$, in the parameter regime such that the abandonment rate is strictly less than the minimum exponential service rate.  This is the first such result for multi-server queues with abandonments when service times are non-Markovian.  As this exponent is strictly less than that conjectured by Dai and He, and further does not only depend on the first two moments of the service distribution, we resolve the conjecture of Dai and He in the negative.  Our main approach is to extend the stochastic comparison framework of Gamarnik and Goldberg (2013) to the setting of abandonments, which requires several novel and non-trivial contributions (especially in light of the non-monotonicities introduced by abandonments), and to combine with results from the theory of Gaussian processes.  Furthermore, our approach sheds light on several novel ways to think about multi-server queues with abandonments in the Halfin-Whitt regime, e.g. by formally connecting to so-called multi-dimensional Ornstein-Uhlenbeck processes (a family of tractable diffusions), and we believe that these insights should hold in considerable generality and provide a novel set of tools for analyzing these systems.
}
\KEYWORDS{many-server queues, Halfin-Whitt regime, stochastic comparison, weak convergence, large deviations, Gaussian process, abandonments}
\maketitle
\newpage
\hypersetup{
  colorlinks,
  linkcolor=blue,
  linktoc=page
}
\renewcommand*\contentsname{Table of contents.\\}
\tableofcontents
\newpage
\section{Introduction.}\label{IntroSec}
\subsection{Multi-server queues with abandonments in the Halfin-Whitt regime.}
It is by now well-known that from a queueing perspective, a service system can operate in several different asymptotic regimes, parametrized by the number of servers and traffic intensity.  Recently there has been a significant interest in the so-called Halfin-Whitt (i.e. HW, Quality-and-Efficiency Driven, QED) regime, in which one can easily trade-off between the number of servers and the quality of service experienced by arrivals to the system (as parametrized through metrics including the probability that an arriving job must wait for service).  This scaling regime was studied originally by \citet{E.48} and \citet{J.74}, and formally introduced by \citet{HW.81}, who studied the First-come-first-serve (FCFS) $G/M/n$ system (for large $n$) when the traffic intensity $\rho$ scales like $1 - B n^{-\frac{1}{2}}$ for some strictly positive excess parameter $B$.  We refer the interested reader to \cite{GG13,G16,G17a,AR16,Leeuw.17} for an overview of recent progress on our understanding of multi-server queues in the HW regime, which we make no attempt to survey here, and to \cite{G17a} for an overview of how this scaling regime relates to other regimes and related results.  
\\\indent In many such service systems customer abandonment is an important feature \cite{Ward.12} (see also several relevant empirical works including \cite{Zeltyn.13,Brown.05,Zelt.04b,Batt.15}) and there is also considerable work studying multi-server queueing models with abandonments in the HW regime.  For an overview of work on multi-server queues with abandonments in the HW regime, as well as other scaling regimes, we refer the interested reader to the recent surveys \cite{Ward.12,Dai.12}.  When studying heavy-traffic limits of queueing systems, there are generally two different levels of granularity at which one can study the associated limits: the fluid scale (fluctuations of order $n$) and the diffusion scale (fluctuations around the fluid scaling, often of order $n^{\frac{1}{2}}$).  For multi-server queues with abandonments, questions regarding convergence at the fluid scale have been studied in several works, including \cite{Kang.12, Kang.15c, Liu.12, Liu.11b, Liu.11f, Whitt.06e, Zhang.13}, where several of these works allow for systems with time-varying parameters, and we refer the interested reader to \cite{Dai.12} for additional references.  We will instead focus on the second level of granularity, i.e. the relevant diffusion-scaled process, where we note that these works often implicitly demonstrate a fluid scaling as well to provide an appropriate centering for the relevant diffusion limit.  Furthermore, we will focus exclusively on models with abandonments in the HW scaling regime, but note in passing that there is a vast literature (beyond the scope of our work) on closely related models and topics pertaining to abandonments, including (but by no means limited to): non-asymptotic analysis and numerical approaches \cite{Adan.17,Massey.14b,Massey.13c,Kat.15,Sar.13}; other scaling regimes \cite{Whitt.04d,Jennings.12b,Reed.12b}; single-server models \cite{Ward.03f,Glynn.05,Reed.08c}; network models \cite{Reed.04,Huang.13b,Massey.17d}; models with more complex dependency structure \cite{Wu.17,Ev.18}; control and dynamic staffing \cite{Zohar.08,Ari.17b,Liu.17c}; simulation \cite{Ni.15}; delay announcements \cite{Armony.09,Ibrahim.09}; strategic behavior \cite{Mandelbaum.00,Ata.17}; and ticket queues \cite{Pender.16d,Ding.15c}.
\subsection{Diffusion approximations.}
It seems the first work to address the question of diffusion-scaled multi-server queues with abandonments in the HW scaling was \cite{Flem.94}, which used the available closed-form expressions for the $M/M/n + M$ queue to work out the relevant asymptotics for the $M/M/n + M$ case, in the steady-state.  In \cite{Garn.02}, these results were also extended to the transient setting for the $M/M/n + M$ queue, and the results gained renewed popularity.  To be slightly more precise about the notion of a diffusion limit in the HW regime, and the relevant concerns, suppose we have fixed an inter-arrival distribution $A$ and service distribution $S$, both with mean 1, a patience distribution $R$ with continuous density function which is stritly positive and finite at 0, and a spare capacity parameter $B$ (which may be negative, positive, or zero).  Let $Q^n$ denote the queueing system with inter-arrival times i.i.d. distributed as $\frac{A}{n + B n^{\frac{1}{2}}}$, service times drawn i.i.d. from the unscaled distribution $S$, and patience times drawn i.i.d. from the unscaled distribution $R$.  For this system, let $\lbrace W^n(t), t \geq 0 \rbrace$ be the queue-length process (i.e. the number of jobs waiting in queue), where for now we leave the initial distribution unspecified.  It is well-known that under mild technical conditions such a system will have a unique steady-state distribution (see e.g. \cite{DDG14}), and we let $W^n(\infty)$ denote a random variable (r.v) with the corresponding steady-state distribution (on occasion we will also use r.v to denote a random vector).  In that case, the results of \cite{Flem.94} showed that for an appropriate sequence of $M/M/n + M$ queues, $\lbrace n^{-\frac{1}{2}} W^n(\infty), n \geq 1 \rbrace$ converges weakly to a limiting r.v $W^{\infty}(\infty)$, while the results of \cite{Garn.02} showed in addition that for each fixed $T \geq 0$, $\lbrace n^{-\frac{1}{2}} W^n(t)_{0 \leq t \leq T} \rbrace$ converges weakly to a limiting process $W^{\infty}(t)_{0 \leq t \leq T}$ (in the appropriate space and topology, i.e. a process-level convergence).  \cite{Garn.02} also showed that in this case a so-called interchange-of-limits holds, namely that $\lbrace W^{\infty}(t), t \geq 0 \rbrace$ converges, as $t \rightarrow \infty$, to $W^{\infty}(\infty)$.  Namely, 
the limit of the steady-states coincides with the steady-state of the limit, which is in general non-trivial to verify.  Related higher-order asymptotic expansions were derived in \cite{Zhang.12b}.  In \cite{Zelt.05b}, these results are extended to show a steady-state convergence in the $M/M/n + G$ setting, building on the non-asymptotic work of \cite{Bac.81}.  In \cite{Whitt.05g}, these results are extended to show a process-level convergence in the $G/M/n  + G$ setting.  Process-level weak convergence and strong approximation results were proven for time-varying networks of $M/M/n$ queues in \cite{Man.98}, with several related results summarized in \cite{Man.02b}, and results for related periodic Markovian systems proven in \cite{Puhalskii.08b}.  
\subsection{Diffusion approximations with non-Markovian service times.}
Note that all the aforementioned works essentially restrict to the setting of Markovian service times, where we note that in the asymptotic regimes considered generalizing the service distribution is typically far more challenging than generalizing the inter-arrival or patience distributions.  However, as empirically studied in \cite{Brown.05}, in many relevant applications of queueing models in service systems the assumption of Markovian service times does not hold, and may lead to poor policy choices.  Significant progress on this front was made in \cite{Dai.10}, in which a process-level weak-convergence result was proven for the case of phase-type (i.e. PH) service distributions (which are dense within the set of all distributions), building on a state-space collapse relating the number of abandoning jobs to a certain integral of the queue-length process \cite{Dai.10b}.  The associated limit process was a complex multi-dimensional Markovian diffusion with piece-wise linear drift and no explicit solution, called a piecewise Ornstein-Uhlenbeck (OU) process in later work, and could also be described as a certain continuous mapping of a simpler process.  We note that a related state-space collapse result was attained in \cite{Tal.09}, in which an asymptotic equivalence between the system with balking (i.e. in which jobs that will eventually abandon do so upon arrival) and with reneging (i.e. in which jobs remain in queue until they abandon and are counted towards the queue-length) is proven.  Process-level weak-convergence results very similar to those of \cite{Dai.10}, albeit for even more general service distributions, were attained in \cite{Man.12b}.  In the same paper, the question of proving analogous results for the associated sequence of steady-state distributions, and a corresponding interchange-of-limits, was explicitly stated as an open problem.  
\\\indent This was resolved in the sequence of papers \cite{Dieker.13,DDG14}.  In \cite{Dieker.13}, the authors proved that the piecewise OU processes which arose in \cite{Dai.10} had unique invariant measures, and were exponentially ergodic (i.e. converged exponentially quickly under an appropriate metric to this invariant measure from general initial conditions), where we note that such exponential ergodicity has since been extended to more general processes \cite{Ari.17b}.  In \cite{DDG14}, the results of \cite{Dieker.13} were combined with several additional arguments to prove a formal interchange-of-limits for $G/PH/n + M$ queues in the HW regime, namely that the sequence of appropriately normalized steady-state distributions is tight, and converges weakly to the unique invariant measure of the associated piecewise OU process identified in \cite{Dai.10}.
\subsection{Large deviations of the steady-state diffusion approximation.}
However, outside the case of Markovian service times, all of the aforementioned works left open any concrete understanding of the associated complex limit processes, beyond their exponential ergodicity \cite{Dieker.13}.  There have been several works which attempted to shed some light on these processes, especially their large deviations behavior, as the probability of rare events is typically a quantity of significant interest in such models, and we refer the interested reader to \cite{Sch.95} for an overview of such results, and to \cite{GG13} for an overview of such results for multi-server queues in the HW scaling.  In \cite{GS12}, the authors study multi-class Markovian queues with abandonments in the HW regime under general non-idling class-scheduling policies.  For such systems, they prove several results regarding existence of exponential moments and sub-Gaussian behaviors.  Their results include as a special case FCFS $M/H_k/n + M$ queues, i.e. queues in which service times are hyper-exponentially distributed, i.e. distributed as a finite mixture of exponential distributions.  For a sequence of such FCFS $M/H_k/n + M$ queues, the results of \cite{GS12} imply that, letting $W^{\infty}(\infty)$ denote the weak limit of $\lbrace n^{-\frac{1}{2}} W^n(\infty), n \geq 1 \rbrace$, $\limsup_{x \rightarrow \infty} x^{-2} \log\bigg( P\big( W^{\infty}(\infty) > x \big) \bigg) < 0$.  As noted in \cite{Dai.13b}, such a Gaussian decay is to be expected due to the mean-reverting nature of a system with abandonments, and based on explicit formulas associated with the Markovian case.  However, the associated bounds identified in \cite{GS12} involve both stochastic comparison and non-explicit Lyapunov-function-type arguments, and do not seem to shed much insight on the precise set of $z > 0$ such that (s.t.) $\E\bigg[ \exp\bigg( z \big(W^{\infty}(\infty)\big)^2 \bigg) \bigg] < \infty$, other than that this set is non-empty.  We note that related Lyapunov-type arguments and bounds also appear in \cite{Dieker.13}, but again it seems unlikely that the approach taken there could be used to identify the exact tail behavior.  We will loosely refer to $- \sup\bigg\lbrace z: \E\bigg[ \exp\bigg( z \big(W^{\infty}(\infty)\big)^2 \bigg) \bigg] < \infty \bigg\rbrace$ as the associated large deviations exponent, equivalently (should the limit exist) $\lim_{x \rightarrow \infty} x^{-2} \log\bigg( P\big( W^{\infty}(\infty) > x \big) \bigg)$, where we note that this concept is seemingly much more developed in the literature for queues without abandonments, in which the associated r.v typically exhibit exponential (as opposed to Gaussian) decay \cite{Sad.91}.  
\subsection{A conjecture of Dai and He \cite{Dai.13b}.}
To gain further insight into $W^{\infty}(\infty)$, and the quality of approximation that it provides to $W^n(\infty)$ for any fixed $n$, various numerical methods were developed in \cite{Dai.13b} for evaluating the probabilities of $W^{\infty}(\infty)$, building on the numerical methods for diffusions developed earlier by Dai and Harrison \cite{Dai.92a}.  As discussed in-depth in \cite{Dai.92a}, an important step in these numerical methods is selecting a reference density which is used to define an inner product, which is then used explicitly by the algorithm to determine the probabilities of $W^{\infty}(\infty)$ by defining the stationary measure through a so-called basic adjoint relationship with respect to this inner product, and then using certain appropriate approximations.  The authors in \cite{Dai.13b} devote a section to discussing the properties that make certain reference densities more desirable than others, and note that ideally the reference density would have similar tail behavior to the stationary measure itself, where in some sense the greater the similarity the faster their algorithm will converge.  However, the problem faced by the authors was that for the diffusion processes they were considering, i.e. the limit processes arising from $G/PH/n + G$ queues in the HW scaling, the tail behavior was unknown.  As such, to motivate their choice of reference density, the authors made a formal conjecture regarding the large deviations exponent for $W^{\infty}(\infty)$.  Namely, inspired by known results for the case of Markovian processing times, as well as an analogous result for the setting without abandonments \cite{GG13}, Dai and He made the following formal conjecture in Conjecture 2 of \cite{Dai.13b}.  Throughout we suppose without loss of generality (w.l.o.g) that the expected service time equals 1, and let $\sigma^2_A (\sigma^2_S) \in (0,\infty)$ denote the variance of $A (S)$, where we recall that in the $n$th system the inter-arrival distribution is $\frac{A}{n + B n^{\frac{1}{2}}}$ for some fixed real excess parameter $B$, the service distribution is $S$, and the patience distribution $R$ has a density which evalutes to $\theta \in (0,\infty)$ at $0$. We note that the precise conjecture of \cite{Dai.13b} makes some additional technical assumptions on $S$, and as we are only here interested in the case that $S$ is hyper-exponentially distributed and the inter-arrival and patience distributions are both exponentially distributed, which is easily verified to satisfy all these technical conditions, we simply refer the reader to \cite{Dai.13b} for the precise assumptions. 
\begin{conjecture}[Conjectured large deviations exponent \cite{Dai.13b}]\label{maincon1}
$$\lim_{x \rightarrow \infty} x^{-2} \log\bigg( P\big( W^{\infty}(\infty) > x \big) \bigg) = - \frac{\theta}{\sigma^2_A + \sigma^2_S}.$$
Customized to the $M/H_2/n + M$ setting in which service times are with probability (w.p.) $p$ exponentially distributed with rate $\mu_1$ and with probability $1 - p$ exponentially distributed with rate $\mu_2$, the canonical example of non-exponential service times explicitly used in \cite{Dai.12,DDG14,Dai.13b} to illustrate the piecewise OU process, it is easily verified (since $\E[S] = 1$ implies $p = \frac{\mu_1 (\mu_2 - 1)}{\mu_2 - \mu_1}$) that this conjecture is equivalent to
$$\lim_{x \rightarrow \infty} x^{-2} \log\bigg( P\big( W^{\infty}(\infty) > x \big) \bigg) = - \frac{\mu_1 \mu_2 \theta}{2 (\mu_1 + \mu_2 - 1)}.$$
Furthermore, the authors of \cite{Dai.13b} note that if this conjecture is true, it provides an example of a so-called insensitivity result, as the exponent is a function only of $\theta$ and the first two moments of the inter-arrival and service distributions.
\end{conjecture}
The authors of \cite{Dai.13b} implemented their numerical algorithm with a reference density strongly based on Conjecture\ \ref{maincon1}.  The demonstrated good numerical performance of their algorithm seemed to provide support for the conjecture.  Of course, as the authors of \cite{Dai.13b} note, their algorithm will in general converge as long as the reference density satisfies certain basic technical conditions, i.e. it need not capture the correct large deviations behavior, i.e. they could have gotten good numerical results even if Conjecture\ \ref{maincon1} was false.  Since publication of \cite{Dai.13b}, \cite{Dieker.13}, and \cite{DDG14}, no further progress has been made towards resolving Conjecture\ \ref{maincon1}, and in general the associated limit processes remain poorly understood.  We do note the works \cite{Gur13,Gur14,BD15} which use Stein's method and related approaches to further our understanding of how good an approximation the associated limit process provides for the rescaled queue length for any fixed $n$ (i.e. bounds the error in using the diffusion to approximate the pre-limit system), and provides some bounds for certain moments of the limiting diffusions.  However, as discussed in \cite{G17a}, developing further qualitative insight into the properties of the limiting processes themselves could shed considerable light on the pre-limit systems, where we do note that some care must be taken when transferring statements proven about the limit process to the pre-limit system, and refer the interested reader to \cite{GS12} for further discussion along these lines.  Prior to this work, Conjecture\ \ref{maincon1} remained open.
\subsection{Stochastic comparison methods for queues with abandonments.}
As it will play a central role in our analysis, we note that for systems without abandonments, a novel stochastic comparison approach for multi-server queues was developed in \cite{GG13}, and has since been used to prove results in several settings related to the HW scaling regime without abandonments \cite{G16,Goldberg.17h,G17a}.  Although there is a vast literature on stochastic comparison results for queues without abandonments \cite{Stoyan.83,Chen.13b}, the results for systems with abandonments are much more limited, seeming to consist of only a few works \cite{Bhat.91,Jouini.07, Boxma.93, Dai.10b,GS12,Moyal.17}.  One of the reasons for this is that certain basic monotonicities that hold in great generality for queueing systems without abandonments no longer hold in the presence of abandonments.  One non-monotonicty of particular relevance to the approach of \cite{GG13} is that in systems with abandonments, adding extra arrivals may actually cause the queue length to decrease along any given sample path.  For example, consider the following simple queueing system.  Suppose the queue is initially empty, and that there is a single arrival at time 1, where that single job has service time 10 and patience 4.  Note that in this queue, the number of jobs in system at time 10 equals 1.  However, consider the modified queuing system in which we add an extra arrival at time 0, with service time 6 and patience 1.  In this modified system, the original arrival at time 1 will abandon the system at time 5, and the extra arrival will depart the system at time 6, leaving 0 jobs in system at time 10.  Thus in systems with abandonments, adding jobs can actually decrease the number in system, which cannot happen in systems without abandonments.  This and related subtleties makes it considerably more challenging to make comparative statements for queues with abandonments.
\subsection{Related diffusion processes.}
Although the piecewise OU processes analyzed in \cite{Dai.10,Dieker.13,DDG14} remain poorly understood, as noted in \cite{Dieker.13}, there is a well-studied family of related processes called multi-dimensional OU processes \cite{Sato.94,Masuda.04}.  Indeed, as noted in \cite{Dieker.13}, one of the main hurdles to understanding the piecewise OU processes arising in \cite{Dai.10,Dieker.13,DDG14} is that the infinitesimal drift is a non-linear function 
of the multi-dimensional state of the system, involving terms of the form $\max\big(0, \overline{x}_1 + \overline{x}_2 \big)$, with $\overline{x}$ the state of the diffusion.  If within the system of stochastic differential equations (SDE) defining these processes, one simply replaces the terms of the form $\max\big(0, \overline{x}_1 + \overline{x}_2 \big)$ with $\overline{x}_1 + \overline{x}_2$, the drift becomes a linear function of the state, and the unique strong solution of the resulting system is a so-called multi-dimensional OU process.  Indeed, these processes are defined as the unique strong solutions to systems of SDE in which the drift is a linear function of the state, and the driving random process is Levy.  These processes have essentially closed-form solutions and much is known about them, and we refer the reader to \cite{Sato.94,Masuda.04} for formal definitions and further details.  We also note that there is some work on certain modified OU processes also reminiscient of the piecewise OU processes arising in \cite{Dai.10,Dieker.13}, including reflected OU processes \cite{Ward.03b,Xing.09} and OU processes driven by OU processes \cite{Berc.14}.  There is also a literature on the large deviation properties of reflections and Skorokhod maps of general diffusion processes \cite{Dup.87,DaiMi.11,Kang.14}.  We also note that one tool that has proven useful in analyzing such reflected processes is a generalization of the traditional Skorokhod reflection map machinery to processes with state-dependent drift, and we refer the interested reader to \cite{Reed.04,Reed.13} for details, especially \cite{Reed.04} which proved that a reflected multi-dimensional OU process can be used to model a certain queueing network in heavy-traffic, in a setting very different from the QED regime.
\subsection{The one-dimensional (i.e. Markovian) case.}
As for the $M/M/n + M$ queue our approach becomes straightforward, we briefly comment in more depth on certain facts either known or straightforward to prove (from known results) in this setting, as they set the stage for our own results and approach.  It is well-known that for any fixed $n$, such a queueing system can be modeled as a birth-death process, and is time-reversible in steady-state.  Furthermore, in this setting it is easy to see that if one attempts to ``mimic" the stochastic comparison approach of \cite{GG13}, i.e. adding an extra arrival any time a server would have gone idle, one ends up with essentially the same 1-dimensional birth-death process, the only difference being that there is now a reflection at 0.  In this case the steady-state distributions of both systems can be computed in essentially closed form using well-known properties of birth-death processes (as in \cite{Flem.94,Garn.02}).  Furthermore, it is similarly well-known that here the steady-state probabilities for the number in queue in the modified system will be equivalent to those of the first system conditioned on having a positive queue.  It follows that: 1. the steady-state number waiting in queue in the modified system stochastically dominates that in the original system (i.e. the non-monotonocities associated with systems with abandonments can be overcome due to the memoryless property) and 2. under the HW scaling both systems will have the same steady-state large-deviations behavior, where this tail exponent is explicitly given in \cite{Flem.94,Garn.02} and equals $-\frac{\theta}{2},$ consistent with the conjecture of Dai and He.  
\\\indent We note that in this case, it is easily verified that the resulting limiting 1-dimensional piecewise OU process $Q^{\infty}$ corresponding to the process-level weak limit of the true queueing system in the HW regime has non-linear drift function $b(x) = B - \theta x - (\theta - 1)(x^+ - x)$, with $x$ the state of the process and $x^+ \stackrel{\Delta}{=} \max(0,x)$ for any $x \in {\mathcal R}$.  Illustrating our aforementioned notion of ``relaxing" the nonlinearity $(x^+ - x)$, i.e. replacing $x^+$ with $x$, we arrive at the defining equations for a one-dimensional OU process with drift function $\overline{b}(x) = B - \theta x$.  As it is well-known that this one-dimensional OU process has the same large deviations exponent $-\frac{\theta}{2}$, we find that in the Markovian (i.e. one-dimensional) case, the upper-bounding process analogous to that defined for systems without abandonments in \cite{GG13}, the true queueing system, and the ``relaxed" OU process all have the same large deviations behavior.  Of course, when $S$ is non-Markovian, these arguments fundamentally break down and essentially nothing is known regarding the questions of interest.
\subsection{Overview of contribution.}
In this paper, we resolve Conjecture\ \ref{maincon1} in the negative.  We do so by explicitly computing the large deviations exponent for the setting of $M/H_2/n + M$ queues, in the parameter regime s.t. the abandonment rate is strictly less than the minimum exponential service rate.  We find that this true exponent is strictly less than the conjectured exponent, and does not depend only on the first two moments of the service distribution (i.e. does not exhibit an insensitivity-type result).  This implies that the probability of seeing a large queue length decays faster than conjectured.  Our main approach is to extend the stochastic comparison framework of \cite{GG13} to the setting of abandonments, which requires several novel and non-trivial contributions, and to combine with results from the theory of Gaussian processes.  Furthermore, we believe that our approach sheds light on several novel ways to think about multi-server queues with abandonments in heavy traffic, e.g. by formally connecting to so-called multi-dimensional OU processes (a family of tractable multi-dimensional diffusions), and should extend to much more general settings.
\subsection{Outline of rest of paper.}
The rest of the paper proceeds as follows.  We present our main results, and discuss their context and significance, in Section\ \ref{MainSec}.  We present a stochastic comparison result for $M/H_2/n + M$ queues in Section\ \ref{UpperSec}, showing that a certain modified queueing system in which an extra arrival is added any time a server would otherwise go idle provides an upper bound.  In Section\ \ref{SkorSec}, we prove that the upper bound from Section\ \ref{UpperSec} can be expressed as the solution to a certain generalized Skorokhod problem.  In Section\ \ref{AsupSec}, we provide an asymptotic analysis of the solution to the aforementioned Skorokhod problem in the HW regime, and combine with several results from the theory of Gaussian processes to complete the proof of the upper-bound of our main result.  We complete the proof of our main result by demonstrating a matching lower bound in Section\ \ref{LowerSec}.  We summarize our main results, and provide concluding remarks and some interesting directions for future research, in Section\ \ref{ConcSec}.  We include a technical appendix containing the proofs of several results from throughout the paper in Section\ \ref{AppSec}.  For clarity of exposition, we put an asterisk next to each statement whose proof appears in the appendix.  For similar reasons the details of certain straightforward proofs are omitted, and we put a dash next to each statement whose proof is omitted.
\section{Main results.}\label{MainSec}
\subsection{Brief review of model, notations, and assumptions.}
Recall that for $z > 0$, $Expo(z)$ refers to a generic exponentially distributed r.v. with rate $z$ (mean $z^{-1}$).  Suppose that a spare capacity parameter $B \in {\mathcal R}$, exponential rates $\mu_1,\mu_2 \in (0,\infty)$, mixing probability $p \in (0,1)$, and abandonment rate $\theta \in (0,\infty)$ are fixed.  For $n > B^2$,  let $\lambda_n \stackrel{\Delta}{=} n +  B n^{\frac{1}{2}}$, and $Q^n$ denote the $M/H_2/n + M$ queue with inter-arrival times distributed as $Expo(\lambda_n)$, service times hyper-exponentially distributed being $Expo(\mu_1)$ with probability $p$ and $Expo(\mu_2)$ with probability $1 - p$, and patience times distributed as $Expo(\theta)$.  W.l.o.g. (by a simple scaling argument), and for clarity of exposition, suppose throughout that the expected service time equals 1, i.e. $\frac{p}{\mu_1} + \frac{1 - p}{\mu_2} = 1$.  With a slight abuse of notation, we will sometimes refer to sequences such as $\lbrace Q^n, n \geq 1 \rbrace$ even though the relevant system is only defined for $n > B^2$, and in all cases this should be read as a sequence starting at some sufficiently large $n$ to ensure all quantities are well-defined.  For general feasible initial conditions (i.e. having a finite number of jobs in system, and having positive queue only if all servers are busy), let $W^n(t)$ denote the number of jobs waiting in queue in $Q^n$ at time $t$.  Then it follows from the results of \cite{DDG14} that $\lbrace W^n(t), t \geq 0 \rbrace$ converges in distribution (as $t \rightarrow \infty$) to a steady-state r.v. $W^n(\infty)$, independent of initial conditions.  Let $\hat{W}^n(\infty) \stackrel{\Delta}{=} n^{-\frac{1}{2}} W^n(\infty)$.  Then it further follows from the results of \cite{DDG14} that $\lbrace \hat{W}^n(\infty), n \geq 1 \rbrace$ converges in distribution (as $n \rightarrow \infty$) to a limiting non-defective (i.e. almost surely finite) r.v. $W^{\infty}(\infty)$ with support on ${\mathcal R}^+$.  Then the conjecture of Dai and He from \cite{Dai.13b} is equivalent to: $\lim_{x \rightarrow \infty} x^{-2} \log\bigg( P\big( W^{\infty}(\infty) > x \big) \bigg) = - \frac{\mu_1 \mu_2 \theta}{2 (\mu_1 + \mu_2 - 1)}.$
\subsection{Main results.}
We now state our main results.  As mentioned in Section\ \ref{IntroSec}, our proofs require that we restrict the range of allowed parameters as dictated by the following assumption.
\begin{Assumption}[Hazard-rate ordering assumption.]\label{A1}
$$0 < \theta < \min(\mu_1, \mu_2).$$
\end{Assumption}
Intuitively, Assumption\ \ref{A1} implies that the infinitesimal rate at which a job leaves the system is always greater while that job is in service, as opposed to while waiting in queue, irregardless of how long that job has already been in service or waiting in queue.  As we will discuss in-depth later, this assumption is needed to formally prove that our upper-bounding system is indeed an upper bound, and we leave the question of what happens outside this parameter range as an interesting question for future research.  We also note that similar assumptions appeared in several arguments of \cite{GS12}, and that although our results still hold if 
$\theta = \min(\mu_1, \mu_2)$ (i.e. the second part of the inequality is not strict), for clarity of exposition we assume strict inequality throughout.
\\\indent Then our main result is as follows.
\begin{theorem}\label{main1}
Under Assumption\ \ref{A1}, 
$$\lim_{x \rightarrow \infty} x^{-2} \log \bigg( P\big( W^{\infty}(\infty) > x \big) \bigg)  = - \frac{\theta( \theta + \mu_1 \mu_2 ) }{2( \theta + \mu_1 + \mu_2 - 1)}.$$
Furthermore, for all $\theta,\mu_1,\mu_2$ s.t. $\E[S] = 1$ and $\mu_1 \neq \mu_2$, it is true that $- \frac{\theta( \theta + \mu_1 \mu_2 ) }{2( \theta + \mu_1 + \mu_2 - 1)}$ is strictly less than the exponent $- \frac{\mu_1 \mu_2 \theta}{2 (\mu_1 + \mu_2 - 1)}$ conjectured by Dai and He in \cite{Dai.13b}.  Furthermore, the true exponent $- \frac{\theta( \theta + \mu_1 \mu_2 ) }{2( \theta + \mu_1 + \mu_2 - 1)}$ is not a function of only $\theta$ and the first two moments of the service distribution, and as such the insensitivity phenomena conjectured in \cite{Dai.13b} does not hold.  Namely, there exist two hyper-exponential distributions with the same first two moments for which the true exponent differs (for the same value of $\theta$).
\end{theorem}
\subsection{Discussion}
We now present a brief discussion of some implications of our main results, as well as some auxiliary results and aspects of our proof which will be discussed formally in later sections.  Our results represent the first concrete insights into the qualitative behavior of the limiting processes which arise from queues with non-Markovian service times in the HW regime, beyond their ergodicity \cite{Dieker.13}, the fact that the tails are in general sub-Gaussian \cite{GS12}, and the rate-of-convergence and moment bounds provided in \cite{BD15}.  We view the main implications to be four-fold.  
\subsubsection{Implications for Dai and He's numerical methods.}
First, the results have direct implications for the selection of a reference density for the numerical algorithms of \cite{Dai.13b}, and indeed suggest that one could get even better performance by using reference densities based on the true exponent which we identify.  Our  results also strongly suggest that the reference density used in \cite{Dai.13b} should indeed satisfy the technical requirements laid forth in \cite{Dai.13b} that any valid reference density should satisfy, which require that certain integrals with respect to the true stationary measure are finite.  We note that \cite{Dai.13b} was not able to formally verify this, as the true scaling of the stationary measure was not known.  
\subsubsection{Qualitative insights from explicit exponent.}
Second, the explicit form of the large deviations exponent yields some concrete insights.  For example, one finds that as $\theta \downarrow 0$, the large deviations exponent scales to first order like 
$- \frac{\mu_1 \mu_2}{2(\mu_1 + \mu_2 - 1)} \theta$, i.e. converges to 0 linearly in $\theta$ at that particular rate.  Another interesting regime is the setting in which one of the service rates (say $\mu_2$) diverges to $\infty$, in which case the large deviations exponent converges to $-\frac{\theta \mu_1}{2}$, i.e. converges to a constant independent of $\mu_2$.  We note that such a limit is conceptually related to the behavior of the $M/H^*_2/n + M$ queue, in which service times are with some probability 0 and with some probability exponentially distributed, whose limiting behavior in the HW regime was stated as an open question in \cite{Whitt.05g}.
\subsubsection{Why no insensitivity?}
Third, the fact that Dai and He's conjecture is false illustrates a more subtle manner in which abandonments change the qualitative behavior of the limiting processes, beyond the obvious incorporation of some sort of mean-reverting behavior.  There are several ways to understand this difference.  One can reason as follows, somewhat heuristically.  Without abandonments large deviations are only exponentially unlikely (not having a Gaussian decay).  However, as the underlying limit processes are still continuous mappings of Gaussian processes, this exponential behavior is in some sense driven by the following phenomena: when considering the probability of a rare event (e.g. exceeding a large $x$), one must aggregate these underlying Gaussian processes over very large time intervals whose length depends on $x$ itself, allowing the accumulation of a variance which itself depends on $x$, which leads (after standard Gaussian calculations) to the exponential (as opposed to Gaussian) decay.  But the fact that one must aggregate these Gaussian processes over large intervals to induce this phenomena means that certain universality results from probability theory, e.g. the central limit theorem for renewal processes which plays an important role in the covariance structure of these underlying Gaussian processes, will have washed away many of the particulars of the underlying Gaussian processes, leading to the kind of insensitivity phenomena conjectured by Dai and He.  However, in the setting of abandonments large deviations have a Gaussian decay, hence are not driven by such long-range aggregations, and hence do not exhibit such insensitivity.  A second reason is as follows.  For $z \in \lbrace 1,2 \rbrace$, we will call a job assigned an $Expo(\mu_z)$ service time a type-z job.
When the large deviations of the queue length have a Gaussian decay, the fluctuations of the number of type-1 and type-2 jobs in service are of the same qualitative magnitude (Gaussian) as the large deviations in the queue length, and hence more subtle interactions between the number of jobs waiting in queue and the type-composition of jobs must be considered when analyzing how rare events occur.
\subsubsection{Connection to multi-dimensional OU processes.}
Fourth, and closely related to the third point, is that our proofs ultimately suggest a close connection between the rough behavior (at least up to the large deviations properties) of the piecewise OU process limits arising from queues with abandonments in the HW regime, and the much simpler multi-dimensional OU processes mentioned in Section\ \ref{IntroSec}.  Although we will defer a formal description of these processes and their central role in our proofs to later sections, we now sketch the connection at a high level.  As noted in Section\ \ref{IntroSec}, one of the main hurdles to understanding the piecewise OU processes arising in \cite{Dai.10,Dieker.13,DDG14} is that the infinitesimal drift is a non-linear function of the multi-dimensional state of the system, involving terms of the form $\max\big(0, \overline{x}_1 + \overline{x}_2 \big)$, with $\overline{x}$ the state of the diffusion.  If within the system of SDE defining these processes, one simply replaces the terms of the form $\max\big(0, \overline{x}_1 + \overline{x}_2 \big)$ with $\overline{x}_1 + \overline{x}_2$, the drift becomes a linear function of the state, and the unique strong solution of the resulting system is a so-called multi-dimensional OU process, a well-understood family of Markov processes with essentially closed-form solutions.  Intuitively, this multi-dimensional OU process corresponds to a relaxation of the true queueing system in which the queue length is allowed to go negative while maintaining the linear relationships (e.g. as regards abandonments) which typically hold only when the queue is positive.  Thus, for example, part of the dynamics are a mean-reversion back to zero at a rate containing $\theta$ while the queue is negative, a sort of ``negative abandonment".  We note that through a simple linear transformation, we are further able to focus on a conceptually even simpler process in which the first component represents the ``net infinitesimal drift" induced by the discrepancy between the number of type-1 and type-2 jobs in service, while the second component represents the total number of jobs in queue (which can be negative), with the dynamics again linear (i.e. leading to another multi-dimensional OU process) but with the time-derivative of the first component a function only of the first component (i.e. the coefficient in the appropriate linear combination dictating the contribution of the second component equals zero).  Interestingly, our analysis shows that the large deviations exponent of the true limit process agrees with that of this ``relaxed" multi-dimensional OU process.
\subsubsection{Overview of insights from proofs.}
To better understand the connection and our proof approach, we now give a very brief overview of how our proof relates the true limit process to the aforementioned multi-dimensional OU process.  Note that when there is a positive queue, the dynamics of the true limit process and the multi-dimensional OU process agree.  Thus both our upper and lower bounds relate the true process to this multi-dimensional OU process by conceptually forcing a restriction to states in which the queue is positive.  For our upper bound, this is accomplished by considering a sequence of modified queueing systems in which such a positivity is enforced by introducing extra arrivals whenever a server would otherwise go idle.  The dynamics of this bounding system, in the limit, are essentially the same as that of the aforementioned multi-dimensional OU process under a certain type of reflection captured through a relevant Skorokhod-type mapping of a more primitive driving Gaussian process.  We note that this asymptotic upper bound is actually very explicit, i.e. the supremum of a Gaussian process with completely explicit covariance function.  
\\\indent In our lower bound, this is accomplished by considering the true limit process initialized with its unknown stationary measure, picking a fixed $T$, and taking as a lower bound the probability that the stationary process exceeds the desired large value at $T$, AND happened to stay positive on $[0,T]$.  As the true piecewise OU limit has the same dynamics as a multi-dimensional OU process over any such interval (as it will have never gone negative), it thus suffices to lower bound the corresponding probability for a multi-dimensional OU process initialized with not its own stationary measure, but with that of the complex piecewise OU limit.  Although this measure is non-explicit (and in fact what we are trying to understand), it suffices to use non-explicit mixing-type arguments to reason that as this unknown measure is non-defective, by picking large enough $T$ and then considering the large-deviations behavior at time $T$ one essentially acheives (modulo some error terms which can be dealt with) as a lower bound the probability that a (nicely initialized) multi-dimensional OU process is positive over an interval $[0,T]$ and exceeds a large value $x$ at $T$.  One then considers the double limit of fixing larger values of $T$, and for each considering the large deviations behavior (i.e. $x \rightarrow \infty$).  We show that ultimately enforcing positivity on $[0,T]$ does not change the lower-bound on the large deviations behavior that one is able to show, which will coincide with that of the associated stationary multi-dimensional OU process, which matches our upper bound (i.e. the reflection considered in our upper bound does not change the large deviations behavior).  Intuitively, this follows from the fact that the probability of remaining positive on $[0,T]$ is a function only of $T$, while the large deviations involve a limit $x \rightarrow \infty$, combined with the fact that enforcing positivity on $[0,T]$ does not decrease the probability of seeing a large deviation at time $T$ (which we show using Slepian's inequality).  
\\\indent Indeed, we prove that both these upper-and-lower bounding processes/approaches lead to the same large deviations exponent, matching that of the aforementioned multi-dimensional OU process one gets by linearizing all relevant non-linearities which appear in the definition of the limit process.  Although we only formalize these connections for the case of $H_2$ services, and for space considerations leave several of these connections implicit in our arguments and proofs, we believe that the aforementioned relationships should be quite robust, holding in much greater generality, and that the study of such multi-dimensional OU relaxations may yield a new set of tools for analyzing multi-server queues with abandonments in the HW regime broadly.  
\section{Upper bound.}\label{UpperSec}
In this section we present an upper bound for a generic $M/H_2/n + M$ queue $Q$ (under Assumption\ \ref{A1}), by considering (at least conceptually) a queueing system equivalent to $Q$ but with an extra arrival added any time a server would otherwise have gone idle.  By never transitioning to that part of the state-space in which there are idle servers, the non-linearities which make the analysis of the true queueing system difficult fade away.  This methodology was pioneered in \cite{GG13} for systems without abandonments, and has been applied to various queueing models without abandonments in \cite{G16,G17a,Goldberg.17h}.  However, several fundamental challenges arise when attempting to implement this framework in systems with abandonments.  In particular, as noted in Section\ \ref{IntroSec} and illustrated with a simple example, in systems with abandonments adding extra arrivals may actually cause the queue length to decrease along any given sample path.  Although there is a small literature on stochatic comparison results for systems with abandonments, see e.g \cite{Bhat.91, Jouini.07, Boxma.93, GS12, Dai.10b}, none of these results seems strong enough to determine whether or not the modified system which we will consider yields an upper bound, not to mention pin down the precise large deviations behavior of the original queueing system.  
\subsection{CTMC representation for generic $M/H_2/n + M$ queue $Q$.}
Due to the extra complexity introduced by abandonments, namely that our stochastic comparison results will no longer hold sample-path-wise (and instead require a more sophisticated coupling of various r.v), it will be simplest to formally define our modified upper-bounding system, and formalize our stochastic comparison arguments, in the language of CTMC.  We begin by reminding the reader of the natural CTMC descriptor for a generic $M/H_2/n + M$ queue $Q$ (not necessarily in the HW regime), which is standard in the literature (see e.g. \cite{Adan.17,Dai.10}).  When referring to such a generic queue $Q$, we suppose the arrival rate is $\lambda$; the service times are hyper-exponentially distributed with mean 1 and parameters $\mu_1,\mu_2,p$; and patience times are $Expo(\theta)$.  In particular, the state of the system is a 3-dimensional vector of non-negative integers, with the first component $N_1$ tracking the number of type-1 jobs in service, the second component $N_2$ tracking the number of type-2 jobs in service, and the third component $W$ tracking the number of jobs waiting in queue.  Here we recall that a type-z job is a job assigned an $Expo(\mu_z)$ service time.  We note that under this description, a job's type is assigned only upon its entering service (if it ever actually enters service, as opposed to abandoning the system), which allows us to track only the total number of jobs waiting in queue, and not the types of those jobs.  We state the relevant dynamics in terms of an appropriate infinitesimal generator, where for a countable-state CTMC and any two states $x,y$, the infinitesimal generator $q$ evaluated at $(x,y)$, i.e. $q(x,y)$, corresponds to the infinitesimal rate at which the CTMC transitions to state $y$ while in state $x$, as is standard in the literature.  Here we denote the infinitesimal generator of $Q$ by $q$.   In that case, $q$ has the following dynamics.
\\\begin{definition}[Generator for $Q$]
\ \begin{itemize}
\item $q\big( (N_1, N_2, W) , (N_1 + 1, N_2, W) \big) = I(N_1 + N_2 \leq n -1) \lambda p$;
\item $q\big( (N_1, N_2, W) , (N_1, N_2 + 1, W) \big) = I(N_1 + N_2 \leq n -1) \lambda (1 - p)$;
\item $q\big( (N_1, N_2, W) , (N_1, N_2, W + 1) \big) = I(N_1 + N_2 = n) \lambda$;
\item $q\big( (N_1, N_2, W) , (N_1 - 1, N_2, W) \big) =  I(W = 0) \mu_1 N_1$;
\item $q\big( (N_1, N_2, W) , (N_1, N_2 - 1, W) \big) =  I(W = 0) \mu_2 N_2$;
\item $q\big( (N_1, N_2, W) , (N_1 - 1, N_2 + 1, W - 1) \big) = I(W \geq 1) \mu_1 N_1 (1 - p)$; 
\item $q\big( (N_1, N_2, W) , (N_1 + 1, N_2 - 1, W - 1) \big) = I(W \geq 1) \mu_2 N_2 p$;
\item $q\big( (N_1, N_2, W) , (N_1, N_2, W - 1) \big) = I(W \geq 1) \big( \mu_1 N_1 p + \mu_2 N_2 (1 - p) + \theta W\big)$.
\end{itemize}
\end{definition}
\ \\\indent Note that in the final transition, the contribution $\mu_1 N_1 p + \mu_2 N_2 (1 - p)$ corresponds to the occurence of a departure when the job at the head of the queue is of the same type as the departing job, while the contribution $\theta W$ corresponds to an abandonment, with both of these events resulting in $N_1$ and $N_2$ remaining unchanged and $W$ being reduced by 1.
\\\indent 
\subsection{CTMC representation for bounding system $\overline{Q}$.}
We next formalize our upper-bounding system as a different CTMC.  For ease of exposition, it will be simpler to introduce our upper-bounding system as its own CTMC, as opposed to initially constructing it on the same probability space as $Q$ by adding extra arrivals.  Thus we will never formalize that this bounding CTMC is a type of modified queueing system with abandonments, in which extra arrivals occur any time a server would otherwise have gone idle, as our formal proofs will not require it, and doing so requires a significant amount of additional notation and definitions.  However, we will at times informally refer to this modified system as if it were a queueing system, referring to e.g. jobs, departures, abandonments, etc., to help build intuition.
\\\indent Our upper-bounding system, denoted $\overline{Q}$, will similarly be a 3-dimensional CTMC on $Z^+ \times Z^+ \times Z^+$, where the components will have the same intuitive meaning as in $Q$.  In particular, we now define the CTMC $\overline{Q} = (\overline{N}_1, \overline{N}_2, \overline{W}).$  The dynamics of $\overline{Q}$ as a CTMC are given by the infinitesimal generator $\overline{q}$ as follows.
\\\begin{definition}[Generator for $\overline{Q}$]
\ \begin{itemize}
\item $\overline{q}\big( (\overline{N}_1, \overline{N}_2, \overline{W}) , (\overline{N}_1, \overline{N}_2, \overline{W} + 1) \big) = \lambda$;
\item $\overline{q}\big( (\overline{N}_1, \overline{N}_2, \overline{W}) , (\overline{N}_1 - 1, \overline{N}_2 + 1, \overline{W}) \big) = I(\overline{W} = 0) \mu_1 \overline{N}_1 (1 - p)$;
\item $\overline{q}\big( (\overline{N}_1, \overline{N}_2, \overline{W}) , (\overline{N}_1 + 1, \overline{N}_2 - 1, \overline{W}) \big) = I(\overline{W} = 0) \mu_2 \overline{N}_2 p$;
\item $\overline{q}\big( (\overline{N}_1, \overline{N}_2, \overline{W}) , (\overline{N}_1 - 1, \overline{N}_2 + 1, \overline{W} - 1) \big) = I(\overline{W} \geq 1) \mu_1 \overline{N}_1 (1 - p)$;
\item $\overline{q}\big( (\overline{N}_1, \overline{N}_2, \overline{W}) , (\overline{N}_1 + 1, \overline{N}_2 - 1, \overline{W} - 1) \big) = I(\overline{W} \geq 1) \mu_2 \overline{N}_2 p$;
\item $\overline{q}\big( (\overline{N}_1, \overline{N}_2, \overline{W}) , (\overline{N}_1, \overline{N}_2, \overline{W} - 1) \big) = I(\overline{W} \geq 1) \big( \mu_1 \overline{N}_1 p + \mu_2 \overline{N}_2 (1 - p) + \theta \overline{W} \big)$.
\end{itemize}
\end{definition}
\ \\\indent Note that in $\overline{Q}$, when there are 0 jobs waiting in queue and a departure occurs (say of type-1), the transition which occurs is equivalent to what would happen if at that moment the departing job leaves, but in addition an arrival occurs, being type-1 w.p. $p$ (in which case no transition occurs whatsoever as the state is unchanged), and type-2 w.p. $1-p$ (in which case $\overline{N}_1$ decreases by 1 and $\overline{N}_2$ increases by 1).  Otherwise the dynamics are the same as those of $Q$.  
\subsection{Proof that $\overline{W}$ yields an upper bound for $W$.}
We now prove that $\overline{W}$ stochastically dominates $W.$  For a general stochastic process $Z = (Z^1,\ldots,Z^k)$, let $Z(t) = \big( Z^1(t), \ldots, Z^k(t) \big)$ denote the state of the process at time $t$.
\begin{theorem}\label{mainub1}
For all $n \geq 1$, all $\theta,\mu_1,\mu_2$ satisfying Assumption\ \ref{A1}, and all non-negative integer triplets $(n_1,n_2,w)$ s.t. $n_1 + n_2 = n$ and $w = 0$, one may construct $Q$ and $\overline{Q}$ on a common probability space s.t. $Q(0) = \overline{Q}(0) = (n_1,n_2,w)$, and w.p.1 $\overline{W}(t) \geq W(t)$ for all $t \geq 0$.
\end{theorem}
We note that under Assumption\ \ref{A1} such a result is intuitive, as the addition of extra arrivals to the system causes certain jobs to wait longer in queue, and assumption\ \ref{A1} ensures that (in an appropriate sense) the longer a job waitings in queue the ``slower" it leaves the system.  We also note that very similar proofs would allow us to prove an analogous result allowing for more general initial conditions for $Q$ and $\overline{Q}$, e.g. having strictly positive initial queue-length, but as this is not needed for our main results, for clarity of exposition we restrict to the setting in which there are exactly $n$ jobs initially in system.  There are many ways to prove Theorem\ \ref{mainub1}.  We proceed by considering a certain expanded Markovian representation for both systems, to be able to keep track of inidividual jobs, and then formally construct both systems on a common probability space as projections of a larger CTMC.  Alternative proofs proceed by e.g. using the virtual waiting time representation for the time-in-system \cite{Glynn.05}.  The complexity with such alternative approaches is that certain ``folk-theorem" type results for $M/G/n + M$ queues involving various representations and conditional independences related to Poisson processes and the system dynamics are not immediate for the modified system with extra arrivals, and rederiving these facts for our modified system (which no doubt hold) seems an even more tedious and cumbersome task than formally constructing the appropriate coupling as a CTMC (which, admittedly, is somewhat cumbersome).  We note that such alternative approaches to proving the desired stochastic comparison may be important for extending the proof to general service and patience distributions, under an assumption analogous to Assumption\ \ref{A1} for the supremum of the hazard rate of the patience distribution and the infinum of the hazard rate of the service distribution.  
\subsubsection{Expanded Markovian representations $Q^0$ and $\overline{Q}^0$ for $Q$ and $\overline{Q}$.}
Before proceeding with the proof of Theorem\ \ref{mainub1}, it will be convenient to provide expanded Markovian representations for both $Q$ and $\overline{Q}$, which we will denote by $Q^0$ and $\overline{Q}^0$ respectively, and which will enable us to explicitly couple the arrivals to both systems, as well as the abandonments from one system to service completions in the other.  In these expanded representations, we will keep track of particular customers by identifying each customer with an index equal to either its arrival number (i.e. its position in the total arrival order, if it is a standard arrival), or a dummy index zero (if it is a so-called special/dummy/extra arrival to be defined below) as well as a type (either 1 or 2) specifying whether its exponential rate of service is $\mu_1$ or $\mu_2$.  Jobs with index $i \geq 1$ will be colloquially referred to as standard/regular jobs/arrivals/customers, while jobs with index 0 will be colloquially referred to as special/dummy/extra jobs/arrivals/customers.  In summary, a given customer will be specified by a 2-dimensional vector (i.e. tuple), the first component being the index, the second component being the type.
\\\indent We will present our expanded Markovian framework at a level of generality to encompass both $Q$ and $\overline{Q}$, but note that certain features (e.g. the presence of dummy jobs in service) will only be used meaningfully in our expanded representation of $\overline{Q}$.  Our expanded Markovian representation will consist of three components, i.e. the state descriptor will have three components.
\\\indent First, there will be an unordered multi-set of 2-dimensional vectors, ${\mathcal S}$, 
which will correspond to the set of jobs in service.  Each of these vectors will have first component a non-negative integer, and second component belonging to $\lbrace 1,2 \rbrace$.  For any strictly positive integer $i$, ${\mathcal S}$ will by construction contain at most one 2-dimensional vector with first component $i$, while ${\mathcal S}$ may contain multiple copies of both $(0,1)$ and $(0,2)$.  If $(0,1)$ and/or $(0,2)$ have multiplicity greater than one, then each copy will correspond to a different job in service.  Let $|{\mathcal S}|$ denote the cardinality of this unordered multi-set, where we note that $|{\mathcal S}|$ is not necessarily the number of distinct elements, as the vectors $(0,1)$ and $(0,2)$ contribute their respective multiplicities to this cardinality.  Our construction will ensure that $|{\mathcal S}| \leq n$.  Given such an unordered multi-set ${\mathcal S}$ of tuples, and $i \geq 0$, $z \in \lbrace 1,2 \rbrace$, let $N({\mathcal S}, i, z)$ denote the multiplicity of $(i,z)$ in ${\mathcal S}$.  Given $i \geq 0$ and $z \in \lbrace 1,2 \rbrace$, let ${\mathcal S} \setminus (i,z)$ denote the unordered multi-set of tuples equivalent to ${\mathcal S}$, but with one copy of $(i,z)$ removed.  If $N\big( {\mathcal S}, i , z \big) = 0$, let ${\mathcal S} \setminus (i,z)$ simply denote ${\mathcal S}$.
Similarly, let ${\mathcal S} \bigcup (i,z)$ denote the unordered multi-set of tuples equivalent to ${\mathcal S}$, but with an additional copy of $(i,z)$ added.  We note that by construction, all operations we perform on such unordered multi-sets will never attempt to add a tuple $(i,z)$ if $i \geq 1$ and $\sum_{z=1}^2 N({\mathcal S},i,z)$ is already strictly positive.  In general we allow these operations to be composed, e.g. $\big( {\mathcal S} \bigcup (i_1,z_1) \big) \setminus (i_2,z_2)$.
\\\indent Second, there will be an ordered set (as opposed to multi-set) of 2-dimensional vectors, ${\mathcal W}$, which will correspond to the ordered set of jobs waiting in queue, with each vector corresponding to a particular job in queue.  Each of these vectors will have first component a strictly positive integer (i.e. by construction dummy customers with index 0 never wait in queue, as they only arrive to the system when a server would otherwise have gone idle), and second component belonging to $\lbrace 1,2 \rbrace$.  Let $|{\mathcal W}|$ denote the cardinality of this ordered set.  The job in position 1 will correspond to the job at the front of the queue (i.e. the next job to enter service).  The job in position $|{\mathcal W}|$ will be the job at the back of the queue.  Sometimes we will colloquially refer to the jobs as being in left-to-right order, with the job in position 1 being the left-most job, and the job in position $|{\mathcal W}|$ being the right-most job.  Our construction will ensure that $|{\mathcal S}| \leq n - 1$ implies ${\mathcal W} = \emptyset$, and ${\mathcal W} \neq \emptyset$ implies $|{\mathcal S}| = n$.  Given $i \geq 0$ and $z \in \lbrace 1,2 \rbrace$, let $N({\mathcal W}, i, z)$ denote the multiplicity of $(i,z)$ in ${\mathcal W}$, where we note that $\sum_{z=1}^2 N({\mathcal W},i,z) \in \lbrace 0,1 \rbrace$ for all $i \geq 0$, i.e. here the multiplicity reduces to the appropriate indicator.  Given such a ${\mathcal W}$, and $i \geq 1$, $z \in \lbrace 1,2 \rbrace$, let ${\mathcal W} \setminus (i,z)$ denote the ordered set of tuples equivalent to ${\mathcal W}$, but with element $(i,z)$ removed, and the relative order of all remaining tuples the same as in ${\mathcal W}$, if $N\big( {\mathcal W}, i , z \big) = 1$.  Otherwise, i.e. if $N\big( {\mathcal W}, i , z \big) = 0$, let ${\mathcal W} \setminus (i,z)$ simply denote ${\mathcal W}$.  Also, given $i \geq 1$ and $z \in \lbrace 1,2 \rbrace$, let ${\mathcal W} \bigcup (i,z)$ denote the ordered set of tuples equivalent to ${\mathcal W}$, but with the tuple $(i,z)$ appended to the right of all tuples in ${\mathcal W}$, with the relative order of the tuples already in ${\mathcal W}$ unaltered.  We note that by construction, all operations we perform on such ordered sets will never attempt to add a tuple $(i,z)$ if $i \geq 1$ and $\sum_{z=1}^2 N({\mathcal W},i,z)$ is already strictly positive, nor attempt to add a tuple with index 0.  Also, for such a set ${\mathcal W}$ and $j \in [1, |{\mathcal W}|]$, let ${\mathcal W}_j$ denote the tuple in position $j$ in ${\mathcal W}$, and let ${\mathcal W}_{j,i}$ denote the $i$th component of that tuple, $i \in \lbrace 1,2 \rbrace$.  
Given an event $A$, let $I(A)$ denote the corresponding indicator function.  Also, given a 2-dimensional vector $z$, let $I(z \in {\mathcal W})$ evaluate to 1 if $z$ appears in the ordered set ${\mathcal W}$ and 0 otherwise.
\\\indent Third, there will be a non-negative integer counter $\nu$ which keeps track of how many standard (i.e. non-dummy) arrivals there have been to the system, which will be necessary to properly assign indices to customers.
\\\indent Thus the state of a queue is such a triplet $\big( {\mathcal S}, {\mathcal W}, \nu \big)$.
Given such a triplet ${\mathcal Q}$, $i \geq 0$, and $z \in \lbrace 1,2 \rbrace$, let $N({\mathcal Q},i,z) \stackrel{\Delta}{=} N({\mathcal S},i,z) + N({\mathcal W},i,z)$.
\\\indent We now use the above definitions and operations to define CTMC $Q^{0}$ (equivalent to $Q$ in our expanded representation) and $\overline{Q}^{0}$ (equivalent to $\overline{Q}$ in our expanded representation), which will provide an upper bound for $Q^0$ (in an appropriate sense).  We first define the CTMC $Q^{0} = \big( {\mathcal S}^0, {\mathcal W}^0, \nu^0 \big)$, denoting the corresponding generator by $q^{0}$.  As a notational convenience, we define $p_1 \stackrel{\Delta}{=} p, p_2 \stackrel{\Delta}{=} 1 - p$.  
\\\begin{definition}[Generator for $Q^0$]
\ \begin{enumerate}[(i).]
\item \label{g1} If $|{\mathcal S}| \leq n -1$, for $z \in \lbrace 1,2 \rbrace$: 
$$q^{0}\bigg( \big( {\mathcal S}, {\mathcal W}, \nu \big) , \big( {\mathcal S} \bigcup( \nu + 1, z), {\mathcal W}, \nu + 1 \big) \bigg) = \lambda p_z.$$
\item \label{g2} If $|{\mathcal S}| = n$, for $z \in \lbrace 1,2 \rbrace$:
$$q^{0}\bigg( \big( {\mathcal S}, {\mathcal W}, \nu \big) , \big( {\mathcal S}, {\mathcal W} \bigcup( \nu + 1, z), \nu + 1 \big) \bigg) = \lambda p_z.$$
\item If ${\mathcal W} = \emptyset$:
\begin{enumerate}
\item \label{g3} For all $i \geq 1$ and $z \in \lbrace 1,2 \rbrace$ s.t. $N\big( {\mathcal S}, i, z \big) = 1$:
$$q^{0}\bigg( \big( {\mathcal S}, {\mathcal W}, \nu \big) , \big( {\mathcal S} \setminus(i,z), {\mathcal W}, \nu \big) \bigg) = \mu_z.$$
\end{enumerate}
\item If ${\mathcal W} \neq \emptyset$:
\begin{enumerate}
\item \label{g4} For all $i \geq 1$ and $z \in \lbrace 1,2 \rbrace$ s.t. $N\big( {\mathcal S}, i, z \big) = 1$:
$$q^{0}\Bigg( \bigg( {\mathcal S}, {\mathcal W}, \nu \bigg) , \bigg( \big( {\mathcal S} \setminus (i,z) \big) \bigcup {\mathcal W}_1, {\mathcal W} \setminus {\mathcal W}_1, \nu \big) \bigg) \Bigg) = \mu_z.$$
\item \label{g5} For all $i \geq 1$ and $z \in \lbrace 1,2 \rbrace$ s.t. $N\big( {\mathcal W}, i, z \big) = 1$:
$$q^{0}\bigg( \big( {\mathcal S}, {\mathcal W}, \nu \big) , \big( {\mathcal S}, {\mathcal W} \setminus (i,z), \nu \big) \bigg) = \theta.$$
\end{enumerate}
\end{enumerate}
\end{definition}
\ \\We next define the generator of the CTMC $\overline{Q}^{0} = \big( \overline{\mathcal S}^{0}, \overline{\mathcal W}^{0}, \overline{\nu}^{0}\big)$, denoting the corresponding generator by $\overline{q}^{0}.$ 
\\\begin{definition}[Generator for $\overline{Q}^0$]
\ \begin{enumerate}[(i).]
\item \label{h1} For all $z \in \lbrace 1,2 \rbrace$:
$$\overline{q}^{0}\bigg( \big( {\mathcal S}, {\mathcal W}, \nu \big) , \big( {\mathcal S}, {\mathcal W} \bigcup( \nu + 1, z), \nu + 1 \big) \bigg) = \lambda p_z.$$
\item If ${\mathcal W} = \emptyset$:
\begin{enumerate}
\item \label{h2} For all $i \geq 1$ and $z_1,z_2 \in \lbrace 1,2 \rbrace$ s.t. $N\big( {\mathcal S}, i, z_1 \big) = 1$:
$$\overline{q}^{0}\Bigg( \bigg( {\mathcal S}, {\mathcal W}, \nu \bigg) , \bigg( \big( {\mathcal S} \setminus(i,z_1) \big) \bigcup (0,z_2), {\mathcal W}, \nu \bigg) \Bigg) = \mu_{z_1} p_{z_2}.$$
\item \label{h3} For all $z_1,z_2 \in \lbrace 1,2 \rbrace$ s.t. $z_1 \neq z_2$:
$$\overline{q}^{0}\Bigg( \bigg( {\mathcal S}, {\mathcal W}, \nu \bigg) , \bigg( \big( {\mathcal S} \setminus(0,z_1) \big) \bigcup (0,z_2), {\mathcal W}, \nu \bigg) \Bigg) = N\big( {\mathcal S}, 0, z_1 \big) \mu_{z_1} p_{z_2}.$$
\end{enumerate}
\item If ${\mathcal W} \neq \emptyset$:
\begin{enumerate}
\item \label{h4} For all $i \geq 1$ and $z \in \lbrace 1,2 \rbrace$ s.t. $N\big( {\mathcal S}, i, z \big) = 1$:
$$\overline{q}^{0}\Bigg( \bigg( {\mathcal S}, {\mathcal W}, \nu \bigg) , \bigg( \big( {\mathcal S} \setminus (i,z) \big) \bigcup {\mathcal W}_1, {\mathcal W} \setminus {\mathcal W}_1, \nu \bigg) \Bigg) = \mu_z.$$
\item \label{h5} For all $i \geq 1$ and $z \in \lbrace 1,2 \rbrace$ s.t. $N\big( {\mathcal W}, i, z \big) = 1$:
$$\overline{q}^{0}\Bigg( \bigg( {\mathcal S}, {\mathcal W}, \nu \bigg) , \bigg( {\mathcal S}, {\mathcal W} \setminus (i,z), \nu \bigg) \Bigg) = \theta.$$
\item \label{h6} For all $z \in \lbrace 1,2 \rbrace$: 
$$\overline{q}^{0}\Bigg( \bigg( {\mathcal S}, {\mathcal W}, \nu \bigg) , \bigg( \big( {\mathcal S} \setminus (0,z) \big) \bigcup {\mathcal W}_1, {\mathcal W} \setminus {\mathcal W}_1, \nu \bigg) \Bigg) = N\big( {\mathcal S}, 0, z \big) \mu_z.$$
\end{enumerate}
\end{enumerate}
\end{definition}
\ \\We note that such non-standard Markovian representations for Markovian queues with abandonments are common in the literature, see e.g. \cite{BD15}.  We now define several notions which will act as a proxy for a state being ``reasonable".
\\\begin{definition}[Good state, 0-good state, $\overline{0}$-good state, very-good state]
\ \\In our expanded Markovian representation, a triplet ${\mathcal Q} = ({\mathcal S}, {\mathcal W}, \nu)$ is good if:
\begin{itemize}
\item $I(|{\mathcal S}| \leq n - 1, |{\mathcal W}| \geq 1) = 0$; 
\item For $i \geq 1$, $\sum_{z=1}^2 N({\mathcal Q},i,z) \leq 1$;
\item $N({\mathcal W},0) = 0$;
\item For $1 \leq i < j \leq  |{\mathcal W}|$, ${\mathcal W}_{i,1} < {\mathcal W}_{j,1}$; 
\item If $|{\mathcal W}| \geq 1$, then for all $(i,z) \in {\mathcal S}$, it holds that ${\mathcal W}_{1,1} > i$;
\item For all $i$ s.t. $\sum_{z=1}^2 N({\mathcal Q},i,z) = 1$, it holds that $i \leq \nu$.
\end{itemize}
Let us say that the triplet is:
\begin{itemize} 
\item 0-good if: in addition to being good, $\sum_{z=1}^2 N({\mathcal S},0,z) = 0$;
\item $\overline{0}$-good if: in addition to being good, $|{\mathcal S}| = n$;
\item very-good if: it is both 0-good and $\overline{0}$-good, and in addition $\nu = n$, ${\mathcal W} = \emptyset$, and $\sum_{z=1}^2 N({\mathcal S},i,z) = 1$ for all $i \in [1,n]$ (i.e. initially there is one job in service with index $i$ for each $i \in [1,n]$, there are no jobs waiting in queue, and the initial counter equals $n$).  
\end{itemize}
\end{definition}
\ \\\indent Note that the fourth and fifth requirements of being a good state ensure that the FCFS ordering is preserved, where we note that even in the presence of abandonments the relative order of those jobs in the system must obey the FCFS order.
\\\indent We now formally state the precise connection betweeen $Q$ and $Q^0$, and $\overline{Q}$ and $\overline{Q}^0$.  

\begin{lemma}--\label{qtoq0}
For any very-good triplet $({\mathcal S}, {\mathcal W}, \nu)$, one may construct $Q$ and $Q^0$ on a common probability space s.t.: $Q^0(0) = ({\mathcal S}, {\mathcal W}, \nu), Q(0) = \bigg( \sum_{i = 1}^{n} N({\mathcal S},i,1) , \sum_{i = 1}^{n} N({\mathcal S},i,2) , 0 \bigg)$, and w.p.1, for all $t \geq 0$, $N_1(t) = \sum_{i = 1}^{\infty} N\big( {\mathcal S}^0(t), i, 1 \big)$, $N_2(t) = \sum_{i = 1}^{\infty} N\big( {\mathcal S}^0(t), i, 2 \big),$ and $W(t) = \big|{\mathcal W}^0(t)\big|$.  
\\\indent Similarly, for any very-good triplet $({\mathcal S}, {\mathcal W}, \nu)$, one may construct $\overline{Q}$ and $\overline{Q}^0$ on a common probability space s.t.: $\overline{Q}^0(0) = ({\mathcal S}, {\mathcal W}, \nu), \overline{Q}(0) = \bigg( \sum_{i = 1}^{n} N({\mathcal S},i,1) , \sum_{i = 1}^{n} N({\mathcal S},i,2) , 0 \bigg)$, and w.p.1, for all $t \geq 0$, $\overline{N}_1(t) = \sum_{i = 0}^{\infty} N\big( \overline{\mathcal S}^0(t), i, 1 \big)$, $\overline{N}_2(t) = \sum_{i = 0}^{\infty} N\big( \overline{\mathcal S}^0(t), i, 2 \big),$ and $\overline{W}(t) = \big|\overline{\mathcal W}^0(t)\big|$.  
\end{lemma}

The proof of Lemma\ \ref{qtoq0} is elementary and follows from a standard and well-known style of coupling argument, using the fact that in both $Q^0$ and $\overline{Q}^0$ the job indices do not impact the abandonment or service dynamics, and the fact that there is similarly no impact on the abandonment and service dynamics if each arriving job generates its type upon arrival, as opposed to generating its type upon its beginning service (if it ever begins service).  We omit the details, and refer the interested reader to e.g. \cite{Tezcan.06} for an example of related arguments.
\\\indent Next, it will be useful to note that w.p.1, both $Q^{0}$ and $\overline{Q}^{0}$ are always in a ``reasonable" state, if they are initialized in a ``reasonable" state, as dictated by our previous definitions of 0-good and $\overline{0}$-good states. 
In all cases the stated results follow from a straightforward induction, and we omit the details.  

\begin{lemma}--\label{staysgood0}
If $Q^{0}(0)$ is 0-good, then w.p.1 $Q^{0}(t)$ is 0-good for all $t \geq 0$.  Similarly, if $\overline{Q}^{0}(0)$ is $\overline{0}$-good, then w.p.1 $\overline{Q}^{0}(t)$ is $\overline{0}$-good for all $t \geq 0$.
\end{lemma}

We note that several of our later arguments will use these properties implicitly, especially e.g. the fact that the indices must obey the FCFS ordering and that there is at most one job with any given strictly positive index in the system at any given time, and for clarity of exposition sometimes we will use these properties without explicitly referring back to Lemma\ \ref{staysgood0}.

\subsubsection{Formal statement of dominance and description of coupling.}
In light of Lemma\ \ref{qtoq0}, to prove the main upper bound result Theorem\ \ref{mainub1}, it suffices to prove an analogous statement for $Q^0$ and $\overline{Q}^0$, which we now formalize.  
\\\begin{definition}[$\geq$ comparator for expanded Markovian representations]
Given a 0-good triplet ${\mathcal Q}^1 = ({\mathcal S}^1, {\mathcal W}^1, \nu^1)$ and $\overline{0}$-good triplet ${\mathcal Q}^2 = ({\mathcal S}^2, {\mathcal W}^2, \nu^2)$, let us say that ${\mathcal Q}^2 \geq {\mathcal Q}^1$ if: 
\begin{itemize}
\item For all $i \geq 1$ and $z \in \lbrace 1,2 \rbrace$, $N({\mathcal W}^2,i,z) \geq N({\mathcal W}^1,i,z)$;
\item For all $i \geq 1$ and $z \in \lbrace 1,2 \rbrace$, $N({\mathcal Q}^2,i,z) \geq N({\mathcal S}^1,i,z)$;
\item $\nu^1 = \nu^2$.
\end{itemize}
\end{definition}
\ \\\indent Intuitively, ${\mathcal Q}^2 \geq {\mathcal Q}^1$ iff every job waiting in queue in ${\mathcal Q}^1$ is also waiting in queue in ${\mathcal Q}^2$, every job in system in ${\mathcal Q}^1$ is also in system in ${\mathcal Q}^2$ (note that a job in service in ${\mathcal Q}^1$ may be in service or waiting in queue in ${\mathcal Q}^2$), and the arrival counters are equal.  Then we will prove the following dominance result for $Q^0$ and $\overline{Q}^0$.

\begin{theorem}\label{ubmain1}
For any very-good triplet $({\mathcal S}, {\mathcal W}, \nu)$, one may construct $Q^0$ and $\overline{Q}^0$ on a common probability space s.t. $Q^0(0) = \overline{Q}^0(0) = ({\mathcal S}, {\mathcal W}, \nu)$, and w.p.1, for all $t \geq 0$, 
$\overline{Q}^0(t) \geq Q^0(t)$.
\end{theorem}

Combining with Lemmas\ \ref{qtoq0} and \ref{staysgood0}, Theorem\ \ref{ubmain1} is easily seen to imply the desired result Theorem\ \ref{mainub1}, so we now focus on proving Theorem\ \ref{ubmain1}.  We proceed by explicitly exhibiting such a construction/coupling as a 6-dimensional CTMC $Q'$, with the first 3 components corresponding to $Q^0$, and the second 3 components corresponding to $\overline{Q}^0$.  In addition to coupling the arrival times and types of regular jobs, the key idea of the coupling is as follows.  If a given job is in service in both $Q^0$ and $\overline{Q}^0$, its departure is coupled between the two systems.  However, if (e.g. due to the extra dummy jobs) a job is in service in $Q^0$ but still waiting in queue in $\overline{Q}^0$, we couple its departure from $Q^0$ to its abandonment from $\overline{Q}^0$ by splitting the infinitesimal rate $\mu_z$ at which it departs from $Q^0$ into two parts: an infinitesimal rate of $\theta$ (which is synced to the abandonment process of the job from $\overline{Q}^0$) and an infinitesimal rate $\mu_z - \theta$ (which is independent from $\overline{Q}^0$).  In this way, we are able to ensure that the job cannot depart from $\overline{Q}^0$ before it departs from $Q^0$, preserving the desired dominance.  Of course, here Assumption\ \ref{A1} is needed for this coupling to make sense.
\\\indent We first define the relevant CTMC $Q' = \big( {\mathcal S}^{0'}, {\mathcal W}^{0'}, \nu^{0'}, \overline{S}^{0'}, \overline{\mathcal W}^{0'}, \overline{\nu}^{0'} \big)$, denoting the corresponding generator by $q'$.  When describing a generator q, in general we have written out q(x,y) for every possible current state x and next state y.  However, note that this is in some sense redundant, as the form of x is generic, so the same state x is written repeatedly.  As here the state x is 6-dimensional and somewhat cumbersome, we will instead suppose that the state x is fixed at 
$$\bigg( {\mathcal S}^{0'} , {\mathcal W}^{0'}, \nu^{0'} , \overline{\mathcal S}^{0'}, \overline{\mathcal W}^{0'} , \overline{\nu}^{0'} \bigg),$$
and when we need reference the entire state simply use the shorthand ``$\cdot$".  In that case, the generator $q'$ is as follows.  
\\\begin{definition}[Generator for $Q'$]
\ \begin{enumerate}
\item \label{z1} If $|{\mathcal S}^{0'}| \leq n - 1$, for $z \in \lbrace 1,2 \rbrace$: 
$$q'\Bigg( \cdot, \bigg( {\mathcal S}^{0'} \bigcup( \nu^{0'} + 1, z), {\mathcal W}^{0'}, \nu^{0'} + 1, \overline{\mathcal S}^{0'}, \overline{\mathcal W}^{0'} \bigcup( \overline{\nu}^{0'} + 1, z), \overline{\nu}^{0'} + 1 \bigg) \Bigg) = \lambda p_z.$$
\item \label{z2} If $|{\mathcal S}^{0'}| = n$ , for $z \in \lbrace 1,2 \rbrace$:
$$q'\Bigg( \cdot, \bigg( {\mathcal S}^{0'}, {\mathcal W}^{0'} \bigcup( \nu^{0'} + 1, z) , \nu^{0'} + 1, \overline{\mathcal S}^{0'}, \overline{\mathcal W}^{0'} \bigcup( \overline{\nu}^{0'} + 1, z), \overline{\nu}^{0'} + 1 \bigg) \Bigg) = \lambda p_z.$$
\item If ${\mathcal W}^{0'} = \emptyset$ and $\overline{\mathcal W}^{0'} = \emptyset$:
\begin{enumerate}
\item \label{z3} For all $i \geq 1$ and $z_1,z_2 \in \lbrace 1,2 \rbrace$ s.t. $N\big( {\mathcal S}^{0'}, i, z_1 \big) + N\big( {\mathcal W}^{0'}, i, z_1 \big)  = 0, N\big( \overline{\mathcal S}^{0'}, i, z_1 \big) = 1$: 
$$q'\Bigg( \cdot, \bigg( {\mathcal S}^{0'} , {\mathcal W}^{0'} , \nu^{0'}, \big( \overline{\mathcal S}^{0'} \setminus (i,z_1) \big) \bigcup (0,z_2), \overline{\mathcal W}^{0'} , \overline{\nu}^{0'} \bigg) \Bigg) = \mu_{z_1} p_{z_2}.$$
\item \label{z4} For all $i \geq 1$ and $z_1,z_2 \in \lbrace 1,2 \rbrace$ s.t. $N\big( {\mathcal S}^{0'}, i, z_1 \big) = 1 , N\big( \overline{\mathcal S}^{0'}, i, z_1 \big) = 1$:
$$q'\Bigg( \cdot, \bigg( {\mathcal S}^{0'} \setminus (i,z_1), {\mathcal W}^{0'}, \nu^{0'}, \big( \overline{\mathcal S}^{0'} \setminus (i,z_1) \big) \bigcup (0,z_2), \overline{\mathcal W}^{0'}, \overline{\nu}^{0'} \bigg) \Bigg) = \mu_{z_1} p_{z_2}.$$
\item \label{z5} For all $z_1,z_2 \in \lbrace 1,2 \rbrace$ s.t. $z_1 \neq z_2$: 
$$q'\Bigg( \cdot, \bigg( {\mathcal S}^{0'} , {\mathcal W}^{0'} , \nu^{0'} , \big( \overline{\mathcal S}^{0'} \setminus (0,z_1) \big) \bigcup (0,z_2) , \overline{\mathcal W}^{0'} , \overline{\nu}^{0'} \bigg) \Bigg) = N\big( \overline{\mathcal S}^{0'}, 0, z_1 \big) \mu_{z_1} p_{z_2}.$$
\end{enumerate}
\item If ${\mathcal W}^{0'} = \emptyset$ and $\overline{\mathcal W}^{0'} \neq \emptyset$:
\begin{enumerate}
\item \label{z6} For all $i \geq 1$ and $z \in \lbrace 1,2 \rbrace$ s.t. $N\big( {\mathcal S}^{0'}, i, z \big) + N\big( {\mathcal W}^{0'}, i, z \big) = 0 , N\big( \overline{\mathcal S}^{0'}, i, z \big) = 1$:
$$q'\Bigg( \cdot , \bigg( {\mathcal S}^{0'} , {\mathcal W}^{0'} , \nu^{0'}, \big( \overline{\mathcal S}^{0'} \setminus (i,z) \big) \bigcup \overline{\mathcal W}^{0'}_1, \overline{\mathcal W}^{0'} \setminus \overline{\mathcal W}^{0'}_1 , \overline{\nu}^{0'} \bigg) \Bigg)  = \mu_z.$$
\item \label{z7} For all $i \geq 1$ and $z \in \lbrace 1,2 \rbrace$ s.t. $N\big( {\mathcal S}^{0'}, i, z \big) + N\big( {\mathcal W}^{0'}, i, z \big) = 0 , N\big( \overline{\mathcal W}^{0'}, i, z \big) = 1$:
$$q'\Bigg( \cdot, \bigg( {\mathcal S}^{0'} , {\mathcal W}^{0'} , \nu^{0'}, \overline{\mathcal S}^{0'} , \overline{\mathcal W}^{0'} \setminus (i,z), \overline{\nu}^{0'} \bigg) \Bigg) = \theta.$$
\item \label{z8} For all $i \geq 1$ and $z \in \lbrace 1,2 \rbrace$ s.t. $N\big( {\mathcal S}^{0'}, i, z \big) = 1 , N\big( \overline{\mathcal S}^{0'}, i, z \big) = 1$:
$$q'\Bigg( \cdot, \bigg( {\mathcal S}^{0'} \setminus (i,z), {\mathcal W}^{0'}, \nu^{0'}, \big( \overline{\mathcal S}^{0'} \setminus (i,z) \big) \bigcup \overline{\mathcal W}^{0'}_1, \overline{\mathcal W}^{0'} \setminus \overline{\mathcal W}^{0'}_1, \overline{\nu}^{0'} \bigg) \Bigg) = \mu_z.$$
\item For all $i \geq 1$ and $z \in \lbrace 1,2 \rbrace$ s.t. $N\big( {\mathcal S}^{0'}, i, z \big) = 1 , N\big( \overline{\mathcal W}^{0'}, i, z \big) = 1$:
\begin{enumerate}
\item \label{z9} $$q'\Bigg( \cdot , \bigg( {\mathcal S}^{0'} \setminus (i,z) , {\mathcal W}^{0'} , \nu^{0'}, \overline{\mathcal S}^{0'} , \overline{\mathcal W}^{0'} \setminus (i,z), \overline{\nu}^{0'} \bigg) \Bigg) = \theta;$$
\item \label{z10} $$q'\Bigg( \cdot, \bigg( {\mathcal S}^{0'} \setminus (i,z) , {\mathcal W}^{0'} , \nu^{0'}, \overline{\mathcal S}^{0'} , \overline{\mathcal W}^{0'} , \overline{\nu}^{0'} \bigg) \Bigg) = \mu_z - \theta.$$
\end{enumerate}
\item \label{z12} For all $z \in \lbrace 1,2 \rbrace$: 
$$q'\Bigg( \cdot , \bigg( {\mathcal S}^{0'} , {\mathcal W}^{0'} , \nu^{0'} , \big( \overline{\mathcal S}^{0'} \setminus (0,z) \big) \bigcup \overline{\mathcal W}^{0'}_1 , \overline{\mathcal W}^{0'} \setminus \overline{\mathcal W}^{0'}_1 , \overline{\nu}^{0'} \bigg) \Bigg) = 
N\big( \overline{\mathcal S}^{0'}, 0, z \big) \mu_z.$$
\end{enumerate}
\item If ${\mathcal W}^{0'} \neq \emptyset$ and $\overline{\mathcal W}^{0'} \neq \emptyset$:
\begin{enumerate}
\item \label{z14} For all $i \geq 1$ and $z \in \lbrace 1,2 \rbrace$ s.t. $N\big( {\mathcal S}^{0'}, i, z \big) + N\big( {\mathcal W}^{0'}, i, z \big) = 0 , N\big( \overline{\mathcal S}^{0'}, i, z \big) = 1$:
$$q'\Bigg( \cdot , \bigg( {\mathcal S}^{0'} , {\mathcal W}^{0'} , \nu^{0'}, \big( \overline{\mathcal S}^{0'} \setminus (i,z) \big) \bigcup \overline{\mathcal W}^{0'}_1, \overline{\mathcal W}^{0'} \setminus \overline{\mathcal W}^{0'}_1 , \overline{\nu}^{0'} \bigg) \Bigg)  = \mu_z.$$
\item \label{z15} For all $i \geq 1$ and $z \in \lbrace 1,2 \rbrace$ s.t. $N\big( {\mathcal S}^{0'}, i, z \big) + N\big( {\mathcal W}^{0'}, i, z \big) = 0 , N\big( \overline{\mathcal W}^{0'}, i, z \big) = 1$:
$$
q'\Bigg( \cdot, \bigg( {\mathcal S}^{0'} , {\mathcal W}^{0'} , \nu^{0'}, \overline{\mathcal S}^{0'} , \overline{\mathcal W}^{0'} \setminus (i,z), \overline{\nu}^{0'} \bigg) \Bigg) = \theta.$$
\item \label{z16} For all $i \geq 1$ and $z \in \lbrace 1,2 \rbrace$  s.t. $N\big( {\mathcal S}^{0'}, i, z \big) = 1 , N\big( \overline{\mathcal S}^{0'}, i, z \big) = 1$:
$$q'\Bigg( \cdot, \bigg( \big( {\mathcal S}^{0'} \setminus (i,z) \big) \bigcup {\mathcal W}^{0'}_1, {\mathcal W}^{0'} \setminus {\mathcal W}^{0'}_1, \nu^{0'}, \big( \overline{\mathcal S}^{0'} \setminus (i,z) \big) \bigcup \overline{\mathcal W}^{0'}_1, \overline{\mathcal W}^{0'} \setminus \overline{\mathcal W}^{0'}_1, \overline{\nu}^{0'} \bigg) \Bigg) = \mu_z.$$
\item For all $i \geq 1$ and $z \in \lbrace 1,2 \rbrace$ s.t. $N\big( {\mathcal S}^{0'}, i, z \big) = 1 , N\big( \overline{\mathcal W}^{0'}, i, z \big) = 1$:
\begin{enumerate}
\item \label{z17} $$q'\Bigg( \cdot , \bigg( \big( {\mathcal S}^{0'} \setminus (i,z) \big) \bigcup {\mathcal W}^{0'}_1 , {\mathcal W}^{0'} \setminus {\mathcal W}^{0'}_1, \nu^{0'}, \overline{\mathcal S}^{0'} , \overline{\mathcal W}^{0'} \setminus (i,z), \overline{\nu}^{0'} \bigg) \Bigg) = \theta;$$
\item \label{z18} $$q'\Bigg( \cdot, \bigg( \big( {\mathcal S}^{0'} \setminus (i,z) \big) \bigcup {\mathcal W}^{0'}_1 , {\mathcal W}^{0'} \setminus {\mathcal W}^{0'}_1 , \nu^{0'}, \overline{\mathcal S}^{0'} , \overline{\mathcal W}^{0'} , \overline{\nu}^{0'} \bigg) \Bigg) = \mu_z - \theta.$$
\end{enumerate}
\item \label{z19} For all $i \geq 1$ and $z \in \lbrace 1,2 \rbrace$ s.t. $N\big( {\mathcal W}^{0'}, i, z \big) = 1 , N\big( \overline{\mathcal W}^{0'}, i, z \big) = 1$:
$$q'\Bigg( \cdot, \bigg( {\mathcal S}^{0'} , {\mathcal W}^{0'} \setminus (i,z) , \nu^{0'}, \overline{\mathcal S}^{0'} , \overline{\mathcal W}^{0'} \setminus (i,z), \overline{\nu}^{0'} \bigg) \Bigg) = \theta.$$
\item \label{z20} For all $z \in \lbrace 1,2 \rbrace$: 
$$q'\Bigg( \cdot , \bigg( {\mathcal S}^{0'} , {\mathcal W}^{0'} , \nu^{0'} , \big( \overline{\mathcal S}^{0'} \setminus (0,z) \big) \bigcup \overline{\mathcal W}^{0'}_1 , \overline{\mathcal W}^{0'} \setminus \overline{\mathcal W}^{0'}_1 , \overline{\nu}^{0'} \bigg) \Bigg) = 
N\big( \overline{\mathcal S}^{0'}, 0, z \big) \mu_z.$$
\end{enumerate}
\end{enumerate}
\end{definition}
In that case, an argument completely analogous to Lemma\ \ref{staysgood0} leads to the following, and we omit the details.
\begin{lemma}--\label{staysgood1}
If $\big( {\mathcal S}^{0'}(0), {\mathcal W}^{0'}(0), \nu^{0'}(0) \big)$ is 0-good and $\big(\overline{\mathcal S}^{0'}(0), \overline{\mathcal W}^{0'}(0), \overline{\nu}^{0'}(0) \big)$ is $\overline{0}$-good, then w.p.1 
$\big( {\mathcal S}^{0'}(t), {\mathcal W}^{0'}(t), \nu^{0'}(t) \big)$ is 0-good and $\big(\overline{\mathcal S}^{0'}(t), \overline{\mathcal W}^{0'}(t), \overline{\nu}^{0'}(t) \big)$ is $\overline{0}$-good for all $t \geq 0$.
\end{lemma}
\subsubsection{Proof of stochastic ordering.}
It is easily verified that to prove Theorem\ \ref{ubmain1}, it suffices to prove the following two lemmas.  
\begin{lemma}\label{staysgood3}
If: 1. $\big( {\mathcal S}^{0'}(0), {\mathcal W}^{0'}(0), \nu^{0'}(0) \big)$ is 0-good, 2. $\big(\overline{\mathcal S}^{0'}(0), \overline{\mathcal W}^{0'}(0), \overline{\nu}^{0'}(0) \big)$ is $\overline{0}$-good, and 3. $\big(\overline{\mathcal S}^{0'}(0), \overline{\mathcal W}^{0'}(0), \overline{\nu}^{0'}(0) \big) \geq \big( {\mathcal S}^{0'}(0), {\mathcal W}^{0'}(0), \nu^{0'}(0) \big)$, then 
$\big(\overline{\mathcal S}^{0'}(t), \overline{\mathcal W}^{0'}(t), \overline{\nu}^{0'}(t) \big) \geq \big( {\mathcal S}^{0'}(t), {\mathcal W}^{0'}(t), \nu^{0'}(t) \big)$ for all $t \geq 0$.
\end{lemma}

\begin{lemma}*\label{marginals1}
Under the same assumptions as Lemma\ \ref{staysgood3}, $\big\lbrace \big( {\mathcal S}^{0'}(t), {\mathcal W}^{0'}(t), {\nu}^{0'}(t) \big), t \geq 0 \big\rbrace$ has the same distribution (at the process level) as 
$\big\lbrace \big( {\mathcal S}^{0}(t), {\mathcal W}^{0}(t), {\nu}^{0}(t) \big), t \geq 0 \big\rbrace$; and $\big\lbrace \big( \overline{\mathcal S}^{0'}(t), \overline{\mathcal W}^{0'}(t), \overline{\nu}^{0'}(t) \big), t \geq 0 \big\rbrace$ has the same distribution (at the process level) as $\big\lbrace \big( \overline{\mathcal S}^{0}(t), \overline{\mathcal W}^{0}(t), \overline{\nu}^{0}(t) \big), t \geq 0 \big\rbrace$, in both cases supposing that corresponding processes have the same initial conditions.
\end{lemma}

\proof{Proof of Lemma\ \ref{staysgood3}:}
It suffices to prove that if $Q'$ is in a state ${\mathcal Q}' = \bigg( {\mathcal S}^{0'} , {\mathcal W}^{0'}, \nu^{0'} , \overline{\mathcal S}^{0'}, \overline{\mathcal W}^{0'} , \overline{\nu}^{0'} \bigg)$ s.t.: 1. $\big( {\mathcal S}^{0'}, {\mathcal W}^{0'}, \nu^{0'} \big)$ is 0-good, 2. $\big(\overline{\mathcal S}^{0'}, \overline{\mathcal W}^{0'}, \overline{\nu}^{0'} \big)$ is $\overline{0}$-good, and 3. $\big(\overline{\mathcal S}^{0'}, \overline{\mathcal W}^{0'}, \overline{\nu}^{0'} \big) \geq \big( {\mathcal S}^{0'}, {\mathcal W}^{0'}, \nu^{0'} \big)$, then w.p.1 it can only transition to another state with those properties.  We note that throughout, we will use the results of Lemma\ \ref{staysgood1} implicitly.  Our assumptions about ${\mathcal Q}'$ are easily seen to imply that
\begin{equation}\label{useit1}
|\overline{\mathcal S}^{0'}| = n \geq |{\mathcal S}^{0'}|\ \ \ ,\ \ \ \textrm{and}\ \ \ |\overline{\mathcal W}^{0'}| \geq |{\mathcal W}^{0'}|;
\end{equation}
and also that for all $i \geq 1$ and $z \in \lbrace 1,2 \rbrace$,
\begin{equation}\label{useit2}
N\big( \overline{\mathcal S}^{0'}, i, z\big) + N\big( \overline{\mathcal W}^{0'}, i, z\big) \geq N\big( {\mathcal S}^{0'}, i, z\big)\ \ \ ,\ \ \ 
N\big( \overline{\mathcal W}^{0'}, i, z\big) \geq N\big( {\mathcal W}^{0'}, i, z\big);
\end{equation}
which we will also use.  We proceed by a case analysis, analyzing each transition.  Let ${\mathcal Q}'' = \bigg( {\mathcal S}^{0''} , {\mathcal W}^{0'}, \nu^{0''} , \overline{\mathcal S}^{0''}, \overline{\mathcal W}^{0''} , \overline{\nu}^{0''} \bigg)$ denote the state to which $Q'$ next transitions.  First, note that the requirement involving the equivalence of $\nu''$ and $\overline{\nu}''$ is easily seen to be maintained under all transitions (i.e. both parameters either stay unchanged or increase by 1, together, under every transition), so we may focus on the other requirements.  
\\\indent Let us first consider transitions (\ref{z1}) and (\ref{z2}).  Under both transitions, our assumptions about ${\mathcal Q}'$ imply that the required inequalities in the definition of $\geq$ (ensuring ${\mathcal Q}^{''}$ has the same desired properties) hold for all $(i,z)$ s.t. $N({\mathcal S}^{0'},i,z) + N({\mathcal W}^{0'},i,z) = 1$ with $i \leq \nu^{0'}$.   By construction, for the one additional job with index $\nu^{0''}$ (first supposing this job has type 1), it holds that $N({\mathcal S}^{0''},\nu^{0''},1) + N({\mathcal W}^{0''},\nu^{0''},1) = 1 = N(\overline{\mathcal W}^{0''},\nu^{0''},1)$, which ensures that the required inequalities hold for all $(i,z)$ and hence that ${\mathcal Q}''$ has the desired properties.  The argument is identical if the new job has type 2.  
\\\indent Next, let us consider the transitions for which ${\mathcal W}^{0'} = \emptyset$ and $\overline{\mathcal W}^{0'} = \emptyset$.  We note that under all cases s.t. ${\mathcal W}^{0'} = \emptyset$, we need not worry about the requirement (on ${\mathcal Q}''$) requiring that $N(\overline{\mathcal W}^{0''},i,z) \geq N({\mathcal W}^{0''},i,z)$, as for the transitions being considered no arrivals occur and thus  ${\mathcal W}^{0''}$ remains equal to $\emptyset$ under the analyed transitions.  Under transition\ (\ref{z3}), the desired properties for ${\mathcal Q}''$ are inhered from ${\mathcal Q}'$, as $N\big( {\mathcal S}^{0'}, i, z_1 \big) + N\big( {\mathcal W}^{0'}, i, z_1 \big)  = 0$ anyways.  Under transition\ (\ref{z4}), the desired properties are inherited from ${\mathcal Q}'$, since the same $(i,z)$ pair is removed from both ${\mathcal S}^{0'}$ and $\overline{\mathcal S}^{0'}$.  Under transition (\ref{z5}), the desired properties are inherited from ${\mathcal Q}'$ as the only job removed from $\overline{\mathcal S}^{0'}$ has index 0.  
\\\indent Next, let us consider the transitions for which ${\mathcal W}^{0'} = \emptyset$ and $\overline{\mathcal W}^{0'} \neq \emptyset$.  Under transitions (\ref{z6}) and (\ref{z7}), the desired properties are inherited from ${\mathcal Q}'$ for the same reason as transition (\ref{z3}).  Under transitions (\ref{z8}) and (\ref{z9}), the desired properties are inherited from ${\mathcal Q}'$ for the same reason as transition (\ref{z4}).  Under transition (\ref{z10}), the desired properties are inherited from ${\mathcal Q}'$ as the job $(i,z)$ is only removed from ${\mathcal S}^{0'}$, while remaining in $\overline{\mathcal W}^{0'}$, preserving the desired property (as it allows for any given $(i,z)$ to appear only in the dominating system).  Under transition (\ref{z12}), the desired properties are inherited from ${\mathcal Q}'$ for the same reason as transition (\ref{z5}).  
\\\indent Finally, let us consider the transitions for which ${\mathcal W}^{0'} \neq \emptyset, \overline{\mathcal W}^{0'} \neq \emptyset$.  Here a bit more care will have to be taken, as we will have to verify that $N({\mathcal W}^{0''},i,z) \leq N(\overline{\mathcal W}^{0''},i,z)$, in addition to $N({\mathcal S}^{0''},i,z) \leq N(\overline{\mathcal S}^{0''},i,z) + N(\overline{\mathcal W}^{0''},i,z)$.  Under transition (\ref{z14}), the reasoning is the same as transition (\ref{z3}), with the caveat that we must verify that $N({\mathcal W}^{0''},i,z) \leq N(\overline{\mathcal W}^{0''},i,z)$.  To do so, we will prove that $\overline{\mathcal W}^{0'}_1 \notin {\mathcal W}^{0'},$ which (combined with our assumptions about ${\mathcal Q}'$) is easily seen to suffice.  Indeed, suppose for contradiction that $\overline{\mathcal W}^{0'}_1 \in {\mathcal W}^{0'}$.  Combined with the fact that ${\mathcal W}^{0'} \neq \emptyset$ and the definition of $\geq$, we conclude that there exist $(i_1,z_1),\ldots,(i_n,z_n) \in {\mathcal S}^{0'}$ s.t.: $i_k < \overline{W}^{0'}_{1,1}$ and $N\big( \overline{\mathcal S}^{0'}, i_k, z_k \big) + 
N\big( \overline{\mathcal W}^{0'}, i_k, z_k \big) = 1$ for all $k \in [1,n]$, and $\lbrace i_k, k \in [1,n] \rbrace$ is a set of $n$ distinct strictly positive integers.  Suppose for contradiction that for some such $k$, it holds that $N\big( \overline{\mathcal W}^{0'}, i_k,z_k \big) = 1$.  Combined with the definition of $\geq$, we would conclude that $i_k \geq \overline{\mathcal W}^{0'}_{1,1}$, a contradiction.  Thus it must be true that $N\big( \overline{\mathcal S}^{0'}, i_k,z_k \big) = 1$ for all $k$.  But as by construction transition (\ref{z14}) implies that there exists $(i*,z*) \in \overline{\mathcal S}^{0'}$ which does not belong to ${\mathcal S}^{0'}$, this (combined with the previous assertion) would imply that $|\overline{\mathcal S}^{0'}| = n + 1$, which contradicts our assumptions on ${\mathcal Q}'$.  Combining the above completes the proof.
\\\indent The proof for transition (\ref{z15}) follows for the same reason as (\ref{z3}).  For transition (\ref{z16}), it is easily verified that the only case of possible concern is that $\overline{\mathcal W}^{0'}_1 \in {\mathcal W}^{0'}$, yet $\overline{\mathcal W}^{0'}_1 \neq {\mathcal W}^{0'}_1$.  But this is impossible, as it would imply the existence of an $(i,z) \in {\mathcal W}^{0'}$ satisfying $i < \overline{\mathcal W}^{0'}_{1,1}$.  But then by the assumed  properties of ${\mathcal Q}'$, it would have to be the case that $(i,z) \in \overline{\mathcal W}^{0'}$, leading to a contradiction as there can exist no $(i,z) \in \overline{\mathcal W}^{0'}$ s.t. $i < \overline{\mathcal W}^{0'}_{1,1}$.  Combining completes the proof for transition (\ref{z16}).  The proofs for transitions (\ref{z17}) and (\ref{z18}) follows for the same reason as (\ref{z9}) and (\ref{z10}).  The proof for transition (\ref{z19}) follows for the same reason as (\ref{z4}), but applied to the queue itself.
\\\indent For the final transition (\ref{z20}), the reasoning is the same as transition (\ref{z5}), with the caveat that we must verify that $N({\mathcal W}^{0''},i,z) \leq N(\overline{\mathcal W}^{0''},i,z)$.  It is easily verified that the only case of possible concern is that $\overline{\mathcal W}^{0'}_1 \in {\mathcal W}^{0'}$.  By the same argument used in analyzing transition (\ref{z14}), we can conclude that if $\overline{\mathcal W}^{0'}_1 \in {\mathcal W}^{0'}$, then it must hold that $\overline{\mathcal W}^{0'}_1 = {\mathcal W}^{0'}_1$.  The fact that ${\mathcal W}^{0'} \neq \emptyset$ then implies that there exist $\lbrace (i_k, z_k), k = 1,\ldots,n \rbrace$ s.t. $N({\mathcal S}^{0'},i_k,z_k) = 1$ for all $k$, $\lbrace i_k, k = 1,\ldots,n \rbrace$ are distinct strictly positive integers, and $i_k < {\mathcal W}^{0'}_{1,1} = \overline{\mathcal W}^{0'}_{1,1}$ for all $k$.  But as transition (\ref{z20}) can only occur if $N\big( \overline{\mathcal S}^{0'}, 0, z \big) \geq 1$ for some $z \in \lbrace 1,2 \rbrace$, it follows that transition (\ref{z20}) can only occur if $\overline{\mathcal S}^{0'}$ contains at most $n-1$ elements from the set $\lbrace (i_k, z_k), k = 1,\ldots,n \rbrace$, and hence there exists at least one such tuple $(i^0,z^0) \in \lbrace (i_k, z_k), k = 1,\ldots,n \rbrace$ which does not belong to $\overline{\mathcal S}^{0'}$.  But our assumptions about ${\mathcal Q}'$ then imply that $(i^0,z^0)$ must belong to $\overline{\mathcal W}^{0'}$.  But this yields a contradiction, since $i^0 < \overline{\mathcal W}^{0'}_{1,1}$, yet all $(i,z) \in \overline{\mathcal W}^{0'}$ satisfy $i \geq \overline{\mathcal W}^{0'}_{1,1}$.   Combining the above completes the proof.  $\qed$
\endproof
\ \\\indent We defer the proof of Lemma\ \ref{marginals1} to the appendix.  Combining the above then completes the proof of Theorem\ \ref{ubmain1}.

\section{Skorokhod representation for upper bound $\overline{W}.$}\label{SkorSec}
To analyze the distribution of $\overline{W}$ (generally, i.e. not necessarily in the HW regime), we will show that $\overline{W}$ can be expressed as the solution to a certain generalized Skorokhod problem, on an appropriate probability space.  We will then use a result of \cite{Reed.13} to express the distribution in terms of a certain (asymptotically) tractable supremum.  We do note that our analysis also yields an upper bound for any fixed $n$, which is roughly the supremum of many coupled mean-reverting random walks, and leave as an interesting direction for future research how best to interpret and use this bound for any fixed $n$ (again, asymptotically it becomes a tractable supremum of a Gaussian process).

\subsection{Review of generalized drift Skorokhod problem.}
We begin by formally reviewing the generalized drift Skorokhod problem (as defined in \cite{Reed.13}), equivalently the linearly generalized regulator mapping (as defined  in \cite{Reed.04,Glynn.05}).  We only describe the relevant problems in the narrow context needed for our purposes, e.g. restricting to one-dimension and linear drift.  Recall that $D[0,\infty)$ denotes the appropriate space of functions, and we refer the interested reader to \cite{Whitt.02} for further details.  In general, for a mapping $Z$ from $D[0,\infty)$ to $D[0,\infty)$, i.e. from functions to functions, we use the notation $Z(X,t)$ to denote the function $Z(X)$ evaluated at time $t$, for $X \in D[0,\infty)$.  Similarly, for a mapping $Z$ from $D[0,\infty) \times {\mathcal R}$ to $D[0,\infty)$, we use the notation $Z(X,s,t)$ to denote the function $Z(X,s)$ evaluated at time $t$.
\ \\\begin{definition}[Generalized drift Skorokhod problem in 1-dimension \cite{Reed.13,Reed.04,Glynn.05}]\label{skordef}
Given a function $X \in D[0,\infty)$ with $X(0) = 0$, there exists a unique pair of functions $\Phi(X): D[0,\infty) \rightarrow D[0,\infty)$ and $\Psi(X): D[0,\infty) \rightarrow D[0,\infty)$ s.t.:
\\\begin{itemize}
\item $\Phi(X,t) = X(t) - \theta \int_0^t \Phi(X,s) ds + \Psi(X,t)$ for all $t \geq 0$;
\item $\Phi(X,t) \geq 0$ and $\Psi(X,t) \geq 0$ for all $t \geq 0$;
\item $\Phi(X,0) = \Psi(X,0) = 0$;
\item $\Psi(X,t)$ is a non-decreasing function of $t$;
\item $\int_0^{t} I\big( \Phi(X,s) > 0 \big) d \Psi(X,s) = 0$ for all $t \geq 0$.
\end{itemize}  
As is standard in the literature, $\int_0^{t} I\big( \Phi(X,s) > 0 \big) d \Psi(X,s)$ should be formally understood as a so-called Lebesgue-Stieltjes integral \cite{Hil.63}, where all such integrals appearing throughout should be similarly interpreted.
\end{definition}
\ \\\indent We next review an important representation property of $\Phi$, proven in \cite{Reed.13}.  We begin by defining an auxiliary map.
\ \\\begin{definition}[\cite{Reed.04}]\label{Reed1}
Given a function $X \in D[0,\infty)$ with $X(0) = 0$, for each $s \geq 0$, there exists a unique function $\zeta(X,s) \in D[0,\infty)$ s.t. for all $t \geq 0$,
\begin{equation}\label{Reedeq1}
\zeta(X,s,t) = X(s + t) - X(s) - \theta \int_0^t \zeta(X,s,y) dy.
\end{equation}
\end{definition}
\ \\\indent In that case, the following representation theorem is proven in \cite{Reed.13}.
\begin{lemma}[\cite{Reed.13}]\label{Reed2}
For any function $X \in D[0,\infty)$ s.t. $X(0) = 0$, for all $t \geq 0$, $\Phi(X,t) = \sup_{0 \leq s \leq t} \zeta(X,s,t-s)$.
\end{lemma}
\ \\\indent We next review an important continuity property of $\Phi$, proven in e.g. \cite{Reed.04,Glynn.05}.
\begin{lemma}\cite{Reed.04,Glynn.05} \label{psicont}
$\Phi$ is Lipschitz continuous, under the uniform norm, as a function from $D[0,\infty)$ to $D[0,\infty)$.  Namely, for each $T \geq 0$, there exists an absolute constant $C_T$ (depending implicitly on $\theta$) s.t. for all $X_1$ and $X_2 \in D[0,T]$ with $X_1(0) = X_2(0) = 0$, $\sup_{0 \leq t \leq T}\big| \Phi(X_1,t) - \Phi(X_2,t) \big| \leq C_T \sup_{0 \leq t \leq T} \big| X_1(t) - X_2(t) \big|$.
\end{lemma}

\subsection{Expressing $\overline{W}$ as a solution to the generalized drift Skorokhod problem.}

We now prove that $\overline{W}$ may be expressed as the solution to an appropriate generalized drift Skorokhod problem.  Intuitively, the arrival of dummy jobs which prevent any servers from going idle will coincide with an appropriate pushing process (appearing in an appropriate Skorokhod problem) which keeps the number in system above $n$, i.e. we will interpret $\overline{W}$ as a certain type of reflection.  Recall that a pooled equilibrium renewal process (of $n$ processes) with renewal distribution $S$ corresponds to the counting process tracking the number of events which have occured collectively in $n$ independent renewal processes, each of which has renewal distribution $S$, and for which the initial residual life is independently drawn from the equilibrium distribution of $S$ for each component renewal process.  Here we recall that the equilibrium distribution of $S$ is itself a hyper-exponential distribution, with parameters $\mu_1,\mu_2, \hat{p} = \frac{p}{\mu_1}$.  Let $\REN$ be such an explicitly constructed pooled equilibrium renewal process.  On the same probability space, let $\RES_1$ and $\RES_2$ be the stochastic processes tracking the number of these renewal processes on which the current residual life has an $Expo(\mu_1)$ and $Expo(\mu_2)$ distribution, respectively.  
It is easy to see (and well-known) that the dynamics of $\big( \RES_1, \RES_2, \REN \big)$ can be represented as a CTMC, with dynamics given by the following generator $q_{Ren}$.
\ \\
\begin{itemize}
\item $q_{Ren}\big( (\RES_1, \RES_2, \REN), (\RES_1 - 1, \RES_2 + 1, \REN + 1) \big) = \mu_1 \RES_1 p_2$;
\item $q_{Ren}\big( (\RES_1, \RES_2, \REN), (\RES_1 + 1, \RES_2 - 1, \REN + 1) \big) = \mu_2 \RES_2 p_1$;
\item $q_{Ren}\big( (\RES_1, \RES_2, \REN), (\RES_1, \RES_2, \REN + 1) \big) = \mu_1 \RES_1 p_1 + \mu_2 \RES_2 p_2$.
\end{itemize}
\ \\We note that under the above construction, letting $Bin(n,\hat{p})$ denote a binomially distributed r.v. with parameters $n$ and $\hat{p}$, $\big( \RES_1(0), \RES_2(0), \REN(0) \big)$ is equivalent in distribution to $\big( Bin(n,\hat{p}), n - Bin(n,\hat{p}), 0 \big)$.  Let $\ARR$ be a rate $\lambda$ Poisson process, independent of $\big( \RES_1, \RES_2, \REN \big)$.  Furthermore, let $\CLOCK$ be a rate 1 Poisson process, whose dependency structure w.r.t. $\ARR,\RES_1, \RES_2, \REN$ we will intentionally leave unspecified.  Supposing that we have constructed $\ARR, \REN, \CLOCK$, and $\overline{Q}$ on a common probability space, let us define $\SKOR$ (on the same probability space) to be the stochastic process s.t. for all $t \geq 0$,
$$\SKOR(t) = n^{-\frac{1}{2}} \big( \ARR(t) - \REN(t) \big) - n^{-\frac{1}{2}} \bigg( \CLOCK\big( \theta \int_0^t \overline{W}(s) ds \big) - \theta \int_0^t \overline{W}(s) ds \bigg).$$
Then our main result expressing $\overline{W}$ as a solution to the generalized drift Skorokhod problem will be as follows.

\begin{lemma}\label{skoro1}
One may construct $\overline{Q}, (\RES_1, \RES_2, \REN)$, $\ARR$, and $\CLOCK$ on the same probability space, with $\ARR$ independent of $(\RES_1, \RES_2, \REN)$, s.t. $\overline{Q}(0) = \big( \RES_1(0), \RES_2(0), \REN(0) \big)$, and w.p.1, for all $t \geq 0$, $n^{-\frac{1}{2}} \overline{W}(t) = \Phi( \SKOR , t).$  
\end{lemma}

To prove Lemma\ \ref{skoro1}, it will be helpful to again consider an expanded Markovian representation for $\overline{Q}$, in which we explicitly keep track of the number of arrivals, departures, abandonments, and dummy job arrivals.  We define a 7-dimensional CTMC $\tilde{Q} = \big( \tilde{\RES}_1, \tilde{\RES}_2, \tilde{W}, \tilde{\ARR}, \tilde{\REN}, \tilde{\ABA}, \tilde{\DMY} \big)$.  The dynamics of $\tilde{Q}$ are given by the generator $\tilde{q}$ as follows, where as before we denote the generic current state
$$\big( \tilde{\RES}_1, \tilde{\RES}_2, \tilde{W}, \tilde{\ARR}, \tilde{\REN}, \tilde{\ABA}, \tilde{\DMY} \big)$$
simply by $``\cdot"$.  Intuitively, $\tilde{\RES}_1 (\tilde{\RES}_2)$ will correspond to the number of jobs in service with residual life distributed as $Expo(\mu_1) \big( Expo(\mu_2) \big)$; $\tilde{W}$ will correspond to the number of jobs waiting in queue; $\tilde{\ARR}$ will correspond to the number of non-dummy arrivals so far; $\tilde{\REN}$ will correspond to the total number of departures so far; $\tilde{\ABA}$ will correspond to the total number of abandonments so far; and $\tilde{\DMY}$ will correspond to the total number of dummy arrivals so far.      
\\\begin{definition}[Generator for $\tilde{Q}$]
\ \\\begin{itemize}
\item $\tilde{q}\bigg( \cdot , \big( \tilde{\RES}_1 - 1, \tilde{\RES}_2 + 1, \tilde{W}, \tilde{\ARR}, \tilde{\REN} + 1, \tilde{\ABA}, \tilde{\DMY} + 1 \big) \bigg) = I(\tilde{W} = 0) \mu_1 \tilde{\RES}_1 p_2$;
\item $\tilde{q}\bigg( \cdot , \big( \tilde{\RES}_1 + 1, \tilde{\RES}_2 - 1, \tilde{W}, \tilde{\ARR}, \tilde{\REN} + 1, \tilde{\ABA}, \tilde{\DMY} + 1 \big) \bigg) = I(\tilde{W} = 0) \mu_2 \tilde{\RES}_2 p_1$;
\item $\tilde{q}\bigg( \cdot , \big( \tilde{\RES}_1, \tilde{\RES}_2, \tilde{W}, \tilde{\ARR}, \tilde{\REN} + 1, \tilde{\ABA}, \tilde{\DMY} + 1 \big) \bigg)
\\\text{\ \ \ \ \ \ \ \ \ \ \ \ } = \text{\ \ \ } I(\tilde{W} = 0) \big( \mu_1 \tilde{\RES}_1 p_1 + \mu_2 \tilde{\RES}_2 p_2 \big)$;
\item $\tilde{q}\bigg( \cdot , \big( \tilde{\RES}_1, \tilde{\RES}_2, \tilde{W} + 1, \tilde{\ARR} + 1, \tilde{\REN}, \tilde{\ABA}, \tilde{\DMY} \big) \bigg) = \lambda$;
\item $\tilde{q}\bigg( \cdot , \big( \tilde{\RES}_1 - 1, \tilde{\RES}_2 + 1, \tilde{W} - 1, \tilde{\ARR}, \tilde{\REN} + 1, \tilde{\ABA}, \tilde{\DMY} \big) \bigg) = I(\tilde{W} \geq 1) \mu_1 \tilde{\RES}_1 p_2$;
\item $\tilde{q}\bigg( \cdot , \big( \tilde{\RES}_1 + 1, \tilde{\RES}_2 - 1, \tilde{W} - 1, \tilde{\ARR}, \tilde{\REN} + 1, \tilde{\ABA}, \tilde{\DMY} \big) \bigg) = I(\tilde{W} \geq 1) \mu_2 \tilde{\RES}_2 p_1$;
\item $\tilde{q}\bigg( \cdot , \big( \tilde{\RES}_1, \tilde{\RES}_2, \tilde{W} - 1, \tilde{\ARR}, \tilde{\REN} + 1, \tilde{\ABA}, \tilde{\DMY} \big) \bigg) 
\\\text{\ \ \ \ \ \ \ \ \ \ \ \ } = \text{\ \ \ } I(\tilde{W} \geq 1) \big( \mu_1 \tilde{\RES}_1 p_1 + \mu_2 \tilde{\RES}_2 p_2 \big)$;
\item $\tilde{q}\bigg( \cdot , \big( \tilde{\RES}_1, \tilde{\RES}_2, \tilde{W} - 1, \tilde{\ARR}, \tilde{\REN}, \tilde{\ABA} + 1, \tilde{\DMY} \big) \bigg) =  \theta \tilde{W}$.
\end{itemize}
\end{definition}
\ \\We next make several observations regarding $\tilde{Q}$, and its relation to $\overline{Q}$.  Then our first observation relates $\overline{W}$ to $\tilde{W}$, and follows from a straightforward comparison of $\overline{q}$ and $\tilde{q}$, the details of which we omit.  

\begin{lemma}--\label{sameagain1}
For all non-negative integer triplets $(n_1,n_2,w)$ s.t. $n_1 + n_2 = n$ and $w = 0$, one can construct $\overline{Q}$ and $\tilde{Q}$ on a common probability space s.t. $\overline{Q}(0) = (n_1,n_2,w),$  $\tilde{Q}(0) = (n_1, n_2, w, 0, 0, 0, 0)$, and w.p.1, for all $t \geq 0$, $\overline{W}(t) = \tilde{W}(t)$.
\end{lemma}

Our second observation relates $(\tilde{\RES}_1, \tilde{\RES}_2, \tilde{\REN}, \tilde{\ARR})$
to $(\RES_1, \RES_2, \REN, {\ARR})$, and follows from a straightforward comparison of $\tilde{q}$ and $q_{ren}$, the details of which we omit.  

\begin{lemma}--\label{sameagain2}
One can construct $\big( {\RES}_1, {\RES}_2, {\REN} \big)$ and $\ARR$ on the same probability space as $\tilde{Q}$ s.t. w.p.1, for all $t \geq 0$,
$\big( \tilde{\RES}_1(t), \tilde{\RES}_2(t), \tilde{\REN}(t), \tilde{\ARR}(t) \big)$ equals
$\big( \RES_1(t), \RES_2(t), \REN(t), \ARR(t) \big)$.
\end{lemma}

Next, we observe that as the infinitesimal rate at which $\tilde{\ABA}$ jumps (always up by 1) is at all times equal to $\theta \tilde{W}$, we may conclude the following from standard random-time-change constructions (see e.g. \cite{Pang.07}), and we omit the details.

\begin{lemma}--\label{sameagain3}
One may construct $\tilde{Q}$ on the same probability space as $\CLOCK$ s.t. w.p.1, for all $t \geq 0$, $\tilde{\ABA}(t) = \CLOCK\big( \theta \int_0^t \tilde{W}(s) ds \big)$, where this integral is w.p.1 well-defined as a Riemann integral for all $t \geq 0$.
\end{lemma}

Finally, we note that a straightforward examination of $\tilde{q}$ and proof by induction, the details of which we omit, reveals the following relationships between $\tilde{W}, \tilde{\ARR}, \tilde{\REN}, \tilde{\ABA}$, and $\tilde{\DMY}$.

\begin{lemma}--\label{sameagain4}
For all non-negative integer triplets $(n_1,n_2,w)$ s.t. $n_1 + n_2 = n$ and $w = 0$, one can construct $\tilde{Q}$ s.t. $\tilde{Q}(0) = (n_1, n_2, w, 0, 0, 0, 0)$, and w.p.1, for all $t \geq 0$:
\begin{itemize}
\item $\tilde{W}(t) = \tilde{\ARR}(t) - \tilde{\REN}(t) - \tilde{\ABA}(t) + \tilde{\DMY}(t)$;
\item $\int_0^t I\big( \tilde{W}(s) > 0 \big) d\tilde{\DMY}(s) = 0$.
\end{itemize}
\end{lemma}

Combining Lemmas\ \ref{sameagain1}\ -\ \ref{sameagain4} and noting that all relevant constructions can themselves be implemented together on a common probability space, the fact that w.p.1 $\tilde{\DMY}$ is non-decreasing and initialized at 0 for appropriate initial conditions, and Definition\ \ref{skordef} then completes the proof of the desired result Lemma\ \ref{skoro1} (after scaling through by $n^{-\frac{1}{2}}$).

\section{Asymptotic analysis of $\overline{W}$ and proof of upper bound for Theorem\ \ref{main1}.}\label{AsupSec}
In this section we provide an asymptotic analysis of $\overline{W}$ in the HW regime, and complete the proof of one direction of our main result Theorem\ \ref{main1}.  We proceed as follows.  First, we prove an asymptotic version of Lemma\ \ref{skoro1} in the HW regime.  Second, we prove that the associated limit may be rewritten in terms of an integrated OU process.  Third, we use this OU representation, along with the supremum representation for the generalized Skorokhod problem from Lemma\ \ref{Reed2}, to prove that in the HW regime $\overline{W}$ may be expressed as the supremum of an explicit Gaussian process.  Finally, we analyze this Gaussian process, and combine with the interchange-of-limits results of \cite{Dieker.13,DDG14} mentioned in Section\ \ref{IntroSec} to complete the proof of one direction of our main result Theorem\ \ref{main1}.

\subsection{Asymptotic version of Lemma\ \ref{skoro1}.}
We now prove an asymptotic version of Lemma\ \ref{skoro1}, in the HW regime.  We begin by introducing a superscript notation to make the dependence of some of the previously defined stochastic processes on $n$ more explicit.  Let $\ARR^n$ denote the corresponding Poisson process when $\lambda = \lambda_n = n + B n^{\frac{1}{2}}$; $\REN^n$ denote the corresponding pooling of $n$ independent equilibrium renewal processes; $\RES^n_1$ and $\RES^n_2$ denote the number of these renewal processes with the corresponding residual life distribution (respectively); $\overline{Q}^n$ and $\overline{W}^n$ denote the corresponding stochastic processes when there are $n$ servers and $\lambda = \lambda_n$; all constructed on a common probability space (with $\CLOCK$) using the construction given in Lemma\ \ref{skoro1}.  Note that under this construction, $\overline{Q}^n(0) = \big( \RES^n_1(0), \RES^n_2(0), \REN^n(0) \big),$ i.e. there are exactly $n$ jobs in system with the initial residual life of each job drawn i.i.d from the equilibrium distribution of $S$, i.e. $Expo(\mu_1)$ with probability $\hat{p} = \frac{p}{\mu_1}$ and $Expo(\mu_2)$ with probability $\frac{1-p}{\mu_2}$.  On this same probability space, let us define (in analogy with our previous definition) $\SKOR^n$ to be the stochastic process s.t. for all $t \geq 0$,
$$\SKOR^n(t) = n^{-\frac{1}{2}} \big( \ARR^n(t) - \REN^n(t) \big) - n^{-\frac{1}{2}} \bigg( \CLOCK\big( \theta \int_0^t \overline{W}^n(s) ds \big) - \theta \int_0^t \overline{W}^n(s) ds \bigg).$$
\\\indent Let $\REN^{\infty}$ denote the unique centered continuous Gaussian process s.t. $\E[\REN^{\infty}(s)\REN^{\infty}(t)] = \E[\big( \REN^1(s) - s \big) \big( \REN^1(t) - t \big)]$ for all $0 \leq s \leq t$, i.e. the covariance structure is the same as that of a single equilibrium renewal process, where existence of such a process is proven in \cite{Whitt.85}.  Let ${\mathcal B}^0$ denote a standard Brownian motion, independent of $\REN^{\infty}$.  Let $\SKOR^{\infty}$ denote the stochastic process s.t. for all $t \geq 0$, $\SKOR^{\infty}(t) = {\mathcal B}^0(t) - \REN^{\infty}(t) + B t$.
\\\indent Then the main asymptotic result of this section will be the following.  

\begin{lemma}\label{xcon1}
For any $T \in [0,\infty)$, the sequence of processes $\lbrace n^{-\frac{1}{2}} \overline{W}^n(t)_{0 \leq t \leq T} , n \geq 1 \rbrace$ converges weakly to $\Phi(\SKOR^{\infty})_{0 \leq t \leq T}$ in the space $D[0,T]$ under the $J_1$ topology.
\end{lemma}

To prove Lemma\ \ref{xcon1}, we begin by recalling a known weak convergence result from \cite{GG13}.

\begin{lemma}[\cite{GG13}]\label{weaklem1}
For any $T \in [0,\infty)$, the sequence of processes $\big\lbrace n^{-\frac{1}{2}} \big( \ARR^n(t) - \REN^n(t) \big)_{0 \leq t \leq T}, n \geq 1 \big\rbrace$ converges weakly to $\SKOR^{\infty}(t)_{0 \leq t \leq T}$
in the space $D[0,T]$ under the $J_1$ topology.
\end{lemma}

We next prove a result showing that the process corresponding to the second term in the definition of $\SKOR^n$ converges weakly to 0.  Namely, let $\ERR^n$ denote the stochastic process s.t. for all $t \geq 0$,
$$\ERR^n(t) = n^{-\frac{1}{2}} \bigg( \CLOCK\big( \theta \int_0^t \overline{W}^n(s) ds \big) - \theta \int_0^t \overline{W}^n(s) ds \bigg).$$
Then we prove the following, deferring the proof to the appendix.

\begin{lemma}*\label{weaklem2}
For any $T \in [0,\infty)$, the sequence of processes $\lbrace \ERR^n(t)_{0 \leq t \leq T}, n \geq 1 \rbrace$ converges weakly to $0_{0 \leq t \leq T}$ in the space $D[0,T]$ under the $J_1$ topology.
\end{lemma}

Combining Lemmas\ \ref{weaklem1} and \ref{weaklem2} with Lemmas\ \ref{skoro1} and \ref{psicont}, the fact that convergence in the $J_1$ metric to a continuous process implies convergence in the uniform metric, and the continuous mapping theorem (see e.g. \cite{Whitt.02}) then completes the proof of Lemma\ \ref{xcon1}.

\subsection{Relating $\SKOR^{\infty}$ to an integrated OU process.}
In this section we prove that $\SKOR^{\infty}$ is equal in distribution to the sum of a so-called time-integrated OU process (to be defined) plus an independent Brownian motion.  Let $\OU$ denote the unique w.p.1 continuous centered Gassian process s.t. $\E[ \OU(s) \OU(t)] = p (1 - p) (\mu_1 \mu_2)^{-1} \exp\big( - \mu_1 \mu_2 (t - s) \big)$ for $0 \leq s \leq t$, equivalently an appropriate stationary OU process, i.e. the strictly stationary process which is the unique such strong solution to the SDE
\begin{equation}\label{thesde1}
d\OU(t) = - \mu_1 \mu_2 \OU(t) dt + \big(2 p (1-p) \big)^{\frac{1}{2}} d{\mathcal B}^0(t),
\end{equation}
with ${\mathcal B}^0$ the previously defined standard Brownian motion.  Here we refer the interested reader to \cite{Doob.42} for further details regarding existence and various properties of this Gaussian process, and to \cite{Masuda.04} for further review of associated SDEs.
As continuous functions are Riemann integrable, it follows that $\OU$ is w.p.1 Riemann integrable over compact time intervals.  Thus we may define $\IOU$ to be the stochastic process s.t. $\IOU(t) = \int_0^t \OU(y) dy$ for all $t \geq 0$.  Such a process is known as a time-integrated OU process, and is well-studied, e.g. in \cite{Barn.03} where the associated variance is computed explicitly (as stated below).  Furthermore, as $\IOU$ has stationary increments (see e.g. \cite{GG13}), $\E[ \IOU^2(t-s) ] = \E[\IOU^2(t)] + \E[\IOU^2(s)] - 2 \E[\IOU(s) \IOU(t)]$ for $0 \leq s \leq t$.  It follows that the variance calculation in \cite{Barn.03}, combined with some straightforward algebra the details of which we omit, yields the following.

\begin{lemma}[\cite{Barn.03}--]\label{intouprop1}
$\IOU$ is a continuous centered Gaussian process with stationary increments s.t. for all $0 \leq s \leq t$, 
$\E\big[\IOU(s)\IOU(t)\big]$ equals
$$p (1 - p) \big(\mu_1 \mu_2)^{-3} \bigg(2 \mu_1 \mu_2 s + \exp( - \mu_1 \mu_2 s) + \exp(-\mu_1 \mu_2 t) - \exp\big(- \mu_1 \mu_2 (t - s) \big) - 1 \bigg).$$
\end{lemma}

Suppose that $\IOU$ and $\OU$ have been constructed on the same probability space as standard Brownian motions ${\mathcal B}^0$ and ${\mathcal B}^1$, with ${\mathcal B}^0$ the driving Brownian motion of $\OU$, $\IOU(t) = \int_0^t \OU(s) ds$ for all $t \geq 0$, and $\OU, \IOU, {\mathcal B}^0$ independent of ${\mathcal B}^1$.  On the same probability space, let us define the stochastic process $\Lambda$ s.t. for all $t \geq 0$, 
$$\Lambda(t) = - (\mu_1 - \mu_2) \IOU(t) + 2^{\frac{1}{2}} {\mathcal B}^1(t) + B t.$$  
Then our next result formalizes the relationship between $\SKOR^{\infty}$ and $\Lambda,$ and we defer the proof to the appendix.

\begin{lemma}*\label{samedist1}
$\lbrace \Lambda(t), t \geq 0 \rbrace$ has the same distribution, at the process level, as $\lbrace \SKOR^{\infty}(t), t \geq 0 \rbrace$.
\end{lemma}

\subsection{Stochastic calculus to analyze $\Lambda$.}

In this section, we use stochastic calculus to rewrite $\zeta\big( \Lambda , s \big)$ in a form amenable to analyzing the supremum arising in Lemma\ \ref{Reed2}, ultimately relating $\Phi( \Lambda, t )$ to the supremum of a certain Gaussian process.  Our main result will be the following.  Let us define ${\mathcal G}$ to be the stochastic process s.t. for all $t \geq 0,$
$${\mathcal G}(t) =  -(\mu_1 - \mu_2) \int_0^t \exp\big( - \theta y \big) \OU(y) dy + 2^{\frac{1}{2}} \int_0^t \exp(- \theta y) d{\mathcal B}^1(y) +  \frac{B}{\theta}\bigg(1 - \exp\big(- \theta t \big) \bigg),$$
with $\OU$ independent of ${\mathcal B}^1$.  Then we prove the following.
\begin{lemma}\label{nicesup1}
For all $t \geq 0$, 
$\Phi( \Lambda, t )$ has the same distribution as $\sup_{0 \leq s \leq t} {\mathcal G}(s)$.
\end{lemma}
In the rest of this subsection we prove Lemma\ \ref{nicesup1}.  We begin by providing some more explicit representations/constructions for several previously defined stochastic processes, which follow immediately from (\ref{thesde1}), standard and well-known results for OU processes \cite{Ross.09}, a few simple reindexings, and our definitions (especially Definition\ \ref{Reed1}), and we omit the details.  We note that in this and remaining sections, we freely use standard notations, definitions, and basic results from the elementary theory of stochastic calculus, and refer the interested reader to \cite{Kar.12} for further details.  Also, for a general stochastic process $X$ and $s \geq 0$, we let $X_s$ denote the stochastic process s.t. $X_s(t) = X(s + t) - X(s)$ for all $t \geq 0$.  Similarly, for a general stochastic process $X$ and $s \geq 0$, we let $X_{+,s}$ denote the stochastic process s.t. $X_{+,s}(t) = X(s + t)$ for all $t \geq 0$.  For a generic r.v $X$, let $Var[X]$ denote the variance of $X$.  Let ${\mathcal B}^3$ and ${\mathcal B}^2$ denote two independent Brownian motions, constructed on the same probability space, which satisfy $Var[{\mathcal B}^3(1)] = 2 p (1 - p), Var[{\mathcal B}^2(1)] = 2.$  Then the desired representations are as follows, where we remind the reader that the mapping $\zeta$, from $D[0,\infty) \times {\mathcal R}^+$ to $D[0,\infty)$, was previously formalized in Definition\ \ref{Reed1}. 
    
\begin{lemma}--\label{thesde45}
One may construct $\Lambda$ and $\OU$ on the same probability space as ${\mathcal B}^3$ and ${\mathcal B}^2$, s.t. w.p.1, for all $s,t \geq 0$,
\begin{equation}\label{thesde5}
\OU_{+,s}(t) = \exp(- \mu_1 \mu_2 t) \OU_{+,s}(0) + \int_0^t \exp\big( - \mu_1 \mu_2 (t - y) \big) d{\mathcal B}^3_s(y);
\end{equation}
\begin{equation}\label{thesde4}
\zeta(\Lambda, s, t) = -(\mu_1 - \mu_2) \int_0^t \OU_{+,s}(y) dy + {\mathcal B}^2_s(t) + B t - \theta \int_0^t \zeta(\Lambda, s, y) dy.
\end{equation}
\end{lemma}

From Lemma\ \ref{thesde45} and basic definitions in stochastic calculus, see e.g. \cite{Kar.12}, we conclude the following.  Let $J \stackrel{\Delta}{=} \begin{bmatrix} \mu_1 \mu_2 & 0 \\ \mu_1 - \mu_2 & \theta \end{bmatrix}$.  Also, for two appropriately dimensioned matrices $M_1,M_2$, let $M_1 \bullet M_2$ denote the standard dot product of the two matrices.

\begin{corollary}--\label{thesdecor}
Under the same construction as Lemma\ \ref{thesde45}, w.p.1, for all $s,t \geq 0$, $\OU_{+,s}$ and $\zeta(\Lambda , s)$ satisfy the system of SDE
\begin{equation}\label{thesde}
\left[ \begin{array}{c} d\OU_{+,s}(t) \\ d\zeta(\Lambda , s , t) \end{array} \right] 
= - J \bullet 
\left[ \begin{array}{c} \OU_{+,s}(t) \\ \zeta(\Lambda, s, t) \end{array} \right]
+
\left[ \begin{array}{c} d{\mathcal B}^3_s(t) \\ d{\mathcal B}^2_s(t) + Bdt \end{array} \right],
\end{equation}
with appropriate boundary conditions.
\end{corollary}

Fortunately, the system of SDE (\ref{thesde}) is well-studied and has an explicit solution, the so-called multi-dimensional OU process (see e.g. \cite{Masuda.04}).  Before recalling the details of such processes, it will be helpful to review some additional matrix-related notations.
Recall that for a matrix $M$ and scalar $z$, $z M$ denotes the same matrix $M$, but with each component multiplied by $z$.  For a square matrix $M$, we let $M^0$ denote the identity matrix of the same dimensions, and $M^k$ denote $M \bullet M^{k-1}$ for $k \geq 1$.  For a matrix $M$, let $\exp(M)$ denote the standard matrix exponential, i.e. the matrix equivalent to $\sum_{k=0}^{\infty} \frac{M^k}{k!}$.  Also, for an $n$ by $m$ matrix $M$ and $i \in [1,n], j \in [1,m]$, let $M_{i,j}$ denote the entry in row $i$, column $j$ of $M$.

\begin{lemma}[\cite{Masuda.04}]\label{thesol1lemma}
Under the same construction as Lemma\ \ref{thesde45}, w.p.1, for all $s,t \geq 0$, the system of SDE (\ref{thesde}) has the following unique strong solution:
\begin{eqnarray}
\ &\ &\ \left[ \begin{array}{c} \OU_{+,s}(t) \\ \zeta(\Lambda, s, t) \end{array} \right] 
= 
\left[ \begin{array}{c} 0 \\ \frac{B}{\theta}\big(1 - \exp(- \theta t) \big) \end{array} \right] 
+
\exp(- t J) \bullet \left[ \begin{array}{c} \OU_{+,s}(0) \\ \zeta(\Lambda, s, 0) \end{array} \right] \label{thesol1}
\\&\ &\ \ \ \ \ \ \ \ \ \ \ \ \ \ \ \ \ \ \ \ +
\left[ \begin{array}{c} \int_0^t \exp\big( - (t - y) J \big)_{1,1} d{\mathcal B}^3_s(y) + \int_0^t \exp\big( - (t - y) J \big)_{1,2} d{\mathcal B}^2_s(y) \\ \int_0^t \exp\big( - (t - y) J \big)_{2,1} d{\mathcal B}^3_s(y) + \int_0^t \exp\big( - (t - y) J \big)_{2,2} d{\mathcal B}^2_s(y) \end{array} \right]. \nonumber
\end{eqnarray}
\end{lemma}

We now simplify (\ref{thesol1}) through a series of lemmas.  First, we compute $\exp(- z J)$ explicitly for $z \geq 0$, deferring the proof to the appendix.  Here we also note that it is easily verified that our assumptions (i.e. $\E[S] = 1$ and Assumption\ \ref{A1}) imply that $\mu_1 \mu_2 > \theta$, as $\frac{p}{\mu_1} + \frac{1-p}{\mu_2} = 1$ implies $\mu_1 \mu_2 = p \mu_2 + (1 - p) \mu_1$.

\begin{lemma}*\label{Qexp1}
For all $z \geq 0$, $$\exp(- z J) = 
\begin{bmatrix} \exp\big( - \mu_1 \mu_2 z \big) & 0 \\ \frac{\mu_1 - \mu_2}{\mu_1 \mu_2 - \theta} \big( \exp( - \mu_1 \mu_2 z) - \exp( - \theta z) \big) & \exp(- \theta z)\end{bmatrix}.$$
\end{lemma}

Combining Lemmas\ \ref{thesol1lemma} and \ref{Qexp1}, we conclude the following.

\begin{lemma}\label{theeq2blemma}
Under the same construction as Lemma\ \ref{thesde45}, w.p.1, for all $0 \leq s \leq t$,
\begin{equation}\label{theeq2b}
\OU_{+,s}(t) = \exp(- \mu_1 \mu_2 t) \OU_{+,s}(0) + \int_0^t \exp\big( - \mu_1 \mu_2 (t - y) \big) d{\mathcal B}^3_s(y),
\end{equation}
and
\begin{eqnarray}
\zeta(\Lambda, s, t - s) &=& \frac{\mu_1 - \mu_2}{\mu_1 \mu_2 - \theta} \bigg( \exp\big(- \mu_1 \mu_2 (t-s)\big) - \exp\big(- \theta (t-s) \big) \bigg) \OU_{+,s}(0)  \label{theeq2}
\\&\ &\ \ \ \ \ +\ \ \ \frac{\mu_1 - \mu_2}{\mu_1 \mu_2 - \theta} \int_0^{t-s} \bigg( \exp\big( - \mu_1 \mu_2 (t - s - y) \big) - \exp\big( - \theta (t - s - y) \big) \bigg) d{\mathcal B}^3_s(y) \nonumber
\\&\ &\ \ \ \ \ +\ \ \ \int_0^{t-s} \exp\big( - \theta (t - s - y) \big) d{\mathcal B}^2_s(y) + \frac{B}{\theta}\bigg(1 - \exp\big(- \theta (t-s) \big) \bigg). \nonumber
\end{eqnarray}
\end{lemma}
Next, we will combine Lemma\ \ref{theeq2blemma} with Ito's lemma to prove the following, deferring the proof to the appendix.
\begin{lemma}*\label{itorewrite1}
Under the same construction as Lemma\ \ref{thesde45}, w.p.1, for all $0 \leq s \leq t$,
\begin{eqnarray*}
\zeta(\Lambda, s, t - s) &=& -(\mu_1 - \mu_2) \int_s^t \exp\big( - \theta(t - y) \big) \OU(y) dy  
\\&\ &\ \ \ \ \ \ +\ \ \ \int_s^t \exp\big( - \theta(t - y) \big) d{\mathcal B}^2(y) + \frac{B}{\theta}\bigg(1 - \exp\big(- \theta(t-s) \big) \bigg).
\end{eqnarray*}
\end{lemma}

The desired result Lemma\ \ref{nicesup1} then follows from Lemma\ \ref{Reed2}, Lemma\ \ref{itorewrite1}, the facts that (by reversibility): 1. $\OU(s)_{0 \leq s \leq t}$ has the same distribution (at the process level) as $\OU(t - s)_{0 \leq s \leq t}$ and 2. ${\mathcal B}^2(s)_{0 \leq s \leq t}$ has the same distribution (at the process level) as $2^{\frac{1}{2}} \big( {\mathcal B}^1(t) - {\mathcal B}^1(t - s) \big)_{0 \leq s \leq t}$, and a simple reindexing, and we omit the details.

\subsection{Gaussian process theory to analyze $\sup_{0 \leq s \leq t} {\mathcal G}(s)$.}

In this section we use Dudley's celebrated entropy bound to prove that the expected all-time supremum of ${\mathcal G}$ is finite.  This will later allow us to relate the large deviations behavior of the all-time supremum of ${\mathcal G}$ to the maximum variance of ${\mathcal G}$ at any one time, which we will explicitly compute, using e.g. Borell's inequality.

We begin by recalling some well-known results from the theory of Gaussian processes, and refer the interested reader to \cite{Ad90,Dud.67,Dud.16,Tal.91} for further details.  Let us define $\overline{\mathcal G}$ to be the stochastic process s.t. for all $t \geq 0$, 
$$\overline{\mathcal G}(t) = {\mathcal G}(t) - \frac{B}{\theta}\bigg(1 - \exp\big(- \theta t \big) \bigg).$$  Then one may easily verify the following.
\begin{observation}.\label{yobs1}
W.p.1, $\overline{\mathcal G}$ is a centered, continuous Gaussian process, where here continuity is on the space ${\mathcal R}^+$ under the standard Euclidean metric, and $\overline{\mathcal G}(0) = 0$.
\end{observation}
For the Gaussian process $\overline{\mathcal G}$, we can define the so-called canonical metric (technically in general a pseudo-metric, although in our special case it will anyways be a metric) for ${\mathcal R}^+$ w.r.t. $\overline{\mathcal G}$, $d_{\overline{\mathcal G}},$ s.t. for $x,y \in {\mathcal R}^+$, $d_{\overline{\mathcal G}}(x,y) = \bigg(\E\bigg[ \big(\overline{\mathcal{G}}(x) - \overline{\mathcal G}(y) \big)^2 \bigg] \bigg)^{\frac{1}{2}}$.  For a given $x \in {\mathcal R}^+$ and $\epsilon > 0$, let $B_{\overline{\mathcal G}}(x,\epsilon) \stackrel{\Delta}{=} \lbrace z \in {\mathcal R}^+: d_{\overline{\mathcal G}}(z,x) \leq \epsilon \rbrace$, i.e. the closed ball of radius $\epsilon$ centered at $x$ under the $\overline{\mathcal G}$ pseudo-metric.  For a given $\epsilon > 0$, let $N_{\overline{\mathcal G}}(\epsilon)$ denote the smallest $k \geq 1$ for which there exist $x_1,\ldots,x_k \in {\mathcal R}^+$ satisfying ${\mathcal R}^+ \subseteq \bigcup_{i=1}^k B_{\overline{\mathcal G}}(x_i, \epsilon)$, i.e. the smallest number of closed $\epsilon$-balls needed to cover ${\mathcal R}^+$ under the $d_{\overline{\mathcal G}}$ pseudo-metric, where we set this value to $\infty$ if no finite number of such balls suffices.  Then the celebrated entropy bound (of Sudakov and Dudley, later generalized by Talagrand and others) for the expected value of the supremum of a centered Gaussian process (as customized to our setting) is as follows.  We note that in general the relevant bounds only hold for the supremum over any fixed countable subset, and/or for an appropriate version of the associated stochastic process.  However, in light of the continuity ensured by Observation\ \ref{yobs1}, those results imply the following in our setting.  We also note that although ${\mathcal R}^+$ will not in general be compact under the metrics we consider, as we will see it will be totally bounded under the pseudo-metric $\overline{\mathcal G}$, which will suffice for our purposes.    
\begin{lemma}[\cite{Ad90,Dud.67,Dud.16,Tal.91}]\label{entropybound}
$\E\big[\sup_{t \in {\mathcal R}^+} \overline{\mathcal G}(t)\big]$ is finite if $\int_0^{\infty} \log^{\frac{1}{2}}\big(N_{\overline{\mathcal G}}(\epsilon)\big) d\epsilon < \infty$.
\end{lemma}
We now prove that the associated entropy integral is indeed finite, i.e. we prove the following.

\begin{lemma}\label{entropybound2}
$\int_0^{\infty} \log^{\frac{1}{2}}\big(N_{\overline{\mathcal G}}(\epsilon)\big) d\epsilon < \infty$.
\end{lemma}

En route to proving Lemma\ \ref{entropybound2}, and for use in later proofs, it will be helpful to explicitly compute the covariance of $\overline{\mathcal G}$.  This is accomplished in the following lemma, whose proof we defer to the appendix.
\begin{lemma}*\label{covlemma1}
For all $0 \leq s < t$, $\E\big[ \big( \overline{\mathcal G}(t) - \overline{\mathcal G}(s) \big)^2 \big]$ equals
\begin{eqnarray*}
\ &\ &\ \bigg( \frac{ p (1 - p) (\mu_1 - \mu_2)^2 }{ \theta \mu_1 \mu_2 (\mu_1 \mu_2 + \theta) } + \theta^{-1} \bigg) \exp(- 2 \theta s) - \bigg( \frac{ p (1 - p) (\mu_1 - \mu_2)^2 }{ \theta \mu_1 \mu_2 (\mu_1 \mu_2 - \theta) } + \theta^{-1} \bigg) \exp(- 2 \theta t)
\\&\ &\ \ \ \ \ \ \ \ +\ \ \ \frac{ 2 p (1 - p) (\mu_1 - \mu_2)^2 }{ \mu_1 \mu_2 (\mu_1 \mu_2 - \theta)(\mu_1 \mu_2 + \theta) } \exp\big( (\mu_1 \mu_2 - \theta) s \big) \exp\big( - (\mu_1 \mu_2 + \theta) t \big).
\end{eqnarray*}
\end{lemma}

We now establish some key monotonicity properties for the covariance of $\overline{\mathcal G}$, for use in later proofs, deferring all proofs to the appendix.  
\begin{corollary}*\label{monovarcor}
For each fixed $s \in {\mathcal R}^+$, $\E\big[ \big( \overline{\mathcal G}(t) - \overline{\mathcal G}(s) \big)^2 \big]$ is a strictly increasing function (of $t$) on $[s,\infty)$.  Also, for each fixed $t \in {\mathcal R}^+$, $\E\big[ \big( \overline{\mathcal G}(t) - \overline{\mathcal G}(s) \big)^2 \big]$ is a strictly decreasing function (of $s$) on $[0,t)$.
\end{corollary}

Note that in light of the strictness established in Corollary\ \ref{monovarcor}, we find that $({\mathcal R}^+, d_{\overline{\mathcal G}})$ is in fact a metric space.  We also note that Lemma\ \ref{covlemma1} implies that $\E\big[ \big( \overline{\mathcal G}(t) - \overline{\mathcal G}(s) \big)^2 \big]$ is jointly continuous in $s$ and $t$, which is a requirement for some of the results from the theory of Gaussian processes which we will be applying.  We now use Corollary\ \ref{monovarcor} to compute the diameter of this metric space, and explicitly bound the associated metric, deferring all proofs to the appendix.
\begin{lemma}*\label{diameter}
$\lim_{t \rightarrow \infty} \E\big[ \overline{\mathcal G}^2(t) \big] = \frac{\theta + \mu_1 + \mu_2 - 1}{\theta(\theta + \mu_1 \mu_2)}$, where in light of Corollary\ \ref{monovarcor} this further equals $\sup_{t \geq 0} \E\big[ \overline{\mathcal G}^2(t) \big].$
In addition, for all $s,t \in {\mathcal R}^+$, it holds that $d_{\overline{\mathcal G}}(s,t) \leq \bigg( \frac{\theta + \mu_1 + \mu_2 - 1}{\theta(\theta + \mu_1 \mu_2)} \bigg)^{\frac{1}{2}},$ i.e. the diameter of the metric space $({\mathcal R}^+, d_{\overline{\mathcal G}})$ equals $\bigg( \frac{\theta + \mu_1 + \mu_2 - 1}{\theta(\theta + \mu_1 \mu_2)} \bigg)^{\frac{1}{2}}.$  In addition, for all $0 \leq s \leq t$, 
$$d_{\overline{\mathcal G}}(s,t) \leq 4 \bigg( \frac{ p (1 - p) (\mu_1 - \mu_2)^2 }{ \mu_1 \mu_2 (\mu_1 \mu_2 - \theta) \theta } + \theta^{-1} \bigg)^{\frac{1}{2}} \bigg( 1 - \exp\big( - (\mu_1 \mu_2 + \theta)(t - s) \big) \bigg)^{\frac{1}{2}} \exp(- \theta s).$$
\end{lemma}

With Lemma\ \ref{diameter} in hand, we now complete the proof of Lemma\ \ref{entropybound2}.
\proof{Proof of Lemma\ \ref{entropybound2}:}
Let $C_0 \stackrel{\Delta}{=} 4 \bigg( \frac{ p (1 - p) (\mu_1 - \mu_2)^2 }{ \mu_1 \mu_2 (\mu_1 \mu_2 - \theta) \theta } + \theta^{-1} \bigg)^{\frac{1}{2}} (1 + \mu_1 \mu_2 + \theta)^{\frac{1}{2}}$.  For any fixed $\epsilon > 0$, let $n_{\epsilon} \stackrel{\Delta}{=} \lceil \frac{C_0^3}{\theta \epsilon^3} \rceil$.  For $k \in \lbrace 0,\ldots,n_{\epsilon} \rbrace$, let $x_{k,\epsilon} \stackrel{\Delta}{=} k \frac{\epsilon^2}{C_0^2}$.  We now prove that for all $\epsilon > 0$,
\begin{equation}\label{iscover}
{\mathcal R}^+ \subseteq \bigcup_{i = 0}^{n_{\epsilon}} B_{\overline{\mathcal G}}(x_{i,\epsilon} , \epsilon).
\end{equation}
We first prove that for all $t \geq x_{n_{\epsilon},\epsilon}$, it holds that $t \in B_{\overline{\mathcal G}}(x_{n_{\epsilon},\epsilon} , \epsilon)$.  Indeed, first note that 
$$
x_{n_{\epsilon},\epsilon}\ \ \ \geq\ \ \ \frac{1}{\theta} \times \frac{C_0}{\epsilon}\ \ \ \geq\ \ \ \frac{1}{\theta} \log\big( \frac{C_0}{\epsilon} \big),
$$
the final inequality following from the fact that $\log(x) < x$ for all $x > 0$.
It follows that for $t \geq x_{n_{\epsilon},\epsilon}$, by Lemma\ \ref{diameter},
$$
d_{\overline{\mathcal G}}\big( x_{n_{\epsilon},\epsilon} , t \big)\ \ \ \leq\ \ \ C_0 \exp\big( - \theta x_{n_{\epsilon},\epsilon} \big)\ \ \ \leq\ \ \ \epsilon,
$$
completing the proof.
\\\indent Next, we prove that for all $k \in [0, n_{\epsilon}]$ and $t \in [x_{k,\epsilon}, x_{k+1,\epsilon}]$, it holds that $t \in B_{\overline{\mathcal G}}(x_{k,\epsilon} , \epsilon)$.  
Indeed, note that as $1 - \exp(- x) \leq x$ for all $x \in {\mathcal R}^+$, it follows from Lemma\ \ref{diameter} that for all $0 \leq s \leq t$, $d_{\overline{\mathcal G}}(s,t) \leq C_0 (t - s)^{\frac{1}{2}}.$
Thus for all $k \in [0,n_{\epsilon}]$ and $t \in [x_{k,\epsilon}, x_{k+1,\epsilon}]$, 
\begin{eqnarray*}
d_{\overline{\mathcal G}}(x_{k,\epsilon},t) &\leq& C_0 \big( x_{k+1,\epsilon} - x_{k,\epsilon} \big)^{\frac{1}{2}}
\\&=& C_0 \big( \frac{\epsilon^2}{C_0^2} \big)^{\frac{1}{2}}\ \ \ =\ \ \ \epsilon,
\end{eqnarray*}
completing the proof.
\\\indent Combining the above completes the proof of (\ref{iscover}), and implies that $N_{\overline{\mathcal G}}(\epsilon) \leq \frac{C_0^3}{\theta \epsilon^3} + 2$ for all $\epsilon > 0$.  Combining with the fact that Lemma\ \ref{diameter} implies that 
$N_{\overline{\mathcal G}}(\epsilon) = 1$ for all $\epsilon \geq \big( \frac{\theta + \mu_1 + \mu_2 - 1}{\theta(\theta + \mu_1 \mu_2)} \big)^{\frac{1}{2}}$, it follows that to complete the proof of Lemma\ \ref{entropybound2}, it suffices to demonstrate that 
$$\int_0^{\big(\frac{\theta + \mu_1 + \mu_2 - 1}{\theta(\theta + \mu_1 \mu_2)}\big)^{\frac{1}{2}}} \log^{\frac{1}{2}}\bigg( \frac{C_0^3}{\theta \epsilon^3} + 2 \bigg) d\epsilon < \infty.$$
Combining with a straightforward exercise in calculus, the details of which we omit, then completes the proof.  $\qed$
\endproof

\subsection{Interchange-of-limits results from \cite{Dieker.13,DDG14}.}
The final ingredient needed in our proofs is the interchange-of-limits results from \cite{Dieker.13,DDG14}, which will allow us to translate results for $\sup_{0 \leq s \leq t} {\mathcal G}(s)$ over finite intervals $[0,t]$ to results for the stationary measure of interest.  Let $\hat{Q}^{n} = \big( \hat{N}^n_1, \hat{N}^n_2, \hat{W}^n \big)$ denote the rescaled CTMC $n^{-\frac{1}{2}} \big( Q^n - ( \frac{p}{\mu_1} n , \frac{1-p}{\mu_2} n, 0 ) \big),$ i.e. the normalized and rescaled (by $n^{-\frac{1}{2}}$) CTMC corresponding to the $M/H_2/n + M$ queue with arrival rate $\lambda_n = n + B n^{\frac{1}{2}}$, service distribution $S$ which is hyper-exponentially distributed with parameters $\mu_1,\mu_2,p$, and abandonment rate $\theta$.  Furthermore, we remind the reader that $\hat{Q}^n(0)$ is equivalent in distribution to $\big( Bin(n,\hat{p}) - n \hat{p}, n - Bin(n,\hat{p}) - n (1 - \hat{p}), 0 \big)$, where $\hat{p} = \frac{p}{\mu_1}$ and $Bin(n,\hat{p})$ is the corresponding binomially distributed r.v.  Note that standard results concerning weak convergence of the binomial to the normal imply that $\lbrace \hat{Q}^n(0), n \geq 1 \rbrace$ converges weakly to $\bigg( \big(p (1- p) \big)^{\frac{1}{2}} {\mathcal B}^0(1) , - \big(p (1- p) \big)^{\frac{1}{2}} {\mathcal B}^0(1) , 0 \bigg)$, where we remind the reader that ${\mathcal B}^0, {\mathcal B}^1$ are two independent standard Brownian motions.  Also, it is well-known that for each fixed $n$, $\lbrace \hat{Q}^n(t), t \geq 0 \rbrace$ converges weakly to a r.v. $\hat{Q}^n(\infty) = \big( \hat{N}^n_1(\infty), \hat{N}^n_2(\infty), \hat{W}^n(\infty) \big)$, corresponding to the steady-state of the corresponding (normalized) $M/H_2/n + M$ queue \cite{DDG14}.

\subsubsection{Review of weak convergence results of \cite{Dai.10}.}
Before reviewing the relevant interchange-of-limits results, we will have to formalize the weak convergence results of \cite{Dai.10} in greater depth.  The relevant weak convergence result is as follows, and is proven in \cite{Dai.10}.
\begin{lemma}[\cite{Dai.10}]\label{funweak1}
There exists a continuous Markov process $Q^{\infty} = (N^{\infty}_1, N^{\infty}_2, W^{\infty})$ s.t. for all $T \geq 0$, $\big\lbrace \hat{Q}^{n}(t)_{0 \leq t \leq T}, n \geq 1 \big\rbrace$ converges weakly to $Q^{\infty}(t)_{0 \leq t \leq T}$ in the space $D^3[0,T]$ under the $J_1$ topology.  Here we note that under the stated weak convergence, $Q^{\infty}(0)$ is distributed as $\bigg( \big(p (1- p) \big)^{\frac{1}{2}} {\mathcal B}^0(1) , - \big(p (1- p) \big)^{\frac{1}{2}} {\mathcal B}^0(1) , 0 \bigg)$.
\end{lemma}
\subsubsection{Review of ergodicity and interchange-of-limits results of \cite{Dieker.13,DDG14}.}
Given a result such as Lemma\ \ref{funweak1}, there is the question of whether an interchange-of-limits hold, i.e. whether the Markov process $Q^{\infty}$ has a unique invariant measure, to which $\lbrace \hat{Q}^{n}(\infty), n \geq 1 \rbrace$ converges.  Such an interchange is not automatic, and there is a large literature studying when such an interchange of limits holds, see the discussion in \cite{DDG14,GS12,Kang.12}.  The sequence of papers \cite{Dieker.13,DDG14} resolved this question for the setting we consider here, i.e. multi-server queues in the HW regime, when services have a phase-type distribution and there is a strictly positive abandonment rate.  First, in \cite{Dieker.13}, the authors proved that the continuous Markov process $Q^{\infty}$ is positive recurrent and has a unique invariant measure $\Upsilon^{\infty}$, and is exponentially ergodic, where we refer the interested reader to \cite{Dieker.13} for precise details regarding the formal definition of exponential ergodicity.  Let $Q^{\infty}(\infty) = \big( N^{\infty}_1(\infty), N^{\infty}_2(\infty), W^{\infty}(\infty) \big)$ denote a generic r.v. distributed as $\Upsilon^{\infty}$.  Then more formally \cite{Dieker.13} proves such a result for a certain ``transformed" Markov process $Y = (Y^1,Y^2),$ where (under analogous initial conditions and on an appropriate probability space) $N^{\infty}_1(t) = Y^1(t) - p \big( Y^1(t) + Y^2(t) \big)^+, 
N^{\infty}_2(t) = Y^2(t) - p \big( Y^1(t) + Y^2(t) \big)^+, W^{\infty}(t) = \big( Y^1(t) + Y^2(t) \big)^+$ for all $t \geq 0$.  Here we refer the interested reader to \cite{Dai.10} for a further discussion of this and related transformations.  We note that the results of \cite{Dieker.13} are easily seen to imply the analogous uniqueness and erodicity results for $Q^{\infty}.$  We also note that although the ergodicity results of \cite{Dieker.13} are stated only when the relevant Markov processes have deterministic initial conditions, when combined with standard results and definitions regarding ergodicity of general Markov processes \cite{Dieker.13,Meyn.92a,Meyn.92b,Meyn.92c,Pes.95}, they are easily seen to imply the analogous results under random initial conditions.
The authors of \cite{Dieker.13} then complete the proof of the desired interchange-of-limits in \cite{DDG14}, where we note that the 
result is actually proven for a certain transformation, but which is easily seen to imply the analogous results for the setting of interest.  These results, including the desired interchange-of-limits, are summarized in the following lemma.
\begin{lemma}[\cite{Dieker.13,DDG14}]\label{posrec1}
The continuous-time Markov process $Q^{\infty}$ is positive recurrent, and has a unique invariant measure $\Upsilon^{\infty}$.  
Furthermore, $\lbrace \hat{Q}^{n}(\infty), n \geq 1 \rbrace$ converges in distribution to $Q^{\infty}(\infty)$, whose associated probability measure is $\Upsilon^{\infty}$.  Namely, an interchange-of-limits holds, which is summarized as follows (after applying the Portmanteau Theorem and basic definitions and results from the theory of weak convergence).  Suppose $Q^{\infty}(0)$ is distributed as $\bigg( \big(p (1- p) \big)^{\frac{1}{2}} {\mathcal B}^0(1) , - \big(p (1- p) \big)^{\frac{1}{2}} {\mathcal B}^0(1) , 0 \bigg)$.  Then for all $x > 0$ which are continuity points of $W^{\infty}(\infty)$ (which will be dense in ${\mathcal R}^+$), 
\begin{equation}\label{interchangeme1}
\liminf_{n \rightarrow \infty} \liminf_{t \rightarrow \infty} P\big( \hat{W}^{n}(t) \geq x \big) \ \ =\ \  \liminf_{n \rightarrow \infty} P\big( \hat{W}^{n}(\infty) \geq x \big) \ \ =\ \ P\big( W^{\infty}(\infty) \geq x \big)\end{equation}
$$=\ \ \lim_{t \rightarrow \infty} P\big( W^{\infty}(t) \geq x \big)\ \ \geq\ \ \limsup_{t \rightarrow \infty} \limsup_{n \rightarrow \infty} P\big( \hat{W}^{n}(t) \geq x \big),$$
and
\begin{equation}\label{interchangeme2}
\limsup_{n \rightarrow \infty} \limsup_{t \rightarrow \infty} P\big( \hat{W}^{n}(t) > x \big)\ \ =\ \ \limsup_{n \rightarrow \infty} P\big( \hat{W}^{n}(\infty) > x \big)\ \ =\ \ P\big( W^{\infty}(\infty) > x \big)
\end{equation}
$$=\ \ \liminf_{t \rightarrow \infty} P\big( W^{\infty}(t) > x \big)\ \ \leq\ \ \liminf_{t \rightarrow \infty} \liminf_{n \rightarrow \infty} P\big( \hat{W}^{n}(t) > x \big).$$
\end{lemma}
We note that the analogous interchange also holds for a general sequence of weakly converging initial conditions (i.e. an analogous interchange-of-limits holds when $Q^{\infty}(0)$ has a different distribution), but for clarity of exposition we do not state at that level of generality.
Also, we note that a study of whether $W^{\infty}(\infty)$ and other related distributions e.g. admit a density (in which case one could avoid the care taken w.r.t open and closed sets above) will be beyond the scope of our paper, and unnecessary for our results.

\subsection{Proof of upper bound part of main result Theorem\ \ref{main1}.}
We now complete the proof of one direction of our main result Theorem\ \ref{main1}.  First, we briefly recall a result from the theory of Gaussian processes, namely the celebrated Borell's inequality, which relates the large deviation properties of the supremum of a centered Gaussian process, which is bounded almost surely, to its maximal variance.  We only state a very special case of the result, as customized to our setting, and refer the interested reader to \cite{Ad90} for further details.

\begin{lemma}[Borell's inequality \cite{Ad90}]\label{borel1}
Suppose that $Z$ is a centered Gaussian process, and that $P\big( \sup_{t \geq 0} Z(t)  < \infty \big) = 1$.  Then $\lim_{x \rightarrow \infty} x^{-2} \log\bigg( P\big( \sup_{t \geq 0} Z(t) > x \big) \bigg) = - \big(2 \sup_{t \geq 0} \E\big[Z^2(t)\big] \big)^{-1}.$
\end{lemma}

With Lemma\ \ref{borel1} in hand, we now complete the proof of the upper bound of our main large deviations result, as formalized in the following lemma.

\begin{lemma}\label{main1up}
Under Assumption\ \ref{A1}, 
$$\limsup_{x \rightarrow \infty} x^{-2} \log P\big( W^{\infty}(\infty) > x \big)  \leq \frac{- \theta( \theta + \mu_1 \mu_2 ) }{2( \theta + \mu_1 + \mu_2 - 1)}.$$
\end{lemma}
\proof{Proof:}
It follows from Lemma\ \ref{posrec1} that for all $x$ which are continuity points of the c.d.f. of $W^{\infty}(\infty)$,
\begin{equation}\label{useinter1}
P\big( W^{\infty}(\infty) > x \big) \leq \liminf_{t \rightarrow \infty} \liminf_{n \rightarrow \infty} P\big( \hat{W}^{n}(t) > x \big).
\end{equation}
By Theorem\ \ref{mainub1}, Lemmas\ \ref{xcon1}, \ref{samedist1}, and \ref{nicesup1}, and the Portmanteau Theorem, for each such continuity point $x$ and $t \geq 0$, 
$$\liminf_{n \rightarrow \infty} P\big( \hat{W}^{n}(t) > x \big) \leq P\bigg( \sup_{0 \leq s \leq t} {\mathcal G}(s) \geq x \bigg).$$
Furthermore, as from definitions and a straightforward calculation w.p.1 $\sup_{t \geq 0}\bigg| {\mathcal G}(t) - \overline{\mathcal G}(t) \bigg| \leq \frac{|B|}{\theta}$, it follows that for all $x$ which are continuity points of the c.d.f. of $W^{\infty}(\infty)$ and all $t \geq 0$,
$$\liminf_{n \rightarrow \infty} P\big( \hat{W}^{n}(t) > x \big) \leq P\bigg( \sup_{0 \leq s \leq t} \overline{\mathcal G}(s) > x - \frac{|B|}{\theta}\bigg),$$
and hence by the monotonicity of the supremum operator that
\begin{equation}\label{useinter2}\liminf_{n \rightarrow \infty} P\big( \hat{W}^{n}(t) > x \big) \leq P\bigg( \sup_{s \geq 0} \overline{\mathcal G}(s) > x - \frac{|B|}{\theta}\bigg).
\end{equation}
Combining (\ref{useinter1}) and (\ref{useinter2}), we conclude that for all such continuity points $x$, 
$$P\big( W^{\infty}(\infty) > x \big) \leq P\bigg( \sup_{t \geq 0} \overline{\mathcal G}(t) > x - \frac{|B|}{\theta} \bigg).$$
Combining with Lemmas\ \ref{entropybound}, \ref{entropybound2}, and \ref{borel1}, it follows that 
$$\limsup_{x \rightarrow \infty} x^{-2} \log\bigg( P\big( W^{\infty}(\infty) > x \big) \bigg) \leq - \big( 2 \sup_{t \geq 0} \E\big[\overline{\mathcal G}^2(t) \big] \big)^{-1}.$$
Combining with Lemma\ \ref{diameter}, which implies that $\sup_{t \geq 0} \E\big[\overline{\mathcal G}^2(t) \big] = \frac{\theta + \mu_1 + \mu_2 - 1}{\theta(\theta + \mu_1 \mu_2)}$, completes the proof.  $\qed$
\endproof

\section{Lower Bound.}\label{LowerSec}
With the upper bound  proven, we now focus on proving lower bounds, to complete the proof of Theorem\ \ref{main1}.
At a high level, the proof of the lower bound proceeds as follows.  First, we even more carefully review the limiting stochastic process $Q^{\infty}$ to which $\lbrace \hat{Q}^n, n \geq 1 \rbrace$ converges as dictated in Lemma\ \ref{funweak1}, explicitly describing the system of SDE to which it is the unique strong solution.  Next, we observe that so long as $W^{\infty}$ has remained strictly positive on $[0,t]$, the dynamics of $Q^{\infty}$ will coincide with those of a much-simpler multi-dimensional OU process, essentially the same as that arising in our previous upper bound.  Finally, we combine these observations with several lower bounds and known results for such multi-dimensional OU processes to complete the proof.

\subsection{Further analysis of the limit process $Q^{\infty}$.}
We begin by establishing some more details about the limiting stochastic process $Q^{\infty}$ to which $\lbrace \hat{Q}^n, n \geq 1 \rbrace$ converges as dictated in Lemma\ \ref{funweak1}.  As mentioned previously, several of the relevant papers describe a certain transformation of this limiting process instead of the limiting process itself.  For completeness, and as we will be using the precise associated system of SDE, we begin by formally defining the relevant transformed limiting process, and more formally clarifying how this limiting process relates to our own.  Although originally studied in \cite{Dai.10} in a more general form (i.e. for phase-type service distributions), a clear exposition of the case of hyper-exponential service times is given in both \cite{Dai.12} and \cite{Dieker.13}, and we now remind the reader of those results.  Let $\lbrace \kappa^i, i \in [1,4] \rbrace$ be four independent standard Brownian motions, constructed on a common probability space.  For a general stochastic process which is denoted using a superscript, e.g. $X^i$, we let $X^{i,+}$ denote the stochastic process s.t. $X^{i,+}(t) = \big( X^i(t) \big)^+$ for all $t \geq 0$.  We now formally define the previously referenced transformed limit processes (i.e. processes which are simple transformations of $Q^{\infty}$), as studied in \cite{Dai.10,Dai.12,Dieker.13}.

\begin{lemma}[\cite{Dai.10,Dai.12,Dieker.13}]\label{lbconverge0}
The following system of SDE (stated for consistency with the relevant literature as a system of integral equations) has a unique strong solution: For all $t \geq 0$,
\begin{eqnarray*}
Y^1(t) &=& Y^1(0) + B p t + p \kappa^3(t) + \big( p(1 - p) \big)^{\frac{1}{2}} \kappa^4(t) - p^{\frac{1}{2}} \kappa^1(t) 
\\&\ &\ \ \ \ \ \ - \mu_1 \int_0^t \bigg( Y^1(s) - p \big( Y^1(s) + Y^2(s) \big)^+ \bigg) ds - p \theta \int_0^t \big( Y^1(s) + Y^2(s) \big)^+ ds;
\end{eqnarray*}
\begin{eqnarray*}
Y^2(t) &=& Y^2(0) + B (1 - p) t + (1 - p) \kappa^3(t) - \big( p(1 - p) \big)^{\frac{1}{2}} \kappa^4(t) - (1 - p)^{\frac{1}{2}} \kappa^2(t) 
\\&\ &\ \ \ \ \ \ - \mu_2 \int_0^t \bigg( Y^2(s) - (1 - p) \big( Y^1(s) + Y^2(s) \big)^+ \bigg) ds - (1 - p) \theta \int_0^t \big( Y^1(s) + Y^2(s) \big)^+ ds.
\end{eqnarray*}
For a given measure $\Upsilon$ on ${\mathcal R}^2$, we let $Y_{\Upsilon} = (Y^1_{\Upsilon}, Y^2_{\Upsilon})$ denote the corresponding solution when $\big( Y^1(0), Y^2(0) \big)$ is distributed as $\Upsilon$, with $Y^1(0), Y^2(0)$ independent of $\lbrace \kappa^i, i \in [1,4] \rbrace$.  In addition, $Y_{\Upsilon}$ is a diffusion process, hence also Markovian with continuous sample paths.
\end{lemma}

Then the following result, proven in \cite{Dai.10}, formalizes the connection between $Y_{\Upsilon}$ and $Q^{\infty}$.  Given a measure $\Upsilon$ on ${\mathcal R}^3$, letting $U = (U_1, U_2, U_3)$ denote a r.v distributed as $\Upsilon$, let $\Upsilon^2$ denote the measure on ${\mathcal R}^2$ corresponding to the distribution of $\big( U_1 + p U_3^+, U_2 + (1 - p) U_3^+ \big)$.  With a slight abuse of notation, given such a measure $\Upsilon$ on ${\mathcal R}^3$, we let $Y_{\Upsilon} = (Y^1_{\Upsilon},Y^2_{\Upsilon})$ denote $Y_{\Upsilon^2} = (Y^1_{\Upsilon^2},Y^2_{\Upsilon^2}).$    Also, for any measure $\Upsilon$ on ${\mathcal R}^3$ s.t. $\Upsilon\bigg( \big\lbrace (x_1,x_2,w) : x_1 + x_2 < 0, w > 0 \big\rbrace \bigg) = \Upsilon\bigg( \big\lbrace (x_1,x_2,w) : x_1 + x_2 > 0 \big\rbrace \bigg) = \Upsilon\bigg( \big\lbrace (x_1,x_2,w) : w < 0 \big\rbrace \bigg) = 0$ , we let $Q^{\infty}_{\Upsilon} = (N^{\infty}_{\Upsilon,1}, N^{\infty}_{\Upsilon,2}, W^{\infty}_{\Upsilon})$ denote the corresponding Markov process when the initial conditions are distributed as $\Upsilon.$  Here we remind the reader that previously we only formally defined $Q^{\infty}$ for the initial conditions distributed as $\bigg( \big(p (1- p) \big)^{\frac{1}{2}} {\mathcal B}^0(1) , - \big(p (1- p) \big)^{\frac{1}{2}} {\mathcal B}^0(1) , 0 \bigg),$ and refer the interested reader to \cite{Dai.10,Dieker.13,DDG14} for further details regarding this Markov process.  Then the connection between $Y_{\Upsilon}$ and $Q^{\infty}_{\Upsilon}$ is as follows, where we note that $Y^i_{\Upsilon}$ intuitively captures the (rescaled and limiting) total number of type-i customers in system (not just in service or waiting in queue).

\begin{theorem}[\cite{Dai.10}]\label{lbconverge1}
For any measure $\Upsilon$ on ${\mathcal R}^3$ s.t. $\Upsilon\bigg( \big\lbrace (x_1,x_2,w) : x_1 + x_2 < 0, w > 0 \big\rbrace \bigg) = \Upsilon\bigg( \big\lbrace (x_1,x_2,w) : x_1 + x_2 > 0 \big\rbrace \bigg) = \Upsilon\bigg( \big\lbrace (x_1,x_2,w) : w < 0 \big\rbrace \bigg) = 0$, one may construct $Y_{\Upsilon}$ and $Q^{\infty}_{\Upsilon}$ on a common probability space s.t. w.p.1, for all $t \geq 0$, $N^{\infty}_{\Upsilon, 1}(t) = Y^1_{\Upsilon}(t) - p \big( Y^1_{\Upsilon}(t) + Y^2_{\Upsilon}(t) \big)^+; N^{\infty}_{\Upsilon, 2}(t) = Y^2_{\Upsilon}(t) - (1- p) \big( Y^1_{\Upsilon}(t) + Y^2_{\Upsilon}(t) \big)^+; W^{\infty}_{\Upsilon}(t) = \big( Y^1_{\Upsilon}(t) + Y^2_{\Upsilon}(t) \big)^+.$
\end{theorem}

In our analysis, it will be convenient to work not with $Y_{\Upsilon}$, but instead with a relevant transformation.  Under this transformation, the first component roughly captures the infinitesimal drift due to service completions, with the second component capturing the rescaled total number in system, where this transformation also has the benefit of resulting in certain transformed terms becoming independent of one-another.  In particular, for a given measure $\Upsilon$ on ${\mathcal R}^2$, let $\overline{Y}_{\Upsilon} = (\overline{Y}^1_{\Upsilon}, \overline{Y}^2_{\Upsilon})$ denote the stochastic process s.t. for all $t \geq 0$, $\overline{Y}^1_{\Upsilon}(t) = (1- p) Y^1_{\Upsilon}(t) - p Y^2_{\Upsilon}(t)$, and $\overline{Y}^2_{\Upsilon}(t) = Y^1_{\Upsilon}(t) + Y^2_{\Upsilon}(t)$.  We note that for such a measure $\Upsilon$ on ${\mathcal R}^2$, $\overline{Y}_{\Upsilon}(0)$ is not distributed as $\Upsilon$, but as $\big( (1 - p) U_1 - p U_2 , U_1 + U_2 \big)$ for a r.v. $(U_1, U_2)$ distributed as $\Upsilon.$  As before, with a slight abuse of notation, given a measure $\Upsilon$ on ${\mathcal R}^3$, we let $\overline{Y}_{\Upsilon}$ denote the analogously defined stochastic process (defined in terms of $Y_{\Upsilon}$ with $Y_{\Upsilon}(0)$ distributed as $\Upsilon^2$).

Then we may conclude the following from Lemma\ \ref{lbconverge0}, by taking appropriate linear combinations of the integral equations appearing in Lemma\ \ref{lbconverge0}, and some straightforward algebra, and we omit the details.

\begin{corollary}--\label{lbconverge2}
The following system of SDE has a unique strong solution: For all $t \geq 0$,
\begin{eqnarray*}\overline{Y}^1(t) &=& \overline{Y}^1(0) + \big( p (1 - p) \big)^{\frac{1}{2}} \kappa^4(t) + \bigg( p (1 - p)^{\frac{1}{2}} \kappa^2(t) - (1 - p) p^{\frac{1}{2}} \kappa^1(t) \bigg) 
\\&\ &\ \ \ \ \ \ - \mu_1 \mu_2 \int_0^t \overline{Y}^1(s) ds + p (1 - p) \mu_1 \mu_2 \int_0^t \big( \overline{Y}^{2,+}(s) - \overline{Y}^2(s) \big) ds;
\end{eqnarray*}
\begin{eqnarray*}
\overline{Y}^2(t) &=& \overline{Y}^2(0) + Bt + \kappa^3(t) - \big( p^{\frac{1}{2}} \kappa^1(t) + (1 - p)^{\frac{1}{2}} \kappa^2(t) \big) - (\mu_1 - \mu_2) \int_0^t \overline{Y}^1(s) ds 
\\&\ &\ \ \ \ \ \ - \theta \int_0^t \overline{Y}^2(s) ds + \big( p \mu_1 + (1 - p) \mu_2 - \theta) \int_0^t \big( \overline{Y}^{2,+}(s) - \overline{Y}^2(s) \big) ds.
\end{eqnarray*}
For any given measure $\Upsilon$ on ${\mathcal R}^2$, if $\big( \overline{Y}^1(0), \overline{Y}^2(0) \big)$ is distributed as 
$\big( (1 - p) U_1 - p U_2 , U_1 + U_2 \big)$ for a r.v. $(U_1, U_2)$ distributed as $\Upsilon$, with $\overline{Y}^1(0), \overline{Y}^2(0)$ independent of $\lbrace \kappa^i, i \in [1,4] \rbrace$, then this unique strong solution has the same distribution (at the process level) as $\overline{Y}_{\Upsilon}$.
\end{corollary}

We now make an additional simplification.  In particular, we note that by a simple covariance calculation, it is easily verified that the Brownian motion $p (1 - p)^{\frac{1}{2}} \kappa^2 - (1 - p) p^{\frac{1}{2}} \kappa^1$ is independent from the Brownian motion 
$p^{\frac{1}{2}} \kappa^1 + (1 - p)^{\frac{1}{2}} \kappa^2$, essentially for the same reason that the sum and difference of two independent Brownian motions are independent Brownian motions.  In that case, it follows from some straightforward variance calculations that Corollary\ \ref{lbconverge2} implies the following.  Recall that ${\mathcal B}^3$ and ${\mathcal B}^2$ are two independent Brownian motions which satisfy $Var[{\mathcal B}^3(1)] = 2 p (1 - p), Var[{\mathcal B}^2(1)] = 2$.

\begin{corollary}\label{lbconverge3}
The following system of SDE has a unique strong solution: For all $t \geq 0$,
$$\overline{Y}^1(t)= \overline{Y}^1(0) - \mu_1 \mu_2 \int_0^t \overline{Y}^1(s) ds + {\mathcal B}^3(t) + p (1 - p) \mu_1 \mu_2 \int_0^t \big( \overline{Y}^{2,+}(s) - \overline{Y}^2(s) \big) ds;$$
\begin{eqnarray*}
\overline{Y}^2(t) &=& \overline{Y}^2(0) + Bt - (\mu_1 - \mu_2) \int_0^t \overline{Y}^1(s) ds - \theta \int_0^t \overline{Y}^2(s) ds + {\mathcal B}^2(t) 
\\&\ &\ \ \ \ \ \ + \big( p \mu_1 + (1 - p) \mu_2 - \theta) \int_0^t \big( \overline{Y}^{2,+}(s) - \overline{Y}^2(s) \big) ds.
\end{eqnarray*}
For any given measure $\Upsilon$ on ${\mathcal R}^2$, if $\big( \overline{Y}^1(0), \overline{Y}^2(0) \big)$ is distributed as 
$\big( (1 - p) U_1 - p U_2 , U_1 + U_2 \big)$ for a r.v. $(U_1, U_2)$ distributed as $\Upsilon$, with $\overline{Y}^1(0), \overline{Y}^2(0)$ independent of ${\mathcal B}^3, {\mathcal B}^2$, then this unique strong solution has the same distribution (at the process level) as $\overline{Y}_{\Upsilon}$.
\end{corollary}

\subsection{Relating $Q^{\infty}$ to a multi-dimensional OU process.}
At this point, it will be helpful to remind ourselves of an SDE we have already encountered in our upper bound arguments, namely that encountered in Lemma\ \ref{thesol1lemma}, whose solution is the multi-dimensional OU process.  In particular, the following is a restatement of Lemma\ \ref{thesol1lemma}, written as an equivalent system of stochastic integral equations (as opposed to differential equations), after applying Lemma\ \ref{Qexp1}, where here the results are used to analyze an appropriate relaxation of $\overline{Y}$ as opposed to $\OU_{+,s}$ and $\zeta\big( \Lambda, s \big),$ and we omit the details (the logic being nearly identical to that already used in our previous analysis).  For a given measure $\Upsilon$ on ${\mathcal R}^2$, let $\hat{Y}_{\Upsilon} = (\hat{Y}^1_{\Upsilon}, \hat{Y}^2_{\Upsilon})$ denote the stochastic process explicitly constructed on the same probability space as ${\mathcal B}^3, {\mathcal B}^2$ s.t. w.p.1, for all $t \geq 0$,
$$\hat{Y}^1_{\Upsilon}(t) = \exp(- \mu_1 \mu_2 t) \hat{Y}^1_{\Upsilon}(0) + \int_0^t \exp\big( - \mu_1 \mu_2 (t - y) d{\mathcal B}^3(y);$$
\begin{eqnarray*}
\hat{Y}^2_{\Upsilon}(t) &=& \frac{\mu_1 - \mu_2}{\mu_1 \mu_2 - \theta} \big( \exp(- \mu_1 \mu_2 t) - \exp(- \theta t) \big) \hat{Y}^1_{\Upsilon}(0) + \exp(- \theta t) \hat{Y}^2_{\Upsilon}(0) 
\\&\ &\ \ \ \ \ \ + \frac{\mu_1 - \mu_2}{\mu_1 \mu_2 - \theta} \int_0^t \bigg( \exp\big(- \mu_1 \mu_2 (t - y) \big) - \exp\big(- \theta (t - y) \big) \bigg) d{\mathcal B}^3(y) 
\\&\ &\ \ \ \ \ \ + \int_0^t \exp\big( - \theta (t - y) \big) d{\mathcal B}^2(y) + \frac{B}{\theta}\bigg(1 - \exp\big(- \theta t \big) \bigg),
\end{eqnarray*}
with $\big( \hat{Y}^1_{\Upsilon}(0), \hat{Y}^2_{\Upsilon}(0) \big)$ distributed as $\Upsilon$, independent of ${\mathcal B}^3, {\mathcal B}^2$.

\begin{corollary}--\label{lbconverge4}
The following system of SDE has a unique strong solution: For all $t \geq 0$,
$$\hat{Y}^1(t)= \hat{Y}^1(0) - \mu_1 \mu_2 \int_0^t \hat{Y}^1(s) ds + {\mathcal B}^3(t);$$
$$\hat{Y}^2(t) = \hat{Y}^2(0) + Bt - (\mu_1 - \mu_2) \int_0^t \hat{Y}^1(s) ds - \theta \int_0^t \hat{Y}^2(s) ds + {\mathcal B}^2(t).$$
For any given measure $\Upsilon$ on ${\mathcal R}^2$, if $\big( \hat{Y}^1(0), \hat{Y}^2(0) \big)$ is distributed as $\Upsilon$, with $\hat{Y}^1(0), \hat{Y}^2(0)$ independent of ${\mathcal B}^3, {\mathcal B}^2$, this unique strong solution has the same distribution (at the process level) as 
$\big( \hat{Y}^1_{\Upsilon}, \hat{Y}^2_{\Upsilon} \big)$.
\end{corollary}

Now, we make the simple observation that the system of SDE from Corollary\ \ref{lbconverge3} and that of Corollary\ \ref{lbconverge4} agree on any interval $[0,T]$ s.t. $\overline{Y}^{2,+}(t) = \overline{Y}^2(t)$ for all $t \in [0,T]$.  Combining with Theorem\ \ref{lbconverge1} and Corollaries\ \ref{lbconverge2} - \ref{lbconverge4}, and the fact that for any measure $\Upsilon$ on ${\mathcal R}^3$ s.t. $\Upsilon\bigg( \big\lbrace (x_1,x_2,w) : x_1 + x_2 < 0, w > 0 \big\rbrace \bigg) = \Upsilon\bigg( \big\lbrace (x_1,x_2,w) : x_1 + x_2 > 0 \big\rbrace \bigg) = \Upsilon\bigg( \big\lbrace (x_1,x_2,w) : w < 0 \big\rbrace \bigg) = 0$ it holds that 
$\big( Y^1_{\Upsilon} + Y^2_{\Upsilon} \big)^+, \overline{Y}^{2,+}_{\Upsilon}$, and $W^{\infty}_{\Upsilon}$ all have the same distribution (at the process level), we conclude the following Corollary\ \ref{sameasOU}, which relates the probability that $W^{\infty}_{\Upsilon}(t)$ exceeds a level $x$, conditioned on $W^{\infty}_{\Upsilon}$ having remained strictly positive on $[0,t]$, to the analogous probability for the second component of the much simpler multi-dimensional OU process similarly conditioned.  Given a measure $\Upsilon$ on ${\mathcal R}^3,$ letting $U = (U_1, U_2, U_3)$ denote a r.v distributed as $\Upsilon$,  let ${\Upsilon}^{2'}$ denote the measure on ${\mathcal R}^2$ corresponding to the distribution of $\big( (1 - p) U_1 - p U_2, U_1 + U_2 + U^+_3 \big)$.  With a slight abuse of notation, given such a measure $\Upsilon$ on ${\mathcal R}^3$, we let $\hat{Y}_{\Upsilon} = (\hat{Y}^1_{\Upsilon}, \hat{Y}^2_{\Upsilon})$ denote $\hat{Y}_{\Upsilon^{2'}} = (\hat{Y}^1_{\Upsilon^{2'}},\hat{Y}^2_{\Upsilon^{2'}})$.  Given a measure $\Upsilon$ on ${\mathcal R}^3$ s.t. $\Upsilon\bigg( \big\lbrace (x_1,x_2,w) : x_1 + x_2 < 0, w > 0 \big\rbrace \bigg) = \Upsilon\bigg( \big\lbrace (x_1,x_2,w) : x_1 + x_2 > 0 \big\rbrace \bigg) = \Upsilon\bigg( \big\lbrace (x_1,x_2,w) : w < 0 \big\rbrace \bigg) = 0,$ let $\tau^{Q^{\infty}}_{\Upsilon} \stackrel{\Delta}{=} \inf\lbrace t \geq 0 : W^{\infty}_{\Upsilon}(t) \leq 0 \rbrace$, and $\tau^{\hat{Y}}_{\Upsilon} \stackrel{\Delta}{=} \inf\lbrace t \geq 0 : \hat{Y}^2_{\Upsilon}(t) \leq 0 \rbrace$.

\begin{corollary}--\label{sameasOU}
For all measures $\Upsilon$ on ${\mathcal R}^3$ s.t. $\Upsilon\bigg( \big\lbrace (x_1,x_2,w) : x_1 + x_2 < 0, w > 0 \big\rbrace \bigg) = \Upsilon\bigg( \big\lbrace (x_1,x_2,w) : x_1 + x_2 > 0 \big\rbrace \bigg) = \Upsilon\bigg( \big\lbrace (x_1,x_2,w) : w < 0 \big\rbrace \bigg) = 0$, and all $t,x > 0$, 
$$\E\big[ I\big( W^{\infty}_{\Upsilon}(t) > x , \tau^{Q^{\infty}}_{\Upsilon} > t \big) \big] = \E\big[ I\big( \hat{Y}^2_{\Upsilon}(t) > x , \tau^{\hat{Y}}_{\Upsilon} > t \big) \big].$$
\end{corollary}

To proceed, we will reason about applying Corollary\ \ref{sameasOU} when $\Upsilon$ is the unique invariant measure of $Q^{\infty}$, i.e. the previously defined measure $\Upsilon^{\infty}$ governing the distribution of $Q^{\infty}(\infty)$, whose existence was established in Lemma\ \ref{posrec1}.  We note that it is easily seen to follow from the results of \cite{Dieker.13} and the definition of $\overline{Y}_{\Upsilon^{\infty}}$ that $\lbrace \overline{Y}_{\Upsilon^{\infty}}(t), t \geq 0 \rbrace$ is a stationary process.  We also note that $\hat{Y}_{\Upsilon^{\infty}}$ is equivalent to a multi-dimensional OU process initialized not with its own stationary measure, but with that of a different process.
\\\indent First, we will need to prove that $\Upsilon^{\infty}$ is sufficiently non-degenerate, and defer the proof to the appendix, noting that a very similar infinite-server lower bound is used in \cite{GS12} towards different ends.

\begin{lemma}*\label{notbad1}
Under Assumption\ \ref{A1}, there exist $C_1, C_2 \in ( 4 \frac{|B|}{\theta} , \infty)$ s.t. $C_1 < C_2$; $C_3, C_4 \in (-\infty,\infty)$ s.t. $C_3 < C_4$; and $\delta \in (0,1)$, all depending only on $p,\mu_1,\mu_2,B,\theta$, s.t. 
$$P\big( \hat{Y}^1_{\Upsilon^{\infty}}(0) \in [C_3, C_4] ,  \hat{Y}^2_{\Upsilon^{\infty}}(0) \in [C_1, C_2] \big) \geq \delta.$$
\end{lemma}

Let $\Upsilon^*$ denote the probability measure governing the distribution of $\big( \hat{Y}^1_{\Upsilon^{\infty}}(0), \hat{Y}^2_{\Upsilon^{\infty}}(0) \big)$ when conditioned on the event $\big\lbrace \hat{Y}^1_{\Upsilon^{\infty}}(0) \in [C_3, C_4] ,  \hat{Y}^2_{\Upsilon^{\infty}}(0) \in [C_1, C_2] \big\rbrace$, for $C_1,C_2,C_3,C_4$ as provided by Lemma\ \ref{notbad1}.  Then by combining Lemma\ \ref{posrec1} and Theorem\ \ref{lbconverge1} with Corollaries\ \ref{lbconverge2} - \ref{sameasOU} and Lemma\ \ref{notbad1}, we may relate $P\big( W^{\infty}(\infty) > x \big)$ to the probability that the second component of a much simpler multi-dimensional OU process exceeds x, with the additional requirement that this second component does not hit 0 over an appropriate time interval.
\begin{lemma}\label{chainit1}
Let $\delta \in (0,1)$ be as in Lemma\ \ref{notbad1}.  Then for all $x,t \geq 0$, 
$$P\big( W^{\infty}(\infty) > x \big)\ \ \geq\ \ \delta \times \E\bigg[ I\bigg( \hat{Y}^2_{\Upsilon^*}(t) > x , \tau^{\hat{Y}}_{\Upsilon^*} > t \bigg) \bigg].$$
\end{lemma}
\proof{Proof: }
By combining Lemma\ \ref{posrec1} and  Theorem\ \ref{lbconverge1} with Corollaries\ \ref{lbconverge2} - \ref{sameasOU} and Lemma\ \ref{notbad1}, we conclude that for all $x > 0$,
\begin{eqnarray*}
P\big( W^{\infty}(\infty) > x \big) &=&  P\big( W^{\infty}_{\Upsilon^{\infty}}(t) > x \big)
\\&=& \E\bigg[ I\bigg( \overline{Y}^2_{\Upsilon^{\infty}}(t) > x \bigg) \bigg]
\\&\geq& \E\bigg[ I\bigg( \overline{Y}^2_{\Upsilon^{\infty}}(t) > x , \overline{Y}^1_{\Upsilon^{\infty}}(0) \in [C_3,C_4], \overline{Y}^2_{\Upsilon^{\infty}}(0) \in [C_1,C_2] \bigg) \bigg]
\\&\geq& \delta \E\big[ I\big( W^{\infty}_{\Upsilon^*}(t) > x \big) \big]
\\&\geq& \delta \E\big[ I\big( W^{\infty}_{\Upsilon^*}(t) > x , \tau^{Q^{\infty}}_{\Upsilon^*} > t \big) \big]\ \ \ =\ \ \ \delta \E\bigg[ I\bigg( \hat{Y}^2_{\Upsilon^*}(t) > x , \tau^{\hat{Y}}_{\Upsilon^*} > t \bigg) \bigg],
\end{eqnarray*}
completing the proof.  $\qed$
\endproof

\subsection{Additional bounds and asymptotics for multi-dimensional OU processes.}
Before proceeding, it will be helpful to apply arguments similar to those used in the proof of Lemma\ \ref{itorewrite1} to derive a simplified representation for (a process which lower bounds) $\hat{Y}^2_{\Upsilon^*}$.  Intuitively, Lemma\ \ref{notbad1}, especially the fact that $C_1 > \frac{4 |B|}{\theta}$, along with the special structure of the relevant multi-dimensional OU processes, will allow us to rewrite the relevant stochastic processes in a form amenable to lower-bounding the hitting times of interest by several relatively simple events occuring simultaneously, e.g. the hitting times to levels involving $B,C_1,C_2,C_3,C_4$ of certain Brownian motions.
\\\indent We begin with some additional definitions.  Let $\OU_{0,a}$ denote an explicit construction of the unique w.p.1 continuous centered Gassian process s.t. $\E[ \OU_{0,a}(s) \OU_{0,a}(t)] = \frac{p (1 - p)}{\mu_1 \mu_2}  \bigg( \exp\big( - \mu_1 \mu_2 (t - s) \big) - \exp\big( - \mu_1 \mu_2 (t + s) \big) \bigg)$, equivalently the appropriate non-stationary OU process initialized at 0.  Similarly, let $\OU_{0,b}$ denote an explicit construction of the unique w.p.1 continuous centered Gassian process s.t. $\E[ \OU_{0,b}(s) \OU_{0,b}(t)] = \theta^{-1} \bigg( \exp\big( - \theta (t - s) \big) - \exp\big( - \theta (t + s) \big) \bigg),$ again an appropriate non-stationary OU process initialized at 0.
Further suppose that $\OU_{0,a}$ and $\OU_{0,b}$ are constructed on the same probability space, and are independent of one-another.  Let us define $\Pi$ to be the stochastic process s.t. for all $t \geq 0$,
$$\Pi(t) = ( \mu_1 - \mu_2 ) \exp(- \theta t) \int_0^t \exp( \theta y ) \OU_{0,a}(y) dy + \OU_{0,b}(t).$$
Also, for $C_1,C_2,C_3,C_4$ as in Lemma\ \ref{notbad1}, let us define
$$\BUF(t) \stackrel{\Delta}{=} \exp( - \theta t ) C_1 - (\mu_1 + \mu_2) \big( |C_3| + |C_4| \big) t \exp(- \theta t) - \frac{|B|}{\theta}.$$
Then the desired representation is as follows, where we defer all proofs to the appendix.
\begin{lemma}*\label{itorewrite2}
One can construct $\Pi$ and $\hat{Y}^2_{\Upsilon^*}$ on a common probability space s.t. w.p.1, for all $t \geq 0$,
$$
\hat{Y}^2_{\Upsilon^*}(t) \geq \Pi(t) + \BUF(t).$$
\end{lemma}
Note that $\BUF$ is strictly positive on some interval containing zero.  This will allow us to bound the probability of certain events of interest by the probability that certain processes which start at 0 don't go very negative over a small interval containing zero, and have already raced up to a large value by the end of said interval.
\\\indent Next, it will be helpful to prove some structural properties of the covariance of $\Pi$, where we again defer all proofs to the appendix.  

\begin{lemma}*\label{covlemma1lower}
$\E\big[ \Pi(s) \Pi(t) \big] \geq 0$ for $0 \leq s \leq t$, and $\lim_{t \rightarrow \infty} \E\big[ \Pi^2(t) \big] = \frac{\theta + \mu_1 + \mu_2 - 1}{\theta (\theta + \mu_1 \mu_2)}.$
\end{lemma}

\subsection{Proof of lower bound part of main result Theorem\ \ref{main1}.}
We now proceed as follows to complete the proof of the lower bound of our main result Theorem\ \ref{main1}.  Lemmas\ \ref{chainit1} and\ \ref{itorewrite2}, suggest to fix a time $t$, and analyze the probability that the explicit process $\Pi$ stays above a certain explicit boundary, related to three things: 1. being above x at  time t, 2. requiring that $\tau^{\hat{Y}}_{\Upsilon^*} > t$, and 3. compensating for the shift involving $\BUF$ appearing in Lemma\ \ref{itorewrite2}.  As the asymptotic variance of $\Pi$, $\frac{\theta + \mu_1 + \mu_2 - 1}{\theta (\theta + \mu_1 \mu_2)}$, matches that appearing in our upper bound, and large deviations are not affected by any fixed shift, the main hurdle will be to show that requiring $\tau^{\hat{Y}}_{\Upsilon^*} > t$ does not somehow change the large deviations behavior of the process at time $t$.  
\subsubsection{Slepian's inequality and bounds for a finite-horizon infimum.}
Fortunately, the fact that Lemma\ \ref{covlemma1lower} proves that $\Pi$ has non-negative covariance will allow us to use the celebrated Slepian's inequaliy to this end.  Let us now recall this well-known result from the theory of Gaussian processes, stating a similar variant to that used in \cite{G16} (where we note that this variant is typically referred to as Gordon's inequality).
\begin{lemma}[Slepian's/Gordon's inequality \cite{Ad90}]\label{Gordon1}
Let $\lbrace X_i, i = 1,\ldots,m \rbrace$ and $\lbrace Y_i, i = 1,\ldots,m \rbrace$ be two centered $m$-dimensional Gaussian vectors s.t. $\E[X^2_i] = \E[Y^2_i]$ for all $i$, and $\E[X_i X_j] \leq \E[Y_i Y_j]$ for all $i,j \in \lbrace 1,\ldots,m \rbrace$.  Then for all $z \in {\mathcal R}^m$, $P\big( \bigcap_{i=1,\ldots,m} \lbrace Y_i > z_i \rbrace \big) \geq P\big( \bigcap_{i=1,\ldots,m} \lbrace X_i > z_i \rbrace \big)$.
\end{lemma}
In light of the fact that $\Pi$ has non-negative covariance by Lemma\ \ref{covlemma1lower}, by fixing $T \in {\mathcal R}^+$ and considering the centered Gaussian process $\Pi'$ s.t. for all $0 \leq s \leq t < T$, 
$\E[\Pi'(s) \Pi'(t)] = \E[\Pi(s) \Pi(t)\big],$ while for all $0 \leq s < T$, 
$\E[\Pi'(s) \Pi'(T)] = 0,$ and for which 
$\E[\Pi'^2(T)] = \E[\Pi^2(T)],$ i.e. the centered Gaussian process whose covariance structure on $[0,T]$ is the same as that of $\Pi$ but with the process at time $T$ independent of the process on $[0,T)$, Lemma\ \ref{Gordon1} is easily seen (after combining the continuity of $\Pi$ with a straightforward approximation argument, the details of which we omit) to imply the following.
\begin{corollary}--\label{Gordon2}
For all $x \in {\mathcal R}$ and $T \geq 0$, 
$$\E\big[ I\bigg( \inf_{t \in [0,T]} \big( \Pi(t) + \BUF(t) \big) > 0, \Pi(T) + \BUF(T) > x \bigg) \big]$$
 is at least 
$$\E\big[ I\bigg( \inf_{t \in [0,T]} \big( \Pi(t) + \BUF(t) \big) > 0 \bigg) \big] \times \E\big[I\bigg( \Pi(T) > x - \BUF(T) \bigg) \big].
$$
\end{corollary}
In light of Lemma\ \ref{itorewrite2}, Corollary\ \ref{Gordon2} essentially accomplishes the goal of allowing us to consider the large deviations at time $T$ and the probability of remaining positive on $[0,T]$ separately.  We now focus on analyzing the two terms appearing in Corollary\ \ref{Gordon2}.  Our general approach to completing the proof of our main result will be to take a double limit, fixing a large value for the time $T$ considered and then for each such large value computing the large deviations behavior (i.e. letting $x \rightarrow \infty$ for $x$ as in Corollary\ \ref{Gordon2}).  We note that since $\E\big[ I\bigg( \inf_{t \in [0,T]} \big( \Pi(t) + \BUF(t) \big) > 0 \bigg) \big]$ does not depend on the choice of $x$ (again as in Corollary\ \ref{Gordon2}), and we are only interested in the relevant large deviations behavior, it will suffice to show that for each fixed $T$ this quantity is strictly positive.  Then the desired result is the following lemma, whose proof we defer to the appendix.
\begin{lemma}*\label{isposlower1}
For all $T > 0$, $\E\big[ I\bigg( \inf_{t \in [0,T]} \big( \Pi(t) + \BUF(t) \big) > 0 \bigg) \big] > 0$.
\end{lemma}
\subsubsection{Proof of lower bound part of main result Theorem\ \ref{main1}.}
Finally, we are in a position to complete the proof of our main result Theorem\ \ref{main1}.  
\proof{Proof of Theorem\ \ref{main1}:}
In light of the already proven upper bound from Lemma\ \ref{main1up}, it suffices to prove the corresponding lower bound.
By Lemma\ \ref{chainit1} and Lemma\ \ref{itorewrite2}, there exists $\delta \in (0,1)$ s.t. for all $T,x > 0$, $P\big( W^{\infty}(\infty) > x \big)$ is at least
\begin{equation}\label{finale1}
\delta \times \E\bigg[ I\bigg( \Pi(t) + \BUF(t) > x , \inf_{t \in [0,T]} \big( \Pi(t) + \BUF(t) \big) > 0 \bigg) \bigg],
\end{equation}
which by Corollary\ \ref{Gordon2} is at least 
\begin{equation}\label{finale2}
\delta \times \E\big[ I\bigg( \inf_{t \in [0,T]} \big( \Pi(t) + \BUF(t) \big) > 0 \bigg) \big] \times \E\big[I\bigg( \Pi(T) + \BUF(T) > x \bigg) \big].
\end{equation}
It follows from Lemma\ \ref{covlemma1lower} that for all $\epsilon > 0$, there exists $T_{\epsilon}$ s.t. 
\begin{equation}\label{finale3}
\E\big[ \Pi^2(T_{\epsilon}) \big] > (1 - \epsilon) \frac{\theta + \mu_1 + \mu_2 - 1}{\theta (\theta + \mu_1 \mu_2)}.
\end{equation}
Combining (\ref{finale1}) - (\ref{finale3}), we conclude that for all $\epsilon > 0$, $\liminf_{x \rightarrow \infty} x^{-2} \log\bigg( P\big( W^{\infty}(\infty) > x \big) \bigg)$ is at least
\begin{eqnarray*}
\ &\ &\ \liminf_{x \rightarrow \infty} x^{-2} \log\Bigg( \delta \times \E\big[ I\bigg( \inf_{t \in [0,T_{\epsilon}]} \big( \Pi(t) + \BUF(t) \big) > 0 \bigg) \big] 
\\&\ &\ \ \ \ \ \ \ \ \ \ \ \ \ \ \ \ \ \ \ \ \ \ \ \ \times\ \ \ \E\big[I\bigg( \Pi(T_{\epsilon}) + \BUF(T_{\epsilon}) > x \bigg) \big] \Bigg),
\end{eqnarray*}
itself equal to
\begin{eqnarray*}
\ &\ &\ \liminf_{x \rightarrow \infty} x^{-2} \log\Bigg( \E\big[I\bigg( \Pi(T_{\epsilon}) + \BUF(T_{\epsilon}) > x \bigg) \big] \Bigg)
\\&\ &\ \ \ \ =\ \ \ \liminf_{x \rightarrow \infty} x^{-2} \log\bigg( P\big( \Pi(T_{\epsilon})  > x \big) \bigg),
\end{eqnarray*}
which by the basic properties of the normal distribution (see \cite{Ad90}) equals  $- \big( 2 \E[ \Pi^2(T_{\epsilon}) \big)^{-1}$, which by (\ref{finale3}) is at least $-\frac{1}{1 - \epsilon} \frac{\theta (\theta + \mu_1 \mu_2)}{2(\theta + \mu_1 + \mu_2 - 1)}$.
Taking the limit $\epsilon \downarrow 0$ then completes the proof of the desired large deviations result.  We defer the proof that the proven exponent does not depend only on the first two moments of $S$, as well as the fact that $- \frac{\theta( \theta + \mu_1 \mu_2 ) }{2( \theta + \mu_1 + \mu_2 - 1)} < - \frac{\mu_1 \mu_2 \theta}{2 (\mu_1 + \mu_2 - 1)},$ to the appendix.  $\qed$
\endproof

\section{Conclusion.}\label{ConcSec}
In this paper, we provided the first precise large-deviations analysis of steady-state multi-server queues with abandonments in the Halfin-Whitt regime when service times are non-Markovian.  In particular, for the setting in which inter-arrival times are Markovian (with spare capacity parameter $B$); service times are hyper-exponentially distributed with mean 1 and parameters $\mu_1,\mu_2,p$; and patience times are Markovian with rate $\theta$ satisfying $0 < \theta < \min(\mu_1,\mu_2)$ (i.e. Assumption\ \ref{A1}), we proved that the weak limit $W^{\infty}$ of the associated sequence of normalized steady-state queue lengths satisfies $\lim_{x \rightarrow \infty} x^{-2} \log \bigg( P\big( W^{\infty}(\infty) > x \big) \bigg)  = - \frac{\theta( \theta + \mu_1 \mu_2 ) }{2( \theta + \mu_1 + \mu_2 - 1)}.$  This is the first result providing such a fine-grained understanding of the behavior of multi-server queues with abandonments in the Halfin-Whitt regime, when service times are non-Markovian.  Furthermore, our results disproved a conjecture of Dai and He, which had posited a particular value for this exponent which we ultimately proved to be incorrect.  Our results also show that the related conjecture of Dai and He that this exponent should depend only on the first two moments of the service distribution, i.e. should exhibit an insensitivity, was similarly incorrect.  This negative resolution shows that in systems with abandonments the large-deviations behavior is considerably more complex than in the setting without abandonments, as the conjectures of Dai an He were in analogy with known results for queues without abandonments, where the true answer is much simpler and indeed insensitive beyond the first two moments.  Our results were also able to provide additional insight into the numerical methods of Dai and He which had inspired their conjecture, as well as how the probability of seeing a large queue length changes as e.g. the abandonment rate goes to 0, or the larger of the two service rates goes to $\infty$.  Perhaps most importantly, our results shed light on a connection between multi-server queues with abandonments in the Halfin-Whitt regime and certain tractable diffusion processes called multi-dimensional Ornstein-Uhlenbeck processes, which can be viewed as relaxations of the dynamics of the true complex limiting processes.  Along the way we also proved several results which may be of independent interest, including a stochastic comparison result for general $M/H_2/n + M$ queues (assuming $\theta < \min(\mu_1,\mu_2)$) which bounds the true system by one with simpler dynamics, in analogy to the stochastic comparison results of \cite{GG13}.  We note that such a comparison may prove useful e.g. in the simulation of multi-server queues with abandonments, in analogy to the utility that the results of \cite{GG13} had for devising simulation algorithms for multi-server queues without abandonments \cite{Blan.12}.  
\\\indent Our work leaves several interesting directions for future research.  The most natural two are: 1. extending our results to $M/PH/n + M$ queues in the Halfin-Whitt regime under a hazard-rate assumption analogous to our Assumption\ \ref{A1}, and 2. extending our results beyond the confines of Assumption\ \ref{A1}.  Although for space considerations we do not here detail precisely how one might acheive 1., i.e. extending our results to the case of phase-type services, we note that for many steps of our approach it is at least clear what one would like to prove (for phase-type services in analogy with our analysis for $H_2$ services), and (roughly) how one might go about it.  As far as 2., a natural first approach would be to attempt to prove stochastic-comparison results directly for the diffusion limit.  Indeed, a careful examination of the relevant diffusion limit $\overline{Y}$ reveals certain (anti-)monotonicities between various terms and their drifts, around which one may be able to build an appropriate coupling.  We note that in the limit certain state-space collapse results have already kicked in \cite{Dai.10}, which may make constructing such couplings an easier task.  Generalizing our pre-limit stochastic-comparison upper bound would seem to be a potentially harder task.  Indeed, it is far less clear whether such a result even holds in the pre-limit beyond the confines of Assumption\ \ref{A1}.  It would also be interesting to develop a deeper understanding of our pre-limit upper bound, whose distribution is framed as the solution to an appropriate Skorokhod problem, and which we only analyze in-depth in the limit.  We note that for more general service distributions, the analogy of Assumption\ \ref{A1} would be to require that the supremum of the hazard rate of the abandonment distribution is at most the infimum of the hazard rate of the service distribution.  On that note, we also believe that our approach should generalize to different inter-arrival and abandonment distributions, and perhaps even to multi-class systems supposing that a state-space-collapse analogous to that which holds in the single-class setting also occurs in the multi-class setting.
\\\indent Another interesting set of questions surround the connections between the true limit processes and the more tractable multi-dimensional OU processes brought to light by our work.  How general is this phenomena, i.e. that taking a complex and non-linear system of SDE and applying an appropriate ``relaxation" (linearization in our case) yields a simpler process which provably shares many of the properties of the true process?  This may relate to questions surrounding precisely how different types of reflections change the qualitative properties of a stochastic process.  Further understanding when the large-deviations exponent of the true process and its relaxation agree might also yield additional insights along these lines, where we note that the framework of \cite{Kang.14} may be especially relevant.  It would also be of interest to investigate connections between our own approach and the Lyapunov-function approaches of previous related work.  We further note that different ``relaxations" may be appropriate for answering different questions, e.g. perhaps the large deviations for the number of idle servers would be amenable to an analysis in which the linear drift appearing in the multi-dimensional OU relaxation corresponded to the drift when some servers are idle.
\\\indent Taking a broader view, such a relaxation technique can be placed in the same philosophical bucket as stochastic comparison, as in both cases one is taking a complex system, making structural modifications to its dynamics to yield a simpler system, and then provably relating the two.  It is natural to ask what are the fundamental limits of such an approach.  By considering more sophisticated comparisons and relaxations, can we bridge the gap between the high-complexity of the true limit processes and the error of any particular bound?  
Is there a principled way to develop heirarchies of such relaxations, yielding tighter-and-tighter bounds at the cost of greater complexity, and formalizing this error-complexity trade-off?

\section*{Acknowledgements.}\label{AckSec}
The authors gratefully acknowledge support from NSF grant no. 1333457, as well as several insightful conversations with Jim Dai, Ton Dieker, Shuangchi He, Shane Henderson, and Jamol Pender.

\section{Appendix.}\label{AppSec}
\subsection{Proof of Lemma\ \ref{marginals1}.}
\proof{Proof of Lemma\ \ref{marginals1}:}
In light of Lemma\ \ref{staysgood3} and the basic properties of CTMC, it suffices to prove that for any fixed state $\bigg( {\mathcal S}^{0'} , {\mathcal W}^{0'}, \nu^{0'} , \overline{\mathcal S}^{0'}, \overline{\mathcal W}^{0'} , \overline{\nu}^{0'} \bigg)$, again referred to generically using ``$\cdot$", s.t.: 1. $( {\mathcal S}^{0'} , {\mathcal W}^{0'}, \nu^{0'} )$ is 0-good, 2. $(\overline{\mathcal S}^{0'}, \overline{\mathcal W}^{0'} , \overline{\nu}^{0'} \big)$ is $\overline{0}$-good, and 3. $(\overline{\mathcal S}^{0'}, \overline{\mathcal W}^{0'} , \overline{\nu}^{0'} \big) \geq ( {\mathcal S}^{0'} , {\mathcal W}^{0'}, \nu^{0'} )$, the following are true.  Let $\TRIP$ denote the set of all possible triplets $({\mathcal S}, {\mathcal W}, \nu)$ which can appear as the state of any queue under our construction.  First, for all triplets $\big( {\mathcal S}^{0''}, {\mathcal W}^{0''}, \nu^{0''} \big)$, 
\begin{equation}\label{allsame1}
\sum_{\big(\overline{\mathcal S}^{0''}, \overline{\mathcal W}^{0''} , \overline{\nu}^{0''}\big) \in \TRIP} q'\Bigg( \cdot , \bigg({\mathcal S}^{0''}, {\mathcal W}^{0''}, \nu^{0''}, \overline{\mathcal S}^{0''}, \overline{\mathcal W}^{0''} , \overline{\nu}^{0''} \bigg) \Bigg) = q^0\bigg( \big( {\mathcal S}^{0'}, {\mathcal W}^{0'}, \nu^{0'} \big) ,\big( {\mathcal S}^{0''}, {\mathcal W}^{0''}, \nu^{0''} \big) \bigg).
\end{equation}
Second, for all $\big( \overline{\mathcal S}^{0''}, \overline{\mathcal W}^{0''}, \overline{\nu}^{0''} \big)$, 
\begin{equation}\label{allsame2}
\sum_{\big( {\mathcal S}^{0''}, {\mathcal W}^{0''} , {\nu}^{0''}\big) \in \TRIP} q'\Bigg( \cdot , \bigg({\mathcal S}^{0''}, {\mathcal W}^{0''}, \nu^{0''}, \overline{\mathcal S}^{0''}, \overline{\mathcal W}^{0''} , \overline{\nu}^{0''} \bigg) \Bigg) = 
\overline{q}^0\bigg( \big( \overline{\mathcal S}^{0'}, \overline{\mathcal W}^{0'}, \overline{\nu}^{0'} \big) ,\big( \overline{\mathcal S}^{0''}, \overline{\mathcal W}^{0''}, \overline{\nu}^{0''} \big) \bigg).
\end{equation}
We begin by proving (\ref{allsame1}).  We proceed by verifying that (\ref{allsame1}) holds for each possible transition in $q^0$, where sometimes we verify (\ref{allsame1}) by treating several cases based on the characteristics of the current state.  For clarity of exposition, we always keep track of both which transition of $q^0$ is being addressed, as well as which transitions in $q'$ (when summed) yield the appropriate rate within each case.
\\\begin{itemize}
\item (\ref{g1}): If $|{\mathcal S}^{0'}| \leq n-1$, for $z \in \lbrace 1,2 \rbrace$:
\begin{eqnarray*}
\lambda p_z &=& q^{0}\bigg( \big( {\mathcal S}^{0'}, {\mathcal W}^{0'}, \nu^{0'} \big) , \big( {\mathcal S}^{0'} \bigcup( \nu^{0'} + 1, z), {\mathcal W}^{0'}, \nu^{0'} + 1 \big) \bigg)
\\&=& q'\Bigg( \cdot, \bigg( {\mathcal S}^{0'} \bigcup( \nu^{0'} + 1, z), {\mathcal W}^{0'}, \nu^{0'} + 1, \overline{\mathcal S}^{0'}, \overline{\mathcal W}^{0'} \bigcup( \overline{\nu}^{0'} + 1, z), \overline{\nu}^{0'} + 1 \bigg)\ \ \ (\ref{z1}).
\end{eqnarray*}
\item (\ref{g2}): If $|{\mathcal S}^{0'}| = n$, for $z \in \lbrace 1,2 \rbrace$:
\begin{eqnarray*}
\lambda p_z &=& q^{0}\bigg( \big( {\mathcal S}^{0'}, {\mathcal W}^{0'}, \nu^{0'} \big) , \big( {\mathcal S}^{0'}, {\mathcal W}^{0'} \bigcup( \nu^{0'} + 1, z) , \nu^{0'} + 1 \big) \bigg)
\\&=& q'\Bigg( \cdot, \bigg( {\mathcal S}^{0'}, {\mathcal W}^{0'} \bigcup( \nu^{0'} + 1, z) , \nu^{0'} + 1, \overline{\mathcal S}^{0'}, \overline{\mathcal W}^{0'} \bigcup( \overline{\nu}^{0'} + 1, z), \overline{\nu}^{0'} + 1 \bigg)\ \ \ (\ref{z2}).
\end{eqnarray*}
\item (\ref{g3}): 
\\\\\indent\indent If ${\mathcal W}^{0'} = \emptyset$ and $\overline{\mathcal W}^{0'} = \emptyset$, for $i \geq 1$ and $z \in \lbrace 1,2 \rbrace$ s.t. $N\big( {\mathcal S}^{0'}, i, z \big) = 1$ and $N\big( \overline{\mathcal S}^{0'}, i, z \big) = 1$:
\begin{eqnarray*}
\mu_z &=& q^{0}\bigg( \big( {\mathcal S}^{0'}, {\mathcal W}^{0'}, \nu^{0'} \big) , \big( {\mathcal S}^{0'} \setminus(i,z), {\mathcal W}^{0'}, \nu^{0'} \big) \bigg) 
\\&=& \sum_{z_2 = 1}^2 q'\Bigg( \cdot, \bigg( {\mathcal S}^{0'} \setminus (i,z), {\mathcal W}^{0'}, \nu^{0'}, \big( \overline{\mathcal S}^{0'} \setminus (i,z) \big) \bigcup (0,z_2), \overline{\mathcal W}^{0'}, \overline{\nu}^{0'} \bigg) \Bigg)\ \ \ (\ref{z4}).
\end{eqnarray*}
\indent\indent If ${\mathcal W}^{0'} = \emptyset$ and $\overline{\mathcal W}^{0'} \neq \emptyset$, for $i \geq 1$ , $z \in \lbrace 1,2 \rbrace$ s.t. $N\big( {\mathcal S}^{0'}, i, z \big) = 1$ and $N\big( \overline{\mathcal S}^{0'}, i, z \big) = 1$:
\begin{eqnarray*}
\mu_z &=& q^{0}\bigg( \big( {\mathcal S}^{0'}, {\mathcal W}^{0'}, \nu^{0'} \big) , \big( {\mathcal S}^{0'} \setminus(i,z), {\mathcal W}^{0'}, \nu^{0'} \big) \bigg) 
\\&=& q'\Bigg( \cdot, \bigg( {\mathcal S}^{0'} \setminus (i,z), {\mathcal W}^{0'}, \nu^{0'}, \big( \overline{\mathcal S}^{0'} \setminus (i,z) \big) \bigcup \overline{\mathcal W}^{0'}_1, \overline{\mathcal W}^{0'} \setminus \overline{\mathcal W}^{0'}_1, \overline{\nu}^{0'} \bigg) \Bigg)\ \ \ (\ref{z8}).
\end{eqnarray*}
\indent\indent If ${\mathcal W}^{0'} = \emptyset$ and $\overline{\mathcal W}^{0'} \neq \emptyset$, for $i \geq 1$, $z \in \lbrace 1,2 \rbrace$ s.t. $N\big( {\mathcal S}^{0'}, i, z \big) = 1$ and $N\big( \overline{\mathcal W}^{0'}, i, z \big) = 1$:
\begin{eqnarray*}
\mu_z &=& q^{0}\bigg( \big( {\mathcal S}^{0'}, {\mathcal W}^{0'}, \nu^{0'} \big) , \big( {\mathcal S}^{0'} \setminus(i,z), {\mathcal W}^{0'}, \nu^{0'} \big) \bigg) 
\\&=& q'\Bigg( \cdot , \bigg( {\mathcal S}^{0'} \setminus (i,z) , {\mathcal W}^{0'} , \nu^{0'}, \overline{\mathcal S}^{0'} , \overline{\mathcal W}^{0'} \setminus (i,z), \overline{\nu}^{0'} \bigg) \Bigg)\ \ \ (\ref{z9})
\\&\ &\ \ \ +\ \ \ q'\Bigg( \cdot, \bigg( {\mathcal S}^{0'} \setminus (i,z) , {\mathcal W}^{0'} , \nu^{0'}, \overline{\mathcal S}^{0'} , \overline{\mathcal W}^{0'} , \overline{\nu}^{0'} \bigg) \Bigg)\ \ \ (\ref{z10}).
\end{eqnarray*}
\item (\ref{g4}): 
\\\\\indent\indent If ${\mathcal W}^{0'} \neq \emptyset$ and $\overline{\mathcal W}^{0'} \neq \emptyset$, for $i \geq 1$ and $z \in \lbrace 1,2 \rbrace$ s.t. $N\big( {\mathcal S}^{0'}, i, z \big) = 1$ and $N\big( \overline{\mathcal S}^{0'}, i, z \big) = 1$:
\begin{eqnarray*}
\ \ \ \ \ \ \ \ \ \ \ \ \mu_z &=& q^{0}\bigg( \big( {\mathcal S}^{0'}, {\mathcal W}^{0'}, \nu^{0'} \big) , \big( \big({\mathcal S}^{0'} \setminus(i,z) \big) \bigcup {\mathcal W}^{0'}_1, {\mathcal W}^{0'} \setminus {\mathcal W}^{0'}_1, \nu^{0'} \big) \bigg) 
\\&=& q'\Bigg( \cdot, \bigg( \big( {\mathcal S}^{0'} \setminus (i,z) \big) \bigcup {\mathcal W}^{0'}_1, {\mathcal W}^{0'} \setminus {\mathcal W}^{0'}_1, \nu^{0'}, \big( \overline{\mathcal S}^{0'} \setminus (i,z) \big) \bigcup \overline{\mathcal W}^{0'}_1, \overline{\mathcal W}^{0'} \setminus \overline{\mathcal W}^{0'}_1, \overline{\nu}^{0'} \bigg) \Bigg)\ \ \ (\ref{z16}).
\end{eqnarray*}
\indent\indent If ${\mathcal W}^{0'} \neq \emptyset$ and $\overline{\mathcal W}^{0'} \neq \emptyset$, for $i \geq 1$ and $z \in \lbrace 1,2 \rbrace$ s.t. $N\big( {\mathcal S}^{0'}, i, z \big) = 1$ and $N\big( \overline{\mathcal W}^{0'}, i, z \big) = 1$:
\begin{eqnarray*}
\mu_z &=& q^{0}\bigg( \big( {\mathcal S}^{0'}, {\mathcal W}^{0'}, \nu^{0'} \big) , \big( \big({\mathcal S}^{0'} \setminus(i,z) \big) \bigcup {\mathcal W}^{0'}_1, {\mathcal W}^{0'} \setminus {\mathcal W}^{0'}_1, \nu^{0'} \big) \bigg) 
\\&=& q'\Bigg( \cdot , \bigg( \big( {\mathcal S}^{0'} \setminus (i,z) \big) \bigcup {\mathcal W}^{0'}_1 , {\mathcal W}^{0'} \setminus {\mathcal W}^{0'}_1, \nu^{0'}, \overline{\mathcal S}^{0'} , \overline{\mathcal W}^{0'} \setminus (i,z), \overline{\nu}^{0'} \bigg) \Bigg)\ \ \ (\ref{z17})
\\&\ &\ \ \ +\ \ \ q'\Bigg( \cdot, \bigg( \big( {\mathcal S}^{0'} \setminus (i,z) \big) \bigcup {\mathcal W}^{0'}_1 , {\mathcal W}^{0'} \setminus {\mathcal W}^{0'}_1 , \nu^{0'}, \overline{\mathcal S}^{0'} , \overline{\mathcal W}^{0'} , \overline{\nu}^{0'} \bigg) \Bigg)\ \ \ (\ref{z18}).
\end{eqnarray*}
\item (\ref{g5}): 
\\\\\indent\indent If ${\mathcal W}^{0'} \neq \emptyset$ and $\overline{\mathcal W}^{0'} \neq \emptyset$, for $i \geq 1$ and $z \in \lbrace 1,2 \rbrace$ s.t. $N\big( {\mathcal W}^{0'}, i, z \big) = 1$ and $N\big( \overline{\mathcal W}^{0'}, i, z \big) = 1$:
\begin{eqnarray*}
\theta &=& q^{0}\bigg( \big( {\mathcal S}^{0'}, {\mathcal W}^{0'}, \nu^{0'} \big) , \big( {\mathcal S}^{0'} , {\mathcal W}^{0'} \setminus (i,z), \nu^{0'} \big) \bigg) 
\\&=& q'\Bigg( \cdot, \bigg( {\mathcal S}^{0'} , {\mathcal W}^{0'} \setminus (i,z) , \nu^{0'}, \overline{\mathcal S}^{0'} , \overline{\mathcal W}^{0'} \setminus (i,z), \overline{\nu}^{0'} \bigg) \Bigg)\ \ \ (\ref{z19}).
\end{eqnarray*}
\end{itemize}
Combining with the fact that under those transitions in $q'$ not noted above, i.e. (\ref{z3}),(\ref{z5}),(\ref{z6}),(\ref{z7}),(\ref{z12}),(\ref{z14}),(\ref{z15}),(\ref{z20}), it holds that $\big( {\mathcal S}^{0''}, {\mathcal W}^{0''}, \nu^{0''} \big) = \big( {\mathcal S}^{0'}, {\mathcal W}^{0'}, \nu^{0'} \big)$ (i.e. these are self-loops w.r.t. the first 3 components), and combining with the fact that $(\overline{\mathcal S}^{0'}, \overline{\mathcal W}^{0'} , \overline{\nu}^{0'} \big) \geq ( {\mathcal S}^{0'} , {\mathcal W}^{0'}, \nu^{0'} ),$ combined with (\ref{useit1}) and (\ref{useit2}), precludes certain configurations, completes the proof of (\ref{allsame1}).
\\\\\indent We now provide a similar analysis for $\overline{q}^0$, proving (\ref{allsame2}).
\\\\\begin{itemize}
\item (\ref{h1}): 
\\\\\indent\indent If $|{\mathcal S}^{0'}| \leq n-1$, for $z \in \lbrace 1,2 \rbrace$:
\begin{eqnarray*}
\lambda p_z &=& \overline{q}^{0}\bigg( \big( \overline{\mathcal S}^{0'}, \overline{\mathcal W}^{0'}, \overline{\nu}^{0'} \big) , \big( \overline{\mathcal S}^{0'}, \overline{\mathcal W}^{0'} \bigcup( \overline{\nu}^{0'} + 1, z), \overline{\nu}^{0'} + 1 \big) \bigg)
\\&=& q'\Bigg( \cdot, \bigg( {\mathcal S}^{0'} \bigcup( \nu^{0'} + 1, z), {\mathcal W}^{0'}, \nu^{0'} + 1, \overline{\mathcal S}^{0'}, \overline{\mathcal W}^{0'} \bigcup( \overline{\nu}^{0'} + 1, z), \overline{\nu}^{0'} + 1 \bigg)\ \ \ (\ref{z1}).
\end{eqnarray*}
\indent\indent If $|{\mathcal S}^{0'}| = n$, for $z \in \lbrace 1,2 \rbrace$:
\begin{eqnarray*}
\lambda p_z &=& \overline{q}^{0}\bigg( \big( \overline{\mathcal S}^{0'}, \overline{\mathcal W}^{0'}, \overline{\nu}^{0'} \big) , \big( \overline{\mathcal S}^{0'}, \overline{\mathcal W}^{0'} \bigcup( \overline{\nu}^{0'} + 1, z), \overline{\nu}^{0'} + 1 \big) \bigg)
\\&=& q'\Bigg( \cdot, \bigg( {\mathcal S}^{0'}, {\mathcal W}^{0'} \bigcup( \nu^{0'} + 1, z) , \nu^{0'} + 1, \overline{\mathcal S}^{0'}, \overline{\mathcal W}^{0'} \bigcup( \overline{\nu}^{0'} + 1, z), \overline{\nu}^{0'} + 1 \bigg)\ \ \ (\ref{z2}).
\end{eqnarray*}
\item (\ref{h2}): 
\\\\\indent\indent If $\overline{\mathcal W}^{0'} = \emptyset$ and ${\mathcal W}^{0'} = \emptyset$, for $i \geq 1$ and $z_1,z_2 \in \lbrace 1,2 \rbrace$ s.t. $N\big( \overline{\mathcal S}^{0'}, i, z_1 \big) = 1$, $N\big( {\mathcal S}^{0'}, i, z_1 \big) +
\\\indent\indent N\big( {\mathcal W}^{0'}, i, z_1 \big)  = 0$:
\begin{eqnarray*}
\mu_{z_1} p_{z_2} &=& \overline{q}^{0}\Bigg( \bigg( \overline{\mathcal S}^{0'}, \overline{\mathcal W}^{0'}, \overline{\nu}^{0'} \bigg) , \bigg( \big( \overline{\mathcal S}^{0'} \setminus(i,z_1) \big) \bigcup (0,z_2), \overline{\mathcal W}^{0'}, \overline{\nu}^{0'} \bigg) \Bigg)
\\&=& q'\Bigg( \cdot, \bigg( {\mathcal S}^{0'} , {\mathcal W}^{0'} , \nu^{0'}, \big( \overline{\mathcal S}^{0'} \setminus (i,z_1) \big) \bigcup (0,z_2), \overline{\mathcal W}^{0'} , \overline{\nu}^{0'} \bigg) \Bigg)\ \ \ (\ref{z3}).
\end{eqnarray*}
\indent\indent If $\overline{\mathcal W}^{0'} = \emptyset$ and ${\mathcal W}^{0'} = \emptyset$, for $i \geq 1$ and $z_1,z_2 \in \lbrace 1,2 \rbrace$ s.t. $N\big( \overline{\mathcal S}^{0'}, i, z_1 \big) = 1$, $N\big( {\mathcal S}^{0'}, i, z_1 \big) = 1$:
\begin{eqnarray*}
\mu_{z_1} p_{z_2} &=& \overline{q}^{0}\Bigg( \bigg( \overline{\mathcal S}^{0'}, \overline{\mathcal W}^{0'}, \overline{\nu}^{0'} \bigg) , \bigg( \big( \overline{\mathcal S}^{0'} \setminus(i,z_1) \big) \bigcup (0,z_2), \overline{\mathcal W}^{0'}, \overline{\nu}^{0'} \bigg) \Bigg)
\\&=& q'\Bigg( \cdot, \bigg( {\mathcal S}^{0'} \setminus (i,z_1), {\mathcal W}^{0'}, \nu^{0'}, \big( \overline{\mathcal S}^{0'} \setminus (i,z_1) \big) \bigcup (0,z_2), \overline{\mathcal W}^{0'}, \overline{\nu}^{0'} \bigg) \Bigg)\ \ \ (\ref{z4}).
\end{eqnarray*}
\item (\ref{h3}): 
\\\\\indent\indent If $\overline{\mathcal W}^{0'} = \emptyset$ and ${\mathcal W}^{0'} = \emptyset$, for $z_1,z_2 \in \lbrace 1,2 \rbrace$ s.t. $z_1 \neq z_2$:
\begin{eqnarray*}
N\big( \overline{\mathcal S}^{0'}, 0, z_1 \big) \mu_{z_1} p_{z_2} &=& \overline{q}^{0}\Bigg( \bigg( \overline{\mathcal S}^{0'}, \overline{\mathcal W}^{0'}, \overline{\nu}^{0'} \bigg) , \bigg( \big( \overline{\mathcal S}^{0'} \setminus(0,z_1) \big) \bigcup (0,z_2), \overline{\mathcal W}^{0'}, \overline{\nu}^{0'} \bigg) \Bigg)
\\&=& q'\Bigg( \cdot, \bigg( {\mathcal S}^{0'} , {\mathcal W}^{0'} , \nu^{0'} , \big( \overline{\mathcal S}^{0'} \setminus (0,z_1) \big) \bigcup (0,z_2) , \overline{\mathcal W}^{0'} , \overline{\nu}^{0'} \bigg) \Bigg)\ \ \ (\ref{z5}).
\end{eqnarray*}
\item (\ref{h4}): 
\\\\\indent\indent If $\overline{\mathcal W}^{0'} \neq \emptyset$ and ${\mathcal W}^{0'} = \emptyset$, for $i \geq 1, z \in \lbrace 1,2 \rbrace$ s.t. $N\big(\overline{\mathcal S}^{0'},i,z\big) = 1 , N\big({\mathcal S}^{0'},i,z\big) + N\big({\mathcal W}^{0'},i,z\big) = 0$:
\begin{eqnarray*}
\mu_z &=& \overline{q}^{0}\Bigg( \bigg( \overline{\mathcal S}^{0'}, \overline{\mathcal W}^{0'}, \overline{\nu}^{0'} \bigg) , \bigg( \big( \overline{\mathcal S}^{0'} \setminus (i,z) \big) \bigcup \overline{\mathcal W}^{0'}_1, \overline{\mathcal W}^{0'} \setminus \overline{\mathcal W}^{0'}_1, \overline{\nu}^{0'} \bigg) \Bigg)
\\&=& q'\Bigg( \cdot , \bigg( {\mathcal S}^{0'} , {\mathcal W}^{0'} , \nu^{0'}, \big( \overline{\mathcal S}^{0'} \setminus (i,z) \big) \bigcup \overline{\mathcal W}^{0'}_1, \overline{\mathcal W}^{0'} \setminus \overline{\mathcal W}^{0'}_1 , \overline{\nu}^{0'} \bigg) \Bigg)
\ \ \ (\ref{z6}).
\end{eqnarray*}
\indent\indent If $\overline{\mathcal W}^{0'} \neq \emptyset$ and ${\mathcal W}^{0'} = \emptyset$, for $i \geq 1, z \in \lbrace 1,2 \rbrace$ s.t. $N\big(\overline{\mathcal S}^{0'},i,z\big) = 1 , N\big({\mathcal S}^{0'},i,z\big) = 1$:
\begin{eqnarray*}
\mu_z &=& \overline{q}^{0}\Bigg( \bigg( \overline{\mathcal S}^{0'}, \overline{\mathcal W}^{0'}, \overline{\nu}^{0'} \bigg) , \bigg( \big( \overline{\mathcal S}^{0'} \setminus (i,z) \big) \bigcup \overline{\mathcal W}^{0'}_1, \overline{\mathcal W}^{0'} \setminus \overline{\mathcal W}^{0'}_1, \overline{\nu}^{0'} \bigg) \Bigg)
\\&=& q'\Bigg( \cdot, \bigg( {\mathcal S}^{0'} \setminus (i,z), {\mathcal W}^{0'}, \nu^{0'}, \big( \overline{\mathcal S}^{0'} \setminus (i,z) \big) \bigcup \overline{\mathcal W}^{0'}_1, \overline{\mathcal W}^{0'} \setminus \overline{\mathcal W}^{0'}_1, \overline{\nu}^{0'} \bigg) \Bigg)
\ \ \ (\ref{z8}).
\end{eqnarray*}
\indent\indent If $\overline{\mathcal W}^{0'} \neq \emptyset$ and ${\mathcal W}^{0'} \neq \emptyset$, for $i \geq 1, z \in \lbrace 1,2 \rbrace$ s.t. $N\big(\overline{\mathcal S}^{0'},i,z\big) = 1 , N\big({\mathcal S}^{0'},i,z\big) + N\big({\mathcal W}^{0'},i,z\big) = 0:$
\begin{eqnarray*}
\mu_z &=& \overline{q}^{0}\Bigg( \bigg( \overline{\mathcal S}^{0'}, \overline{\mathcal W}^{0'}, \overline{\nu}^{0'} \bigg) , \bigg( \big( \overline{\mathcal S}^{0'} \setminus (i,z) \big) \bigcup \overline{\mathcal W}^{0'}_1, \overline{\mathcal W}^{0'} \setminus \overline{\mathcal W}^{0'}_1, \overline{\nu}^{0'} \bigg) \Bigg)
\\&=& q'\Bigg( \cdot , \bigg( {\mathcal S}^{0'} , {\mathcal W}^{0'} , \nu^{0'}, \big( \overline{\mathcal S}^{0'} \setminus (i,z) \big) \bigcup \overline{\mathcal W}^{0'}_1, \overline{\mathcal W}^{0'} \setminus \overline{\mathcal W}^{0'}_1 , \overline{\nu}^{0'} \bigg) \Bigg)
\ \ \ (\ref{z14}).
\end{eqnarray*}
\indent\indent If $\overline{\mathcal W}^{0'} \neq \emptyset$ and ${\mathcal W}^{0'} \neq \emptyset$, for $i \geq 1, z \in \lbrace 1,2 \rbrace$ s.t. $N\big(\overline{\mathcal S}^{0'},i,z\big) = 1 , N\big({\mathcal S}^{0'},i,z\big) = 1$:
\begin{eqnarray*}
\ \ \ \ \ \ \mu_z &=& \overline{q}^{0}\Bigg( \bigg( \overline{\mathcal S}^{0'}, \overline{\mathcal W}^{0'}, \overline{\nu}^{0'} \bigg) , \bigg( \big( \overline{\mathcal S}^{0'} \setminus (i,z) \big) \bigcup \overline{\mathcal W}^{0'}_1, \overline{\mathcal W}^{0'} \setminus \overline{\mathcal W}^{0'}_1, \overline{\nu}^{0'} \bigg) \Bigg)
\\&=& q'\Bigg( \cdot, \bigg( \big( {\mathcal S}^{0'} \setminus (i,z) \big) \bigcup {\mathcal W}^{0'}_1, {\mathcal W}^{0'} \setminus {\mathcal W}^{0'}_1, \nu^{0'}, \big( \overline{\mathcal S}^{0'} \setminus (i,z) \big) \bigcup \overline{\mathcal W}^{0'}_1, \overline{\mathcal W}^{0'} \setminus \overline{\mathcal W}^{0'}_1, \overline{\nu}^{0'} \bigg) \Bigg)
\ \ \ (\ref{z16}).
\end{eqnarray*}
\item (\ref{h5}): 
\\\\\indent\indent If $\overline{\mathcal W}^{0'} \neq \emptyset$ and ${\mathcal W}^{0'} = \emptyset$, for $i \geq 1, z \in \lbrace 1,2 \rbrace$ s.t. $N\big(\overline{\mathcal W}^{0'},i,z\big) = 1 , N\big({\mathcal S}^{0'},i,z\big) + N\big({\mathcal W}^{0'},i,z\big) = 0$:
\begin{eqnarray*}
\theta &=& \overline{q}^{0}\Bigg( \bigg( \overline{\mathcal S}^{0'}, \overline{\mathcal W}^{0'}, \overline{\nu}^{0'} \bigg) , \bigg( \overline{\mathcal S}^{0'}, \overline{\mathcal W}^{0'} \setminus (i,z), \overline{\nu}^{0'} \bigg) \Bigg)
\\&=& q'\Bigg( \cdot, \bigg( {\mathcal S}^{0'} , {\mathcal W}^{0'} , \nu^{0'}, \overline{\mathcal S}^{0'} , \overline{\mathcal W}^{0'} \setminus (i,z), \overline{\nu}^{0'} \bigg) \Bigg)\ \ \ (\ref{z7}).
\end{eqnarray*}
\indent\indent If $\overline{\mathcal W}^{0'} \neq \emptyset$ and ${\mathcal W}^{0'} = \emptyset$, for $i \geq 1, z \in \lbrace 1,2 \rbrace$ s.t. $N\big(\overline{\mathcal W}^{0'},i,z\big) = 1 , N\big({\mathcal S}^{0'},i,z\big) = 1$:
\begin{eqnarray*}
\theta &=& \overline{q}^{0}\Bigg( \bigg( \overline{\mathcal S}^{0'}, \overline{\mathcal W}^{0'}, \overline{\nu}^{0'} \bigg) , \bigg( \overline{\mathcal S}^{0'}, \overline{\mathcal W}^{0'} \setminus (i,z), \overline{\nu}^{0'} \bigg) \Bigg)
\\&=& q'\Bigg( \cdot , \bigg( {\mathcal S}^{0'} \setminus (i,z) , {\mathcal W}^{0'} , \nu^{0'}, \overline{\mathcal S}^{0'} , \overline{\mathcal W}^{0'} \setminus (i,z), \overline{\nu}^{0'} \bigg) \Bigg)\ \ \ (\ref{z9}).
\end{eqnarray*}
\indent\indent If $\overline{\mathcal W}^{0'} \neq \emptyset$ and ${\mathcal W}^{0'} \neq \emptyset$, for $i \geq 1, z \in \lbrace 1,2 \rbrace$ s.t. $N\big(\overline{\mathcal W}^{0'},i,z\big) = 1 , N\big({\mathcal S}^{0'},i,z\big) + N\big({\mathcal W}^{0'},i,z\big) = 0$:
\begin{eqnarray*}
\theta &=& \overline{q}^{0}\Bigg( \bigg( \overline{\mathcal S}^{0'}, \overline{\mathcal W}^{0'}, \overline{\nu}^{0'} \bigg) , \bigg( \overline{\mathcal S}^{0'}, \overline{\mathcal W}^{0'} \setminus (i,z), \overline{\nu}^{0'} \bigg) \Bigg)
\\&=& q'\Bigg( \cdot, \bigg( {\mathcal S}^{0'} , {\mathcal W}^{0'} , \nu^{0'}, \overline{\mathcal S}^{0'} , \overline{\mathcal W}^{0'} \setminus (i,z), \overline{\nu}^{0'} \bigg) \Bigg)\ \ \ (\ref{z15}).
\end{eqnarray*}
\indent\indent If $\overline{\mathcal W}^{0'} \neq \emptyset$ and ${\mathcal W}^{0'} \neq \emptyset$, for $i \geq 1, z \in \lbrace 1,2 \rbrace$ s.t. $N\big(\overline{\mathcal W}^{0'},i,z\big) = 1 , N\big({\mathcal S}^{0'},i,z\big) = 1$:
\begin{eqnarray*}
\theta &=& \overline{q}^{0}\Bigg( \bigg( \overline{\mathcal S}^{0'}, \overline{\mathcal W}^{0'}, \overline{\nu}^{0'} \bigg) , \bigg( \overline{\mathcal S}^{0'}, \overline{\mathcal W}^{0'} \setminus (i,z), \overline{\nu}^{0'} \bigg) \Bigg)
\\&=& q'\Bigg( \cdot , \bigg( \big( {\mathcal S}^{0'} \setminus (i,z) \big) \bigcup {\mathcal W}^{0'}_1 , {\mathcal W}^{0'} \setminus {\mathcal W}^{0'}_1, \nu^{0'}, \overline{\mathcal S}^{0'} , \overline{\mathcal W}^{0'} \setminus (i,z), \overline{\nu}^{0'} \bigg) \Bigg)\ \ \ (\ref{z17}).
\end{eqnarray*}
\indent\indent If $\overline{\mathcal W}^{0'} \neq \emptyset$ and ${\mathcal W}^{0'} \neq \emptyset$, for $i \geq 1, z \in \lbrace 1,2 \rbrace$ s.t. $N\big(\overline{\mathcal W}^{0'},i,z\big) = 1 , N\big({\mathcal W}^{0'},i,z\big) = 1$:
\begin{eqnarray*}
\theta &=& \overline{q}^{0}\Bigg( \bigg( \overline{\mathcal S}^{0'}, \overline{\mathcal W}^{0'}, \overline{\nu}^{0'} \bigg) , \bigg( \overline{\mathcal S}^{0'}, \overline{\mathcal W}^{0'} \setminus (i,z), \overline{\nu}^{0'} \bigg) \Bigg)
\\&=& q'\Bigg( \cdot, \bigg( {\mathcal S}^{0'} , {\mathcal W}^{0'} \setminus (i,z) , \nu^{0'}, \overline{\mathcal S}^{0'} , \overline{\mathcal W}^{0'} \setminus (i,z), \overline{\nu}^{0'} \bigg) \Bigg)\ \ \ (\ref{z19}).
\end{eqnarray*}
\item (\ref{h6}): 
\\\\\indent\indent If $\overline{\mathcal W}^{0'} \neq \emptyset$ and ${\mathcal W}^{0'} = \emptyset$, for $z \in \lbrace 1,2 \rbrace$:
\begin{eqnarray*}
N\big( \overline{\mathcal S}^{0'}, 0, z \big) \mu_z &=& \overline{q}^{0}\Bigg( \bigg( \overline{\mathcal S}^{0'}, \overline{\mathcal W}^{0'}, \overline{\nu}^{0'} \bigg) , \bigg( \big( \overline{\mathcal S}^{0'} \setminus (0,z) \big) \bigcup \overline{\mathcal W}^{0'}_1, \overline{\mathcal W}^{0'} \setminus \overline{\mathcal W}^{0'}_1, \nu \bigg) \Bigg) 
\\&=& q'\Bigg( \cdot , \bigg( {\mathcal S}^{0'} , {\mathcal W}^{0'} , \nu^{0'} , \big( \overline{\mathcal S}^{0'} \setminus (0,z) \big) \bigcup \overline{\mathcal W}^{0'}_1 , \overline{\mathcal W}^{0'} \setminus \overline{\mathcal W}^{0'}_1 , \overline{\nu}^{0'} \bigg) \Bigg)\ \ \ (\ref{z12}).
\end{eqnarray*}
\indent\indent If $\overline{\mathcal W}^{0'} \neq \emptyset$ and ${\mathcal W}^{0'} \neq \emptyset$, for $z \in \lbrace 1,2 \rbrace$:
\begin{eqnarray*}
N\big( \overline{\mathcal S}^{0'}, 0, z \big) \mu_z &=& \overline{q}^{0}\Bigg( \bigg( \overline{\mathcal S}^{0'}, \overline{\mathcal W}^{0'}, \overline{\nu}^{0'} \bigg) , \bigg( \big( \overline{\mathcal S}^{0'} \setminus (0,z) \big) \bigcup \overline{\mathcal W}^{0'}_1, \overline{\mathcal W}^{0'} \setminus \overline{\mathcal W}^{0'}_1, \nu \bigg) \Bigg) 
\\&=& q'\Bigg( \cdot , \bigg( {\mathcal S}^{0'} , {\mathcal W}^{0'} , \nu^{0'} , \big( \overline{\mathcal S}^{0'} \setminus (0,z) \big) \bigcup \overline{\mathcal W}^{0'}_1 , \overline{\mathcal W}^{0'} \setminus \overline{\mathcal W}^{0'}_1 , \overline{\nu}^{0'} \bigg) \Bigg)\ \ \ (\ref{z20}).
\end{eqnarray*}
\end{itemize}
Combining with the fact that under those transitions in $q'$ not noted above, i.e. (\ref{z10}),(\ref{z18}), it holds that $\big( \overline{\mathcal S}^{0''}, \overline{\mathcal W}^{0''}, \overline{\nu}^{0''} \big) = \big( \overline{\mathcal S}^{0'}, \overline{\mathcal W}^{0'}, \overline{\nu}^{0'} \big)$ (i.e. these are self-loops w.r.t. the latter 3 components), and combining with the fact that $(\overline{\mathcal S}^{0'}, \overline{\mathcal W}^{0'} , \overline{\nu}^{0'} \big) \geq ( {\mathcal S}^{0'} , {\mathcal W}^{0'}, \nu^{0'} ),$ combined with (\ref{useit1}) and (\ref{useit2}), precludes certain configurations, completes the proof of (\ref{allsame2}).  Combining the above completes the proof.  $\qed$
\endproof

\subsection{Proof of Lemma\ \ref{weaklem2}.}
\proof{Proof of Lemma\ \ref{weaklem2}:}
Let $\tilde{Q}^n = \big( \tilde{\RES}^n_1, \tilde{\RES}^n_2, \tilde{W}^n, \tilde{\ARR}^n, \tilde{REN}^n, \tilde{\ABA}^n, \tilde{\DMY}^n \big)$ denote the CTMC $\tilde{Q}$ for the given $n$, when $\lambda = \lambda_n = n + B n^{\frac{1}{2}},$ and $\tilde{Q}^n(0) = \big( \overline{N}^n_1(0), \overline{N}^n_2(0), \overline{W}^n(0), 0, 0, 0, 0 \big)$.  It follows from a proof essentially identical to that of Lemma\ \ref{skoro1}, the details of which we omit, that one may construct $\tilde{Q}^n$ and $\CLOCK$ on a common probability space s.t. w.p.1, for all $t \geq 0$,
\begin{itemize}
\item $\tilde{\ABA}^n(t) = \CLOCK\big( \theta \int_0^t \tilde{W}^n(s) ds \big)$;
\item $\tilde{W}^n(t) = \tilde{\ARR}^n(t) - \tilde{\REN}^n(t) - \tilde{\ABA}^n(t) + \tilde{\DMY}^n(t)$;
\item $\int_0^t I\big( \tilde{W}^n(s) > 0 \big) d\tilde{\DMY}^n(s) = 0.$
\end{itemize}
On the same probability space, let $\tilde{\ERR}^n$ denote the stochastic process s.t. for all $t \geq 0$,
$$\tilde{\ERR}^n(t) = n^{-\frac{1}{2}} \bigg( \CLOCK\big( \theta \int_0^t \tilde{W}^n(s) ds \big) - \theta \int_0^t \tilde{W}^n(s) ds \bigg).$$
By Lemma\ \ref{sameagain1}, it suffices to prove the desired weak convergence for $\big\lbrace \tilde{\ERR}^n, n \geq 1 \big\rbrace$.
Note that 
\begin{eqnarray} 
\big|\tilde{\ERR}^n(y)\big| &=& n^{-\frac{1}{2}} \bigg| \theta \int_0^y \tilde{W}^n(s) ds - \CLOCK\big(\theta \int_0^y \tilde{W}^n(s) ds \big) \bigg| \nonumber
\\&\leq& n^{-\frac{1}{2}} \sup_{0 \leq x \leq \int_0^y \tilde{W}^n(s) ds} |\theta x - \CLOCK(\theta x)| \nonumber
\\&\leq& n^{-\frac{1}{2}} \sup_{0 \leq x \leq y \times \sup_{0 \leq s \leq y} \tilde{W}^n(s)} |\theta x - \CLOCK(\theta x)|. \label{usetozero1}
\end{eqnarray}
For each fixed $y > 0$, Let $T(y,n) \stackrel{\Delta}{=} \sup\lbrace s \in [0,y] : \tilde{W}^n(s) = 0 \rbrace$.  Recall that for a stochastic process $Z$, $Z(t^-)$ denotes the appropriate left limit (when it exists).  Note that the definition of $T(y,n)$ and the defining properties of $\tilde{\DMY}^n$ imply that $\tilde{\DMY}^n\big(T(y,n)^-\big) = \tilde{\DMY}^n(y)$.  Also, it follows from the monotonicity of $\tilde{\ABA}^n$ that $\tilde{\ABA}^n\big(T(y,n)^-\big) \leq \tilde{\ABA}^n(y)$.  Combining with Lemma\ \ref{sameagain4}, we conclude that we may construct all relevant processes on a common probability space s.t. w.p.1, for all $y \geq 0$, 
\begin{eqnarray*}
\tilde{W}^n(y) &=& \tilde{W}^n\big( T(y,n)^- \big) + \bigg( \tilde{W}^n\big( y \big) - \tilde{W}^n\big( T(y,n)^- \big) \bigg)
\\&=& \bigg( \tilde{\ARR}^n(y) - \tilde{\REN}^n(y) - \tilde{\ABA}^n(y) + \tilde{\DMY}^n(y) \bigg) 
\\&\ &\ \ \ \ \ \ \ \ - \bigg( \tilde{\ARR}^n\big(T(y,n)^-\big) - \tilde{\REN}^n\big(T(y,n)^-\big) - \tilde{\ABA}^n\big(T(y,n)^-\big) + \tilde{\DMY}^n\big(T(y,n)^-\big) \bigg)
\\&=& \bigg( \tilde{\ARR}^n(y) - \tilde{\REN}^n(y) - \tilde{\ABA}^n(y) \bigg) 
\\&\ &\ \ \ \ \ \ \ \ - \bigg( \tilde{\ARR}^n\big(T(y,n)^-\big) - \tilde{\REN}^n\big(T(y,n)^-\big) - \tilde{\ABA}^n\big(T(y,n)^-\big) \bigg)
\\&\leq& \bigg( \tilde{\ARR}^n(y) - \tilde{\REN}^n(y) \bigg) - \bigg( \tilde{\ARR}^n\big(T(y,n)^-\big) - \tilde{\REN}^n\big(T(y,n)^-\big) \bigg)
\\&\leq& \sup_{0 \leq s \leq y}\Bigg( \big( \tilde{\ARR}^n(y) - \tilde{\REN}^n(y) \big) - \big( \tilde{\ARR}^n(s) - \tilde{\REN}^n(s) \big) \Bigg),
\end{eqnarray*}
implying that for $0 \leq y \leq t$, $\sup_{0 \leq s \leq y} \tilde{W}^n(s)$ is at most
\begin{eqnarray*}
\ &\ &\ \ \sup_{0 \leq r \leq s \leq y} \Bigg( \big( \tilde{\ARR}^n(s) - \tilde{\REN}^n(s) \big) - \big( \tilde{\ARR}^n(r) - \tilde{\REN}^n(r) \big) \Bigg)
\\&\ &\ \ \ \ \leq\ \ \ 2 \sup_{0 \leq s \leq t} \big| \tilde{\ARR}^n(s) - \tilde{\REN}^n(s) \big|.
\end{eqnarray*}
Combining the above, we conclude that for all $t \geq 0$ and $y \in [0,t]$,
$$
\big|\tilde{\ERR}^n(y)\big| \leq n^{-\frac{1}{2}} \sup_{0 \leq x \leq 2 t \times \sup_{0 \leq s \leq t} \big| \tilde{\ARR}^n(s) - \tilde{\REN}^n(s) \big| } |\theta x - \CLOCK(\theta x)|,$$
and hence
$$\sup_{0 \leq y \leq t} \big|\tilde{\ERR}^n(y)\big| \leq n^{-\frac{1}{2}} \sup_{0 \leq x \leq 2 t \times \sup_{0 \leq s \leq t} \big| \tilde{\ARR}^n(s) - \tilde{\REN}^n(s) \big| } |\theta x - \CLOCK(\theta x)|.$$
It then follows from a union bound that for all $T \geq 0$, $P\big( \sup_{0 \leq y \leq T} \big|\tilde{\ERR}^n(y)\big| > n^{-\frac{1}{5}} \big)$ is at most
\begin{eqnarray}
\ &\ &\ P\bigg( 2 T \sup_{0 \leq s \leq T} \big| \tilde{\ARR}^n(s) - \tilde{\REN}^n(s) \big| \geq n^{\frac{11}{20}} \bigg) \nonumber
\\&\ &\ \ \ \ +\ \ \ P\bigg( n^{-\frac{1}{2}} \sup_{0 \leq x \leq n^{\frac{11}{20}} } |\theta x - \CLOCK(\theta x)| \geq n^{-\frac{1}{5}} \bigg) \nonumber
\\&=&\ P\bigg( 2 T \sup_{0 \leq s \leq T} \big| n^{-\frac{1}{2}} \big( \tilde{\ARR}^n(s) - \tilde{\REN}^n(s) \big) \big| \geq n^{\frac{1}{20}} \bigg) \label{tozero1}
\\&\ &\ \ \ \ +\ \ \ P\bigg( \sup_{0 \leq x \leq n^{\frac{11}{20}} } |\theta x - \CLOCK(\theta x)| \geq n^{\frac{3}{10}} \bigg). \label{tozero2}
\end{eqnarray}
It follows from Lemma\ \ref{weaklem1}, the continuity of the supremum and absolute value maps, and the continuous mapping theorem that for all $T \geq 0$,
\begin{equation}\label{tozero1b}
\lim_{n \rightarrow \infty} P\bigg( 2 T \sup_{0 \leq s \leq T} \big| n^{-\frac{1}{2}} \big( \tilde{\ARR}^n(s) - \tilde{\REN}^n(s) \big) \big| \geq n^{\frac{1}{20}} \bigg) = 0.
\end{equation}
It follows from the well-known Doob's $L2$ inequality for continuous time martingales, see e.g. \cite{Cohen.15}, that (\ref{tozero2}) is at most $4 \theta n^{\frac{11}{20}} \big( n^{\frac{3}{10}} \big)^{-2} = 4 \theta n^{-\frac{1}{20}}$, and hence
\begin{equation}\label{tozero2b}
\lim_{n \rightarrow \infty} P\bigg( \sup_{0 \leq x \leq n^{\frac{11}{20}} } |\theta x - \CLOCK(\theta x)| \geq n^{\frac{3}{10}} \bigg) = 0.
\end{equation}
Combining the above, we conclude that for all $T \geq 0$,
$\lim_{n \rightarrow \infty} P\big( \sup_{0 \leq y \leq T} \big|\tilde{\ERR}^n(y)\big| > n^{-\frac{1}{5}} \big) = 0$.  The desired result then follows from the basic properties and definitions associated with weak convergence.  $\qed$
\endproof

\subsection{Proof of Lemma\ \ref{samedist1}.}
\proof{Proof of Lemma\ \ref{samedist1}:}
Let $M(t)$ denote the renewal function associated with renewal distribution $S$ evaluated at $t$, i.e. the expected number of renewals up to time $t$ in the corresponding ordinary renewal process.  It follows from a straightforward computation, the details of which we omit and which can be found in e.g. \cite{Ching.16}, that 
\begin{equation}\label{ord1}
M(t) = t + p (1-p)(\mu_1 - \mu_2)^2 (\mu_1 \mu_2)^{-2} \big(1 - \exp(- \mu_1 \mu_2 t ) \big).
\end{equation}  
Recall that for a generic r.v $X$, $Var[X]$ denotes the variance of $X$.  It then follows from (\ref{ord1}), the well-known representation for the variance of an equilibrium renewal process in terms of the renewal function (see e.g. \cite{Whitt.02} Theorem 7.2.4), and a straightforward calculation, that
\begin{eqnarray*}
Var[\REN^{\infty}(t)] &=& 2 \int_0^t \big( M(s) - s + \frac{1}{2} \big) ds
\\&=& 2 p (1-p)(\mu_1 - \mu_2)^2 (\mu_1 \mu_2)^{-2} \int_0^t \big(1 - \exp(- \mu_1 \mu_2 s ) \big) ds + t
\\&=& \big( 2 p (1-p)(\mu_1 - \mu_2)^2 (\mu_1 \mu_2)^{-2} + 1 \big) t - 2 p (1-p)(\mu_1 - \mu_2)^2 (\mu_1 \mu_2)^{-3} \big( 1 - \exp(- \mu_1 \mu_2 t) \big)
\\&=&\ \big( 2 p(1-p)(\mu_1 - \mu_2)^2 (\mu_1 \mu_2)^{-2} + 1 \big) t - 2 p(1-p)(\mu_1 - \mu_2)^2 (\mu_1 \mu_2)^{-3} 
\\&\ &\ \ \ \ \ + 2 p(1-p)(\mu_1 - \mu_2)^2 (\mu_1 \mu_2)^{-3} \exp(- \mu_1 \mu_2 t ).
\end{eqnarray*}
Furthermore, as $\REN^{\infty}$ has stationary increments (see e.g. \cite{GG13}), $\E[ \big( \REN^{\infty}(t-s) \big)^2 ] = \E[\big( \REN^{\infty}(s) \big)^2] + \E[\big( \REN^{\infty}(t) \big)^2] - 2 \E[ \REN^{\infty}(s) \REN^{\infty}(t)]$ for $0 \leq s \leq t$.
We conclude from our calculation of $Var[\REN^{\infty}(t)]$, after a straightforward calculation the details of which we omit, that for all $0 \leq s \leq t$, $\E[ \REN^{\infty}(s) \REN^{\infty}(t)]$ equals
\begin{eqnarray*}
\ &\ &\ \big( 2 p(1-p)(\mu_1 - \mu_2)^2 (\mu_1 \mu_2)^{-2} + 1 \big) s - p(1-p)(\mu_1 - \mu_2)^2 (\mu_1 \mu_2)^{-3}
\\&\ &\ \ \ \ \ \ +\ \ p(1-p)(\mu_1 - \mu_2)^2 (\mu_1 \mu_2)^{-3} \bigg( \exp\big(- \mu_1 \mu_2 s \big) + \exp\big(- \mu_1 \mu_2 t \big) - \exp\big(- \mu_1 \mu_2 (t-s) \big)\bigg).
\end{eqnarray*}
Combining the above with Lemma\ \ref{intouprop1}, definitions, a straightforward calculation the details of which we omit, and the fact that continuous Gaussian processes are uniquely determined by their mean and covariance functions completes the proof.  $\qed$
\endproof

\subsection{Proof of Lemma\ \ref{Qexp1}.}
\proof{Proof of Lemma\ \ref{Qexp1}:}
We begin by proving that for all $n \geq 0$, $J^n = \begin{bmatrix} (\mu_1 \mu_2)^n & 0 \\ \frac{\mu_1 - \mu_2}{\mu_1 \mu_2 - \theta} \big( (\mu_1 \mu_2)^n - \theta^n \big) & \theta^n \end{bmatrix}$, and proceed by induction.  The base case $n = 0$ is trivial, yielding the identity matrix.  Thus, suppose the induction holds for all $n \leq N$ for some $N \geq 0$.  Then
\begin{eqnarray*}
J^{N+1} &=& J^N \bullet J
\\&=&  \begin{bmatrix} (\mu_1 \mu_2)^N & 0 \\ \frac{\mu_1 - \mu_2}{\mu_1 \mu_2 - \theta} \big( (\mu_1 \mu_2)^N - \theta^N \big) & \theta^N \end{bmatrix} \bullet \begin{bmatrix} \mu_1 \mu_2 & 0 \\ \mu_1 - \mu_2 & \theta \end{bmatrix}
\\&=& \begin{bmatrix} (\mu_1 \mu_2)^{N+1} & 0 \\ \frac{\mu_1 - \mu_2}{\mu_1 \mu_2 - \theta} \big( (\mu_1 \mu_2)^N - \theta^N \big) (\mu_1 \mu_2) + \theta^N (\mu_1 - \mu_2) & \theta^{N+1} \end{bmatrix}
\\&=& \begin{bmatrix} (\mu_1 \mu_2)^{N+1} & 0 \\ \frac{\mu_1 - \mu_2}{\mu_1 \mu_2 - \theta} \big( (\mu_1 \mu_2)^{N+1} - \theta^{N+1} \big) & \theta^{N+1} \end{bmatrix},
\end{eqnarray*} 
completing the proof.  The desired lemma then follows from the definition of the matrix exponential and a straightforward calculation, the details of which we omit.  $\qed$
\endproof

\subsection{Proof of Lemma\ \ref{itorewrite1}.}
\proof{Proof of Lemma\ \ref{itorewrite1}:}
First, note that by a simple reindexing,
\begin{equation}\label{itoeq1}
\int_s^t \exp\big( - \theta(t - y) \big) \OU(y) dy = \exp\big(- \theta(t - s) \big) \int_0^{t-s} \exp(\theta y) \OU_{+,s}(y) dy.
\end{equation}
Applying the celebrated Ito's lemma (see \cite{Kar.12}), we conclude that
\begin{equation}\label{itoeq2}
\exp\big( \theta(t - s) \big) \OU_{+,s}(t - s) - \OU_{+,s}(0) = \int_0^{t - s} \theta \exp(\theta y) \OU_{+,s}(y) dy + \int_0^{t - s} \exp(\theta y) d\OU_{+,s}(y).
\end{equation}
Combining with (\ref{thesde}) and standard techniques from stochastic calculus, we conclude that $\int_0^{t - s} \exp(\theta y) \OU_{+,s}(y) dy$ equals
$$\theta^{-1} \Bigg( \exp\big( \theta(t - s) \big) \OU_{+,s}(t - s) - \OU_{+,s}(0) - \int_0^{t - s} \exp(\theta y) \bigg( - \mu_1 \mu_2 \OU_{+,s}(y) dy + d{\mathcal B}^3_s(y) \bigg) \Bigg),$$
from which it follows that $\int_0^{t - s} \exp(\theta y) \OU_{+,s}(y) dy$ equals
\begin{equation}\label{itoeq3}
(\theta - \mu_1 \mu_2)^{-1} \bigg( \exp\big( \theta(t - s) \big) \OU_{+,s}(t - s) - \OU_{+,s}(0) - \int_0^{t - s} \exp(\theta y) d{\mathcal B}^3_s(y) \bigg).
\end{equation}
Combining with (\ref{theeq2b}), we conclude that $\int_0^{t - s} \exp(\theta y) \OU_{+,s}(y) dy$ equals
\begin{eqnarray*}
\ &\ &\ (\theta - \mu_1 \mu_2)^{-1} \exp\big( \theta(t - s) \big) \bigg( \exp\big(-\mu_1 \mu_2 (t-s)\big) \OU_{+,s}(0) + \int_0^{t-s} \exp\big( - \mu_1 \mu_2(t - s - y) \big) d{\mathcal B}^3_s(y) \bigg)
\\&\ &\ \ \ \ \ \ \ - (\theta - \mu_1 \mu_2)^{-1} \OU_{+,s}(0) - (\theta - \mu_1 \mu_2)^{-1} \int_0^{t - s} \exp(\theta y) d{\mathcal B}^3_s(y).
\end{eqnarray*}
Combining with (\ref{itoeq1}), we find that 
$$\int_s^t \exp\big( - \theta(t - y) \big) \OU(y) dy = \int_0^{t-s} \exp\big( - \theta (t - s - y) \big) \OU_{+,s}(y) dy,$$ 
itself equal to
\begin{eqnarray*}
\ &\ &\ (\theta - \mu_1 \mu_2)^{-1} \bigg( \exp\big(-\mu_1 \mu_2 (t-s)\big) \OU_{+,s}(0) + \int_0^{t-s} \exp\big( - \mu_1 \mu_2(t - s - y) \big) d{\mathcal B}^3_s(y) \bigg)
\\&\ &\ \ \ \ \ \ \ - (\theta - \mu_1 \mu_2)^{-1} \exp\big( - \theta(t - s) \big) \OU_{+,s}(0) - (\theta - \mu_1 \mu_2)^{-1} \int_0^{t - s} \exp\big(-\theta(t - s - y)\big) d{\mathcal B}^3_s(y),
\\&=&\ (\theta - \mu_1 \mu_2)^{-1} \bigg( \exp\big(-\mu_1 \mu_2 (t-s)\big) - \exp\big( - \theta(t - s) \big) \bigg) \OU_{+,s}(0)
\\&\ &\ \ \ \ \ \ \ + (\theta - \mu_1 \mu_2)^{-1} \int_0^{t-s} \bigg( \exp\big( - \mu_1 \mu_2(t - s - y) \big) - \exp\big( - \theta (t - s - y) \big) \bigg) d{\mathcal B}^3_s(y).
\end{eqnarray*}
Combining with (\ref{theeq2}) and another simple reindexing completes the proof.  $\qed$
\endproof

\subsection{Proof of Lemma\ \ref{covlemma1}.}
\proof{Proof of Lemma\ \ref{covlemma1}:}
From definitions and independence, $\E\big[ \big( \overline{\mathcal G}(t) - \overline{\mathcal G}(s) \big)^2 \big]$ equals
\begin{eqnarray}
\ &\ &\ (\mu_1 - \mu_2)^2 \E\big[ \big( \int_s^t \exp\big( - \theta y \big) \OU(y) dy \big)^2 \big] \label{computevar1}
\\&\ &\ \ \ \ +\ \ \ \ 2 \E\big[ \big( \int_s^t \exp(- \theta y) d{\mathcal B}^1(y) \big)^2 \big]. \label{computevar2}
\end{eqnarray}
Furthermore,
\begin{eqnarray*}
\E\big[ \big( \int_s^t \exp\big( - \theta y \big) \OU(y) dy \big)^2 \big] &=& \E\big[ \int_s^t \exp\big( - \theta y_1 \big) \OU(y_1) dy_1 \times \int_s^t \exp\big( - \theta y_2 \big) \OU(y_2) dy_2 \big]
\\&=& \E\big[ \int_s^t \int_s^t \exp\big( - \theta (y_1 + y_2) \big) \OU(y_1) \OU(y_2) dy_1 dy_2 \big].
\end{eqnarray*}
It is easily verified from the basic properties of the OU-process, see e.g \cite{Doob.42}, that one may apply Fubini's theorem to the above expectation, and conclude (after a simple reindexing) that (\ref{computevar1}) equals
$$
2 (\mu_1 - \mu_2)^2 \int_s^t \int_s^{y_2} \exp\big( - \theta (y_1 + y_2) \big) \E\big[ \OU(y_1) \OU(y_2) \big] dy_1 dy_2,
$$
itself equal to (after substituting in the known covariance of $\OU$)
\begin{eqnarray*}
\ &\ &\ 2 p (1 - p) (\mu_1 - \mu_2)^2 (\mu_1 \mu_2)^{-1} \int_s^t \int_s^{y_2} \exp\big( - \theta (y_1 + y_2) \big) \exp\big( - \mu_1 \mu_2 (y_2 - y_1) \big) dy_1 dy_2 \nonumber
\\&\ &\ \ \ \ =\ \ \ \frac{ 2 p (1 - p) (\mu_1 - \mu_2)^2 }{ \mu_1 \mu_2 (\mu_1 \mu_2 - \theta) } \int_s^t \exp\big( -(\mu_1 \mu_2 + \theta) y_2 \big) \bigg( \exp\big( (\mu_1 \mu_2 - \theta) y_2 \big) - \exp\big( (\mu_1 \mu_2 - \theta) s \big) \bigg) dy_2
\\&\ &\ \ \ \ =\ \ \ \frac{ 2 p (1 - p) (\mu_1 - \mu_2)^2 }{ \mu_1 \mu_2 (\mu_1 \mu_2 - \theta) } \int_s^t \exp(- 2 \theta y_2) dy_2
\\&\ &\ \ \ \ \ \ \ \ \ -\ \ \frac{ 2 p (1 - p) (\mu_1 - \mu_2)^2 }{ \mu_1 \mu_2 (\mu_1 \mu_2 - \theta) } \exp\big( (\mu_1 \mu_2 - \theta) s \big) \int_s^t \exp\big( -(\mu_1 \mu_2 + \theta) y_2 \big) dy_2
\\&\ &\ \ \ \ =\ \ \ \frac{ 2 p (1 - p) (\mu_1 - \mu_2)^2 }{ 2 \theta \mu_1 \mu_2 (\mu_1 \mu_2 - \theta) } \exp(- 2 \theta s) - \frac{ 2 p (1 - p) (\mu_1 - \mu_2)^2 }{ 2 \theta \mu_1 \mu_2 (\mu_1 \mu_2 - \theta) } \exp(- 2 \theta t)
\\&\ &\ \ \ \ \ \ \ \ -\ \ \ \frac{ 2 p (1 - p) (\mu_1 - \mu_2)^2 }{ \mu_1 \mu_2 (\mu_1 \mu_2 - \theta)(\mu_1 \mu_2 + \theta) } \exp(- 2 \theta s)
\\&\ &\ \ \ \ \ \ \ \ +\ \ \ \frac{ 2 p (1 - p) (\mu_1 - \mu_2)^2 }{ \mu_1 \mu_2 (\mu_1 \mu_2 - \theta)(\mu_1 \mu_2 + \theta) } \exp\big( (\mu_1 \mu_2 - \theta) s \big) \exp\big( - (\mu_1 \mu_2 + \theta) t \big).
\end{eqnarray*}
Combining with some straightforward algebra, we conclude that (\ref{computevar1}) equals
\begin{eqnarray}
\ &\ &\ \frac{ p (1 - p) (\mu_1 - \mu_2)^2 }{ \theta \mu_1 \mu_2 (\mu_1 \mu_2 + \theta) } \exp(- 2 \theta s) - \frac{ p (1 - p) (\mu_1 - \mu_2)^2 }{ \theta \mu_1 \mu_2 (\mu_1 \mu_2 - \theta) } \exp(- 2 \theta t) \label{compute1varb}
\\&\ &\ \ \ \ \ \ \ \ +\ \ \ \frac{ 2 p (1 - p) (\mu_1 - \mu_2)^2 }{ \mu_1 \mu_2 (\mu_1 \mu_2 - \theta)(\mu_1 \mu_2 + \theta) } \exp\big( (\mu_1 \mu_2 - \theta) s \big) \exp\big( - (\mu_1 \mu_2 + \theta) t \big).\nonumber
\end{eqnarray}
We next analyze (\ref{computevar2}).  By a straightforward application of the Ito isometry, see e.g. \cite{Kar.12}, (\ref{computevar2}) equals $2 \int_s^t \exp(-2 \theta y) dy = \theta^{-1} \big( \exp(- 2 \theta s) - \exp(-2 \theta t) \big)$.  Combining with some straightforward algebra completes the proof.$\qed$ 
\endproof

\subsection{Proof of Corollary\ \ref{monovarcor}.}
\proof{Proof of Corollary\ \ref{monovarcor}:}
It follows from Lemma\ \ref{covlemma1} that 
\begin{eqnarray*}
\partial_t \E\big[ \big( \overline{\mathcal G}(t) - \overline{\mathcal G}(s) \big)^2 \big] &=& \bigg( \frac{ 2 p (1 - p) (\mu_1 - \mu_2)^2 }{ \mu_1 \mu_2 (\mu_1 \mu_2 - \theta) } + 2 \bigg) \exp(-2 \theta t)
\\&\ &\ \ \ \ \ \ -\ \ \frac{ 2 p (1 - p) (\mu_1 - \mu_2)^2 }{ \mu_1 \mu_2 (\mu_1 \mu_2 - \theta)} \exp\big( (\mu_1 \mu_2 - \theta) s \big) \exp\big( - (\mu_1 \mu_2 + \theta) t \big)
\\&\geq& \bigg( \frac{ 2 p (1 - p) (\mu_1 - \mu_2)^2 }{ \mu_1 \mu_2 (\mu_1 \mu_2 - \theta) } + 2 \bigg) \exp(-2 \theta t)
\\&\ &\ \ \ \ \ \ -\ \ \frac{ 2 p (1 - p) (\mu_1 - \mu_2)^2 }{ \mu_1 \mu_2 (\mu_1 \mu_2 - \theta)} \exp\big( (\mu_1 \mu_2 - \theta) t \big) \exp\big( - (\mu_1 \mu_2 + \theta) t \big)
\\&=& \Bigg( \bigg( \frac{ 2 p (1 - p) (\mu_1 - \mu_2)^2 }{ \mu_1 \mu_2 (\mu_1 \mu_2 - \theta) } + 2 \bigg) - \frac{ 2 p (1 - p) (\mu_1 - \mu_2)^2 }{ \mu_1 \mu_2 (\mu_1 \mu_2 - \theta)} \Bigg) \exp(-2 \theta t)
\\&=& 2 \exp(-2 \theta t)\ \ \ >\ \ \ 0.
\end{eqnarray*}
Similarly,
\begin{eqnarray*}
\partial_s \E\big[ \big( \overline{\mathcal G}(t) - \overline{\mathcal G}(s) \big)^2 \big] &=& - 2 \theta \bigg( \frac{ p (1 - p) (\mu_1 - \mu_2)^2 }{ \theta \mu_1 \mu_2 (\mu_1 \mu_2 + \theta) } + \theta^{-1} \bigg) \exp(- 2 \theta s)
\\&\ &\ \ \ \ \ \ \ +\ \ \frac{ 2 p (1 - p) (\mu_1 - \mu_2)^2 }{ \mu_1 \mu_2 (\mu_1 \mu_2 + \theta) } \exp\big( (\mu_1 \mu_2 - \theta) s \big) \exp\big( - (\mu_1 \mu_2 + \theta) t \big)
\\&\leq& - 2 \theta \bigg( \frac{ p (1 - p) (\mu_1 - \mu_2)^2 }{ \theta \mu_1 \mu_2 (\mu_1 \mu_2 + \theta) } + \theta^{-1} \bigg) \exp(- 2 \theta s)
\\&\ &\ \ \ \ \ \ \ +\ \ \frac{ 2 p (1 - p) (\mu_1 - \mu_2)^2 }{ \mu_1 \mu_2 (\mu_1 \mu_2 + \theta) } \exp\big( (\mu_1 \mu_2 - \theta) s \big) \exp\big( - (\mu_1 \mu_2 + \theta) s \big)
\\&=& \Bigg( - 2 \theta \bigg( \frac{ p (1 - p) (\mu_1 - \mu_2)^2 }{ \theta \mu_1 \mu_2 (\mu_1 \mu_2 + \theta) } + \theta^{-1} \bigg) 
+ \frac{ 2 p (1 - p) (\mu_1 - \mu_2)^2 }{ \mu_1 \mu_2 (\mu_1 \mu_2 + \theta) } \Bigg) \exp(- 2 \theta s)
\\&=& - 2 \exp(- 2 \theta s)\ \ \ <\ \ \ 0.
\end{eqnarray*}
Combining the above completes the proof.
$\qed$
\endproof

\subsection{Proof of Lemma\ \ref{diameter}.}
\proof{Proof of Lemma\ \ref{diameter}:}
In light of the monotonicities established in Corollary\ \ref{monovarcor}, the first two parts of the lemma would follow if we could prove that $\lim_{t \rightarrow \infty} \E\big[ \overline{\mathcal G}^2(t) \big] = \frac{\theta + \mu_1 + \mu_2 - 1}{\theta(\theta + \mu_1 \mu_2)}$.  A straightforward asymptotic analysis of Lemma\ \ref{covlemma1} demonstrates that
$$
\lim_{t \rightarrow \infty} \E\big[ \overline{\mathcal G}^2(t) \big] = \frac{ p (1 - p) (\mu_1 - \mu_2)^2 }{ \theta \mu_1 \mu_2 (\mu_1 \mu_2 + \theta) } + \theta^{-1}.$$
But note that since $\E[S] = 1$ implies that $\frac{p}{\mu_1} + \frac{1 - p}{\mu_2} = 1$, and thus 
\begin{equation}\label{p1mp}
p = \frac{\mu_1 (\mu_2 - 1)}{\mu_2 - \mu_1}\ \ \ \textrm{and}\ \ \ 1 - p = \frac{\mu_2(1 - \mu_1)}{\mu_2 - \mu_1},
\end{equation}
it follows that
\begin{eqnarray*}
\lim_{t \rightarrow \infty} \E\big[ \overline{\mathcal G}^2(t) \big] &=& \frac{ \mu_1 (\mu_2 - 1) \mu_2 (1 - \mu_1)}{ \theta \mu_1 \mu_2 (\mu_1 \mu_2 + \theta) } + \theta^{-1}
\\&=& \frac{ (\mu_2 - 1) (1 - \mu_1)}{\theta (\mu_1 \mu_2 + \theta) } + \frac{1}{\theta}\ \ \ =\ \ \ \frac{\theta + \mu_1 + \mu_2 - 1}{\theta(\theta + \mu_1 \mu_2)},
\end{eqnarray*}
completing the proof.  For the final part of the lemma, i.e. the explicit bound on $d_{\overline{\mathcal G}}(s,t)$, note that by Lemma\ \ref{covlemma1}, for $0 \leq s \leq t$, $\E\big[ \big( \overline{\mathcal G}(t) - \overline{\mathcal G}(s) \big)^2 \big]$ equals
\begin{eqnarray*}
\ &\ &\ \bigg( \frac{ p (1 - p) (\mu_1 - \mu_2)^2 }{ \theta \mu_1 \mu_2 (\mu_1 \mu_2 + \theta) } + \theta^{-1} \bigg) \exp(- 2 \theta s) - \bigg( \frac{ p (1 - p) (\mu_1 - \mu_2)^2 }{ \theta \mu_1 \mu_2 (\mu_1 \mu_2 - \theta) } + \theta^{-1} \bigg) \exp(- 2 \theta s) \exp\big(- 2 \theta (t - s) \big)
\\&\ &\ \ \ \ \ \ \ \ +\ \ \ \frac{ 2 p (1 - p) (\mu_1 - \mu_2)^2 }{ \mu_1 \mu_2 (\mu_1 \mu_2 - \theta)(\mu_1 \mu_2 + \theta) } \exp\big( (\mu_1 \mu_2 - \theta) s \big) \exp\big( - (\mu_1 \mu_2 + \theta) s \big) \exp\big( - (\mu_1 \mu_2 + \theta) (t-s) \big)
\\&\ &\ \ \ =\ \ \Bigg( \bigg( \frac{ p (1 - p) (\mu_1 - \mu_2)^2 }{ \theta \mu_1 \mu_2 (\mu_1 \mu_2 + \theta) } + \theta^{-1} \bigg) - \bigg( \frac{ p (1 - p) (\mu_1 - \mu_2)^2 }{ \theta \mu_1 \mu_2 (\mu_1 \mu_2 - \theta) } + \theta^{-1} \bigg)  \exp\big(- 2 \theta (t - s) \big) \Bigg) \exp(- 2 \theta s)
\\&\ &\ \ \ \ \ \ \ \ +\ \ \ \frac{ 2 p (1 - p) (\mu_1 - \mu_2)^2 }{ \mu_1 \mu_2 (\mu_1 \mu_2 - \theta)(\mu_1 \mu_2 + \theta) } \exp\big( (\mu_1 \mu_2 - \theta) s \big) \exp\big( - (\mu_1 \mu_2 + \theta) s \big) \exp\big( - (\mu_1 \mu_2 + \theta) (t-s) \big)
\\&\ &\ \ \ =\ \ \Bigg( \bigg( \frac{ p (1 - p) (\mu_1 - \mu_2)^2 }{ \theta \mu_1 \mu_2 (\mu_1 \mu_2 + \theta) } + \theta^{-1} \bigg) - \bigg( \frac{ p (1 - p) (\mu_1 - \mu_2)^2 }{ \theta \mu_1 \mu_2 (\mu_1 \mu_2 - \theta) } + \theta^{-1} \bigg)  \exp\big(- 2 \theta (t - s) \big) 
\\&\ &\ \ \ \ \ \ \ \ \ \ \ \ \ \ \ +\ \ \ \frac{ 2 p (1 - p) (\mu_1 - \mu_2)^2 }{ \mu_1 \mu_2 (\mu_1 \mu_2 - \theta)(\mu_1 \mu_2 + \theta) } \exp\big( - (\mu_1 \mu_2 + \theta) (t-s) \big) \Bigg) \exp(- 2 \theta s).
\end{eqnarray*}
Noting that
\begin{eqnarray*}
\ &\ &\ \bigg( \frac{ p (1 - p) (\mu_1 - \mu_2)^2 }{ \theta \mu_1 \mu_2 (\mu_1 \mu_2 + \theta) } + \theta^{-1} \bigg) - \bigg( \frac{ p (1 - p) (\mu_1 - \mu_2)^2 }{ \theta \mu_1 \mu_2 (\mu_1 \mu_2 - \theta) } + \theta^{-1} \bigg) + \frac{ 2 p (1 - p) (\mu_1 - \mu_2)^2 }{ \mu_1 \mu_2 (\mu_1 \mu_2 - \theta)(\mu_1 \mu_2 + \theta) }
\\&\ &\ \ \ =\ \ \ \frac{2 p (1 - p) (\mu_1 - \mu_2)^2}{\mu_1 \mu_2} \bigg( \frac{1}{2 \theta (\mu_1 \mu_2 + \theta)} - \frac{1}{2 \theta (\mu_1 \mu_2 - \theta)} + \frac{1}{(\mu_1 \mu_2 - \theta)(\mu_1 \mu_2 + \theta)} \bigg)
\\&\ &\ \ \ =\ \ \ \frac{2 p (1 - p) (\mu_1 - \mu_2)^2}{\mu_1 \mu_2} \frac{ (\mu_2 \mu_2 - \theta) - (\mu_2 \mu_2 + \theta) + 2 \theta }{ 2 \theta (\mu_1 \mu_2 - \theta)(\mu_1 \mu_2 + \theta)}\ \ \ =\ \ \ 0,
\end{eqnarray*}
we conclude that $\E\big[ \big( \overline{\mathcal G}(t) - \overline{\mathcal G}(s) \big)^2 \big]$ equals
\begin{eqnarray*}
\ &\ &\ \Bigg( \bigg(1 - \exp\big( - 2 \theta (t - s) \big) \bigg) \bigg( \frac{ p (1 - p) (\mu_1 - \mu_2)^2 }{ \theta \mu_1 \mu_2 (\mu_1 \mu_2 - \theta) } + \theta^{-1} \bigg) 
\\&\ &\ \ \ \ \ \ \ - \bigg(1 - \exp\big( - (\mu_1 \mu_2 + \theta)(t - s) \big) \bigg) \frac{ 2 p (1 - p) (\mu_1 - \mu_2)^2 }{ \mu_1 \mu_2 (\mu_1 \mu_2 - \theta)(\mu_1 \mu_2 + \theta) } \Bigg) \exp(- 2 \theta s),
\end{eqnarray*}
itself at most
$$8 \bigg( \frac{ p (1 - p) (\mu_1 - \mu_2)^2 }{ \mu_1 \mu_2 (\mu_1 \mu_2 - \theta) \theta } + \theta^{-1} \bigg) \max\bigg( 1 - \exp\big( - 2 \theta (t - s) \big) , 1 - \exp\big( - (\mu_1 \mu_2 + \theta)(t - s) \big) \bigg) \exp(- 2 \theta s).$$
Recalling that as our assumptions imply that $\mu_1 \mu_2 > \theta$, it follows that $2 \theta < \mu_1 \mu_2 + \theta$, and hence
$$
\E\big[ \big( \overline{\mathcal G}(t) - \overline{\mathcal G}(s) \big)^2 \big] \leq 8 \bigg( \frac{ p (1 - p) (\mu_1 - \mu_2)^2 }{ \mu_1 \mu_2 (\mu_1 \mu_2 - \theta) \theta } + \theta^{-1} \bigg) \bigg( 1 - \exp\big( - (\mu_1 \mu_2 + \theta)(t - s) \big) \bigg) \exp(- 2 \theta s).$$
Combining with definitions completes the proof.  $\qed$
\endproof

\subsection{Proof of Lemma\ \ref{notbad1}.}
\proof{Proof of Lemma\ \ref{notbad1}:}
Using the results of \cite{Dieker.13, DDG14}, it is easily verified that to prove the lemma, it suffices to demonstrate that for all $x > 0$, $P\big( W^{\infty}(\infty) > x \big) > 0$.  To establish this fact, we will appeal to a straightforward infinite-server lower bound.  Indeed, for $n \geq 1$, let  $\ddot{Q}^n$ denote the $M/G/\infty$ queue in which arrivals come as a Poisson process with rate $\lambda_n$, services are distributed as $S$, and the system is initially empty; and $\lbrace \ddot{W}^n(t), t \geq 0 \rbrace$ the corresponding total number of jobs in system, with $\ddot{W}^n(\infty)$ a r.v. with the corresponding steady-state distribution (simply a Poisson r.v).  Let $\dot{Q}^n$ denote the $M/G/n + M$ queue in which arrivals come as a Poisson process with rate $\lambda_n$, services are distributed as $S$, the abandonment rate is $\theta$, and the system is initially empty; with $\lbrace \dot{W}^n(t), t \geq 0 \rbrace$ the corresponding total number of jobs in system, and $\dot{W}^n(\infty)$ a r.v. with the corresponding steady-state distribution (see e.g. \cite{DDG14} for a discussion of ergodicity for $M/G/n + M$ queues).  Recall that in an $M/G/n + M$ queue, the virtual waiting time of an arriving job J is the time that job would have to wait in system until it began service, if that job was infinitely patient, and is a measurable function of the inter-arrival, service, and patience times of those jobs which arrived to the system prior to $J$, as well as the inter-arrival time of job $J$.  Our arguments will rely implicitly on the well-known virtual-waiting-time representation for the time-in-system of an arriving job to a $M/G/n + M$ queue, and we refer the intereseted reader to e.g. \cite{Tal.09,Glynn.05} for further details.  We now provide an explicit coupling / construction for $\ddot{Q}^n$ and $\dot{Q}^n$, demonstrating that one may construct both queueing systems on a common probability space s.t. both systems see the same arrival stream, with common arrivals having the same type, in such a way that w.p.1 every job $J$ departs later (i.e. no sooner) from $\dot{Q}^n$ than from $\ddot{Q}^n$.  Intuitively, this follows from the simple fact that (supposing the arrival processess are coupled in the two systems), the infinitesimal rate at which a type-z job J leaves the system (either due to abandonment or service completion) in $\dot{Q}^n$ is either $\theta$ or $\mu_z$, while in $\ddot{Q}^n$ this rate always equals $\mu_z$, in which case the desired dominance follows from Assumption\ \ref{A1}.  
\\\indent More formally, let $P_{1,\lambda}$ be a rate-$(p \lambda_n)$ Poisson process, and $P_{2,\lambda}$ be a rate-$\big( (1-p)\lambda_n)$ Poisson process.  For $z \in \lbrace 1,2 \rbrace$ and $k \geq 1$, let $P_{z,k}$ be a rate $\mu_z$ Poisson process.  For $z \in \lbrace 1,2 \rbrace$ and $k \geq 1$, let $B_{z,k}$ be a geometrically distributed r.v., with support on the strictly positive integers, s.t. $P(B_{z,k} = 1) = \frac{\theta}{\mu_z}$; $E_{z,k}$ be an exponentially distributed r.v. with rate $\theta$; and $S_{z,k}$ be an exponentially distributed r.v. with rate $\mu_z$.  Suppose that all these processes and r.v are constructed together on a common probability space s.t. they are mutually independent.  We proceed to generate the inter-arrival, service, and patience times for the arrivals to both systems as follows.  For $z \in \lbrace 1,2 \rbrace$, both queueing systems witness a type-z arrival at each jump in $P_{z,\lambda}$, and see no other arrivals.  We determine the service and patience times of the arrivals to both systems inductively.  
\\\indent For an arriving job $J$, let us denote the patience time of $J$ in $\dot{Q}^n$ by $\dot{\PAT}_J$, service time in $\dot{Q}^n$ by $\dot{S}_J$, and service time in $\ddot{Q}^n$ by $\ddot{S}_J$.  Also, let us denote the virtual waiting time (in $\dot{Q}^n$) of job $J$ by $\dot{\VRW}_J$ (our inductive definition will ensure this quantity is well-defined, i.e. not circular).  In addition, for a job $J$, let $z_J$ denote the type of $J$.  Also, if job $J$ is the $k$th type-$z_J$ arrival to the system, we set $k_J = k$.  Thus (for example) if job $J$ is type 1 and there have been exactly 3 previous type-1 arrivals to the system, we set $z_J = 1$ and $k_J = 4$.  For the job $J$ which is the first arrival to the system, note that $\dot{\VRW}_J = 0$, as both systems are initially empty, and $k_J = 1$.  Set both $\dot{S}_J$ and $\ddot{S}_J$ equal to the time of the first jump in $P_{z_J,k_J}$, and set $\dot{\PAT}_J$ equal to $E_{z_J,k_J}$.  Now, suppose for induction that we have constructed the patience and service times for the first $i$ arrivals to the system.  Let job $J$ be the $(i+1)$st arrival to the system.  By our description of the arrival processes, our inductive hypothesis, and the measurability properties of the virtual waiting time, $\dot{\VRW}_J$ has thus already been defined (implicitly by the basic definitions associated with the relevant queueing systems).  Set $\ddot{S}_J$ equal to the time of the first jump in $P_{z_J,k_J}$.  If the time of the $B_{z_J,k_J}$-th jump in $P_{z_J,k_J}$ occurs before time $\ddot{W}_J$, set $\dot{\PAT}_J$ equal to the time of this jump, and set $\dot{S}_J$ equal to $S_{z_J,k_J}$.  Otherwise, set $\dot{\PAT}_J$ equal to $\dot{\VRW}_J + E_{z_J,k_J}$, and set $\dot{S}_J$ equal to the difference between the time of the first jump in $P_{z_J,k_J}$ which occurs after time $\dot{\VRW}_J$, and $\dot{\VRW}_J$ itself.  
\\\indent That this construction is correct, i.e. yielding sequences of inter-arrival, patience, and service times with the correct joint distributions, is straightforward and follows from the basic properties of Poisson processes and the relevant queueing systems, and we omit the details.  It follows that we may construct both systems on a common probability space s.t. every job arrives to both systems at the same time, but departs later in $\dot{Q}^n$.  The well-known ergodicity properties of both systems then imply that for all $x \geq 0$, $P\big( \dot{W}^n(\infty) > x \big) \geq P\big( \ddot{W}^n(\infty) > x \big)$.  Rescaling both systems, and applying the well-known convergence of the appropriately scaled sequence of Poisson r.v.s (i.e. $\lbrace \ddot{W}^n(\infty), n \geq 1 \rbrace$) to a limiting normal distribution then completes the proof, and we refer the interested reader to \cite{GG13} for further details regarding the relevant limiting arguments.  $\qed$
\endproof

\subsection{Proof of Lemma\ \ref{itorewrite2}.}
\proof{Proof of Lemma\ \ref{itorewrite2}:}
It follows from an argument nearly identical to that used in the proof of Lemma\ \ref{itorewrite1}, the details of which we omit, that 
$\int_0^t \exp\big( - \theta(t - y) \big) \hat{Y}^1_{\Upsilon^*}(y) dy$ equals
\begin{eqnarray*}
\ &\ &\ \ (\theta - \mu_1 \mu_2)^{-1} \bigg( \exp\big(-\mu_1 \mu_2 t\big) - \exp\big( - \theta t \big) \bigg) \hat{Y}^1_{\Upsilon^*}(0)
\\&\ &\ \ \ \ \ \ \ + (\theta - \mu_1 \mu_2)^{-1} \int_0^t \bigg( \exp\big( - \mu_1 \mu_2(t - y) \big) - \exp\big( - \theta (t - y) \big) \bigg) d{\mathcal B}^3(y).
\end{eqnarray*}
Combining with definitions, we conclude that w.p.1, for all $t \geq 0$,
\begin{eqnarray*}
\hat{Y}^2_{\Upsilon^*}(t) &=& - (\mu_1 - \mu_2) \int_0^t \exp\big( - \theta(t - y) \big) \hat{Y}^1_{\Upsilon^*}(y) dy + \exp(- \theta t) \hat{Y}^2_{\Upsilon^*}(0)
\\&\ &\ \ \ \ \ \ \ + \int_0^t \exp\big( - \theta(t - y) \big) d{\mathcal B}^2(y) + \frac{B}{\theta}\bigg(1 - \exp\big(- \theta t \big) \bigg),
\end{eqnarray*}
from which it follows (by Lemma\ \ref{chainit1}) that 
$$\hat{Y}^2_{\Upsilon^*}(t) \geq - (\mu_1 - \mu_2) \int_0^t \exp\big( - \theta(t - y) \big) \hat{Y}^1_{\Upsilon^*}(y) dy + \int_0^t \exp\big( - \theta(t - y) \big) d{\mathcal B}^2(y) + \exp( - \theta t ) C_1 - \frac{|B|}{\theta}.
$$  
Noting from definitions that w.p.1, for all $y \geq 0$,
$$- (\mu_1 - \mu_2)\hat{Y}^1_{\Upsilon^*}(y) \geq - (\mu_1 + \mu_2) \big( |C_3| + |C_4| \big) \exp(- \mu_1 \mu_2 y) - (\mu_1 - \mu_2) \int_0^y \exp\big( - \mu_1 \mu_2 (y - s) \big) d{\mathcal B}^3(s),$$
we further conclude after some straightforward algebra that w.p.1, for all $t \geq 0$, $\hat{Y}^2_{\Upsilon^*}(t)$ is at least
\begin{eqnarray*}
\ &\ &\ - ( \mu_1 - \mu_2 ) \int_0^t \exp\big( - \theta(t - y) \big) \int_0^y \exp\big( - \mu_1 \mu_2 (y - s) \big) d{\mathcal B}^3(s)
\\&\ &\ \ \ \ \ \ + \int_0^t \exp\big( - \theta(t - y) \big) d{\mathcal B}^2(y) + \exp( - \theta t ) C_1 - (\mu_1 + \mu_2) \big( |C_3| + |C_4| \big) t \exp(- \theta t) - \frac{|B|}{\theta}.
\end{eqnarray*}
The lemma then follows from the well-known stochastic integral representation for the OU process, and the symmetries of the centered OU process, and we refer the reader to \cite{Ross.09} for details.  $\qed$
\endproof

\subsection{Proof of Lemma\ \ref{covlemma1lower}.}
\proof{Proof of Lemma\ \ref{covlemma1lower}:}
We first compute $\E\big[ \Pi(s) \Pi(t) \big]$.  Let us define
$$C(s,t) \stackrel{\Delta}{=} \E\bigg[ \int_0^t \exp( \theta y_2 ) \OU_{0,a}(y_2) dy_2 \times \int_0^s \exp( \theta y_1 ) \OU_{0,a}(y_1) dy_1 \bigg].$$
First, note that by applying Fubini's theorem as in the proof of Lemma\ \ref{covlemma1}, and substituting in the known covariance of the OU process, we conclude that $C(s,t)$ equals
\begin{eqnarray}
\ &\ &\ \E\bigg[ \int_0^t \int_0^s \exp\big( \theta (y_1 + y_2) \big) \OU_{0,a}(y_1) \OU_{0,a}(y_2) dy_1  dy_2 \bigg] \nonumber
\\&\ &\ \ \ \ \ \ =\ \ \ \int_0^t \int_0^s \exp\big( \theta (y_1 + y_2) \big) \E\big[ \OU_{0,a}(y_1) \OU_{0,a}(y_2) \big] dy_1  dy_2 \nonumber
\\&\ &\ \ \ \ \ \ =\ \ \ 2 \int_0^s \int_0^{y_2} \exp\big( \theta (y_1 + y_2) \big) \E\big[ \OU_{0,a}(y_1) \OU_{0,a}(y_2) \big] dy_1  dy_2 \label{exp1lower}
\\&\ &\ \ \ \ \ \ \ \ \ \ \ \ \ + \int_s^t \int_0^s \exp\big( \theta (y_1 + y_2) \big) \E\big[ \OU_{0,a}(y_1) \OU_{0,a}(y_2) \big] dy_1  dy_2. \label{exp2lower}
\end{eqnarray}
(\ref{exp1lower}) equals
\begin{eqnarray*}
\ &\ &\ 2 \frac{p (1 - p)}{\mu_1 \mu_2} \int_0^s \int_0^{y_2} \exp\big( \theta (y_1 + y_2) \big) \bigg( \exp\big( - \mu_1 \mu_2 (y_2 - y_1) \big) - \exp\big( - \mu_1 \mu_2 (y_2 + y_1) \big) \bigg) dy_1  dy_2
\\&\ &\ \ =\ \ \ 2 \frac{p (1 - p)}{\mu_1 \mu_2} \int_0^s \int_0^{y_2} \exp\big( \theta (y_1 + y_2) \big) \exp\big( - \mu_1 \mu_2 (y_2 - y_1) \big) dy_1 dy_2
\\&\ &\ \ \ \ \ \ -\ \ \ 2 \frac{p (1 - p)}{\mu_1 \mu_2} \int_0^s \int_0^{y_2} \exp\big( \theta (y_1 + y_2) \big) \exp\big( - \mu_1 \mu_2 (y_2 + y_1) \big) dy_1 dy_2
\\&\ &\ \ =\ \ \ 2 \frac{p (1 - p)}{\mu_1 \mu_2} \int_0^s \exp\big( (\theta - \mu_1 \mu_2) y_2 \big) \int_0^{y_2} \exp\big( (\theta + \mu_1 \mu_2) y_1 \big) dy_1 dy_2
\\&\ &\ \ \ \ \ \ -\ \ \ 2 \frac{p (1 - p)}{\mu_1 \mu_2} \int_0^s \exp\big( (\theta - \mu_1 \mu_2) y_2 \big) \int_0^{y_2} \exp\big( (\theta - \mu_1 \mu_2) y_1 \big) dy_1 dy_2
\\&\ &\ \ =\ \ \ 2 \frac{p (1 - p)}{\mu_1 \mu_2 (\theta + \mu_1 \mu_2)} \int_0^s \exp\big( (\theta - \mu_1 \mu_2) y_2 \big) \bigg( \exp\big( (\theta + \mu_1 \mu_2) y_2 \big) - 1 \bigg) dy_2
\\&\ &\ \ \ \ \ \ -\ \ \ 2 \frac{p (1 - p)}{\mu_1 \mu_2 (\mu_1 \mu_2 - \theta)} \int_0^s \exp\big( (\theta - \mu_1 \mu_2) y_2 \big) \bigg( 1 - \exp\big( - (\mu_1 \mu_2 - \theta) y_2 \big) \bigg) dy_2
\\&\ &\ \ =\ \ \ 2 \frac{p (1 - p)}{\mu_1 \mu_2 (\theta + \mu_1 \mu_2)} \int_0^s \exp( 2 \theta y_2) dy_2
\\&\ &\ \ \ \ \ \ \ -\ \ \ 2 \frac{p (1 - p)}{\mu_1 \mu_2 (\theta + \mu_1 \mu_2)} \int_0^s \exp\big( - (\mu_1 \mu_2 - \theta) y_2 \big) dy_2
\\&\ &\ \ \ \ \ \ \ -\ \ \ 2 \frac{p (1 - p)}{\mu_1 \mu_2 (\mu_1 \mu_2 - \theta)} \int_0^s \exp\big( - (\mu_1 \mu_2 - \theta) y_2 \big) dy_2
\\&\ &\ \ \ \ \ \ \ +\ \ \ 2 \frac{p (1 - p)}{\mu_1 \mu_2 (\mu_1 \mu_2 - \theta)} \int_0^s \exp\big( - 2 (\mu_1 \mu_2 - \theta) y_2 \big) dy_2
\\&\ &\ \ =\ \ \ \frac{p (1 - p)}{\mu_1 \mu_2 (\theta + \mu_1 \mu_2) \theta} \big( \exp(2 \theta s) - 1 \big)
\\&\ &\ \ \ \ \ \ \ -\ \ \ 2 \frac{p (1 - p)}{\mu_1 \mu_2 (\theta + \mu_1 \mu_2) (\mu_1 \mu_2 - \theta)} \bigg( 1 - \exp\big( - (\mu_1 \mu_2 - \theta) s \big) \bigg)
\\&\ &\ \ \ \ \ \ \ -\ \ \ 2 \frac{p (1 - p)}{\mu_1 \mu_2 (\mu_1 \mu_2 - \theta)^2} \bigg( 1 - \exp\big( - (\mu_1 \mu_2 - \theta) s \big) \bigg)
\\&\ &\ \ \ \ \ \ \ +\ \ \ \frac{p (1 - p)}{\mu_1 \mu_2 (\mu_1 \mu_2 - \theta)^2} \bigg( 1 - \exp\big( - 2 (\mu_1 \mu_2 - \theta) s \big) \bigg)
\\&\ &\ \ =\ \ \ \frac{p (1 - p) }{\mu_1 \mu_2 (\mu_1 \mu_2 + \theta) \theta } \exp(2 \theta s)
\\&\ &\ \ \ \ \ \ \ +\ \ \ \frac{4 p (1 - p)}{(\mu_1 \mu_2 + \theta)(\mu_1 \mu_2 - \theta)^2} \exp\big( - (\mu_1 \mu_2 - \theta) s \big)
\\&\ &\ \ \ \ \ \ \ -\ \ \ \frac{p (1 - p)}{\mu_1 \mu_2 (\mu_1 \mu_2 - \theta)^2} \exp\big( - 2 (\mu_1 \mu_2 - \theta) s \big)
\\&\ &\ \ \ \ \ \ \ -\ \ \ \frac{p(1 - p)}{\theta (\mu_1 \mu_2 - \theta)^2}.
\end{eqnarray*}
Similarly, (\ref{exp2lower}) equals
\begin{eqnarray*}
\ &\ &\ \frac{p (1 - p)}{\mu_1 \mu_2} \int_s^t \int_0^s \exp\big( \theta (y_1 + y_2) \big) \bigg( \exp\big( - \mu_1 \mu_2 (y_2 - y_1) \big) - \exp\big( - \mu_1 \mu_2 (y_2 + y_1) \big) \bigg) dy_1  dy_2
\\&\ &\ \ =\ \ \ \frac{p (1 - p)}{\mu_1 \mu_2} \int_s^t \int_0^s \exp\big( \theta (y_1 + y_2) \big) \exp\big( - \mu_1 \mu_2 (y_2 - y_1) \big) dy_1 dy_2
\\&\ &\ \ \ \ \ \ -\ \ \ \frac{p (1 - p)}{\mu_1 \mu_2} \int_s^t \int_0^s \exp\big( \theta (y_1 + y_2) \big) \exp\big( - \mu_1 \mu_2 (y_2 + y_1) \big) dy_1 dy_2
\\&\ &\ \ =\ \ \ \frac{p (1 - p)}{\mu_1 \mu_2} \int_s^t \exp\big( (\theta - \mu_1 \mu_2) y_2 \big) \int_0^s \exp\big( (\theta + \mu_1 \mu_2) y_1 \big) dy_1 dy_2
\\&\ &\ \ \ \ \ \ -\ \ \ \frac{p (1 - p)}{\mu_1 \mu_2} \int_s^t \exp\big( (\theta - \mu_1 \mu_2) y_2 \big) \int_0^s \exp\big( (\theta - \mu_1 \mu_2) y_1 \big) dy_1 dy_2
\\&\ &\ \ =\ \ \ \frac{p (1 - p)}{\mu_1 \mu_2 (\mu_1 \mu_2 + \theta)} \int_s^t \exp\big( (\theta - \mu_1 \mu_2) y_2 \big)\bigg( \exp\big( (\mu_1 \mu_2 + \theta) s \big) - 1 \bigg) dy_2
\\&\ &\ \ \ \ \ \ -\ \ \ \frac{p (1 - p)}{\mu_1 \mu_2 (\mu_1 \mu_2 - \theta)} \int_s^t \exp\big( (\theta - \mu_1 \mu_2) y_2 \big) \bigg( 1 - \exp\big( - (\mu_1 \mu_2 - \theta) s \big) \bigg) dy_2
\\&\ &\ \ =\ \ \ \frac{p (1 - p)}{\mu_1 \mu_2 (\mu_1 \mu_2 + \theta) (\mu_1 \mu_2 - \theta)} \bigg( \exp\big( (\mu_1 \mu_2 + \theta) s \big) - 1 \bigg) \bigg( \exp\big( - (\mu_1 \mu_2 - \theta) s \big) - \exp\big( - (\mu_1 \mu_2 - \theta) t \big) \bigg)
\\&\ &\ \ \ \ \ \ -\ \ \ \frac{p (1 - p)}{\mu_1 \mu_2 (\mu_1 \mu_2 - \theta)^2} \bigg( 1 - \exp\big( - (\mu_1 \mu_2 - \theta) s \big) \bigg) \bigg( \exp\big( - (\mu_1 \mu_2 - \theta) s \big) - \exp\big( - (\mu_1 \mu_2 - \theta) t \big) \bigg)
\\&\ &\ \ =\ \ \ \frac{p (1 - p)}{\mu_1 \mu_2 (\mu_1 \mu_2 + \theta) (\mu_1 \mu_2 - \theta)} \exp(2 \theta s) 
\\&\ &\ \ \ \ \ \ -\ \ \ \frac{p (1 - p)}{\mu_1 \mu_2 (\mu_1 \mu_2 + \theta) (\mu_1 \mu_2 - \theta)} \exp\big( (\mu_1 \mu_2 + \theta) s - (\mu_1 \mu_2 - \theta) t \big)
\\&\ &\ \ \ \ \ \ -\ \ \ \frac{p (1 - p)}{\mu_1 \mu_2 (\mu_1 \mu_2 + \theta) (\mu_1 \mu_2 - \theta)} \exp\big( - (\mu_1 \mu_2 - \theta) s \big)
\\&\ &\ \ \ \ \ \ +\ \ \ \frac{p (1 - p)}{\mu_1 \mu_2 (\mu_1 \mu_2 + \theta) (\mu_1 \mu_2 - \theta)} \exp\big( - (\mu_1 \mu_2 - \theta) t \big)
\\&\ &\ \ \ \ \ \ -\ \ \ \frac{p (1 - p)}{\mu_1 \mu_2 (\mu_1 \mu_2 - \theta)^2} \exp\big( - (\mu_1 \mu_2 - \theta) s \big)
\\&\ &\ \ \ \ \ \ +\ \ \ \frac{p (1 - p)}{\mu_1 \mu_2 (\mu_1 \mu_2 - \theta)^2} \exp\big( - (\mu_1 \mu_2 - \theta) t \big)
\\&\ &\ \ \ \ \ \ +\ \ \ \frac{p (1 - p)}{\mu_1 \mu_2 (\mu_1 \mu_2 - \theta)^2} \exp\big( - 2 (\mu_1 \mu_2 - \theta) s \big)
\\&\ &\ \ \ \ \ \ -\ \ \ \frac{p (1 - p)}{\mu_1 \mu_2 (\mu_1 \mu_2 - \theta)^2} \exp\big( - (\mu_1 \mu_2 - \theta) (s + t) \big).
\end{eqnarray*}
We conclude, after some straightforward algebra, that $C(s,t)$ equals
\begin{eqnarray}
\ &\ &\ \frac{p (1 - p)}{(\mu_1 \mu_2 + \theta)(\mu_1 \mu_2 - \theta) \theta} \exp( 2 \theta s) \label{thecovlower0}
\\&\ &\ \ \ \ \ \ +\ \ \ \frac{2 p (1 - p)}{(\mu_1 \mu_2 + \theta)(\mu_1 \mu_2 - \theta)^2} \exp\big( - (\mu_1 \mu_2 - \theta) s \big) \nonumber
\\&\ &\ \ \ \ \ \ -\ \ \ \frac{p(1 - p)}{\theta (\mu_1 \mu_2 - \theta)^2} \nonumber
\\&\ &\ \ \ \ \ \ -\ \ \ \frac{p (1 - p)}{\mu_1 \mu_2 (\mu_1 \mu_2 + \theta) (\mu_1 \mu_2 - \theta)} \exp\big( (\mu_1 \mu_2 + \theta) s - (\mu_1 \mu_2 - \theta) t \big) \nonumber
\\&\ &\ \ \ \ \ \ +\ \ \ \frac{2 p (1 - p)}{(\mu_1 \mu_2 + \theta)(\mu_1 \mu_2 - \theta)^2} \exp\big( - (\mu_1 \mu_2 - \theta) t \big) \nonumber
\\&\ &\ \ \ \ \ \ -\ \ \ \frac{p (1 - p)}{\mu_1 \mu_2 (\mu_1 \mu_2 - \theta)^2} \exp\big( - (\mu_1 \mu_2 - \theta) (s + t) \big). \nonumber
\end{eqnarray}
It then follows from definitions and the known covariance of $\OU_{0,b}$, and the independence of $\OU_{0,a}$ and $\OU_{0,b}$, that for all $0 \leq s \leq t$, $\E\big[ \Pi(s) \Pi(t) \big]$ equals
\begin{eqnarray}
\ &\ &\ \bigg( \frac{(\mu_1 - \mu_2)^2 p (1 - p)}{(\mu_1 \mu_2 + \theta)(\mu_1 \mu_2 - \theta) \theta} + \theta^{-1} \bigg) \exp\big( - \theta (t-s) \big) \label{thecovlower}
\\&\ &\ \ \ \ \ \ +\ \ \ \frac{(\mu_1 - \mu_2)^2 2 p (1 - p)}{(\mu_1 \mu_2 + \theta)(\mu_1 \mu_2 - \theta)^2} \exp\big( - \mu_1 \mu_2 s - \theta t \big) \nonumber
\\&\ &\ \ \ \ \ \ -\ \ \ \big( \frac{(\mu_1 - \mu_2)^2 p(1 - p)}{\theta (\mu_1 \mu_2 - \theta)^2} + \theta^{-1} \big) \exp\big( - \theta(s + t) \big) \nonumber
\\&\ &\ \ \ \ \ \ -\ \ \ \frac{(\mu_1 - \mu_2)^2 p (1 - p)}{\mu_1 \mu_2 (\mu_1 \mu_2 + \theta) (\mu_1 \mu_2 - \theta)} \exp\big( - \mu_1 \mu_2 (t - s) \big) \nonumber
\\&\ &\ \ \ \ \ \ +\ \ \ \frac{(\mu_1 - \mu_2)^2 2 p (1 - p)}{(\mu_1 \mu_2 + \theta)(\mu_1 \mu_2 - \theta)^2} \exp\big( - \theta s - \mu_1 \mu_2 t \big) \nonumber
\\&\ &\ \ \ \ \ \ -\ \ \ \frac{(\mu_1 - \mu_2)^2 p (1 - p)}{\mu_1 \mu_2 (\mu_1 \mu_2 - \theta)^2} \exp\big( - \mu_1 \mu_2 (s + t) \big). \nonumber
\end{eqnarray}
Next, let us prove that $\E\big[ \Pi(s) \Pi(t) \big] \geq 0$.  In light of the independence of $\OU_{0,a}$ and $\OU_{0,b}$ and known non-negative covariance of $\OU_{0,b},$ it suffices to prove that $C(s,t)$ is non-negative.  It follows from definitions that $C(s,s) \geq 0$ for all $s \geq 0$.  Thus it suffices to prove that for all $s \geq 0$, $\frac{d}{dt} C(s,t) \geq 0$ on $[s,\infty)$.
Thus let us fix an $s \geq 0$.  Then from (\ref{thecovlower0}), we compute that $\frac{d}{dt} C(s,t)$ equals 
\begin{eqnarray*}
\ &\ &\ -( \mu_1 \mu_2 - \theta) \frac{p (1 - p)}{\mu_1 \mu_2 - \theta} \bigg(
\frac{2}{(\mu_1 \mu_2 + \theta) (\mu_1 \mu_2 - \theta)} - \frac{1}{\mu_1 \mu_2 (\mu_1 \mu_2 + \theta)} \exp\big( (\mu_1 \mu_2 + \theta) s \big) 
\\&\ &\ \ \ \ \ \ \ \ \ \ \ \ \ \ \ \ \ \ \ \ \ \ \ \ \ \ \ \ \ \ - \frac{1}{\mu_1 \mu_2 (\mu_1 \mu_2 - \theta)} \exp\big( -(\mu_1 \mu_2 - \theta) s \big) \bigg)
\exp\big( -(\mu_1 \mu_2 - \theta) t \big).
\end{eqnarray*}
Thus it suffices to prove that
$$\ddot{g}(s) \stackrel{\Delta}{=}
\frac{2}{(\mu_1 \mu_2 + \theta) (\mu_1 \mu_2 - \theta)} - \frac{1}{\mu_1 \mu_2 (\mu_1 \mu_2 + \theta)} \exp\big( (\mu_1 \mu_2 + \theta) s \big) - \frac{1}{\mu_1 \mu_2 (\mu_1 \mu_2 - \theta)} \exp\big( -(\mu_1 \mu_2 - \theta) s \big)$$
is non-positive for all $s \geq 0$.  Noting that a straightforward calcualtion shows that $\ddot{g}(0) = 0$, we now accomplish this by proving that $\frac{d}{ds}\ddot{g}(s) \leq 0$ for all $s \geq 0$.  Indeed, a straightforward calculation shows that
$$
\frac{d}{ds}\ddot{g}(s) = (\mu_1 \mu_2)^{-1} \bigg( \exp\big( -(\mu_1 \mu_2 - \theta) s \big) - \exp\big( (\mu_1 \mu_2 + \theta) s \big) \bigg),$$
which is non-positive since the fact that $\mu_1 \mu_2 > \theta$ implies that $\exp\big( -(\mu_1 \mu_2 - \theta) s \big) \leq \exp\big( (\mu_1 \mu_2 + \theta) s \big)$ for all $s \geq 0$.  Combining the above then completes the proof that $\E\big[ \Pi(s) \Pi(t) \big] \geq 0$ for all $0 \leq s \leq t$.
\\\indent Finally, let us compute $\lim_{t \rightarrow \infty} \E\big[ \Pi^2(t) \big]$.  It follows from (\ref{thecovlower}), a straightforward asymptotic analysis, and (\ref{p1mp}) that the desired limit equals
\begin{eqnarray*}
\ &\ &\ \bigg( \frac{(\mu_1 - \mu_2)^2 p (1 - p)}{(\mu_1 \mu_2 + \theta)(\mu_1 \mu_2 - \theta) \theta} + \theta^{-1} \bigg) - \frac{(\mu_1 - \mu_2)^2 p (1 - p)}{\mu_1 \mu_2 (\mu_1 \mu_2 + \theta) (\mu_1 \mu_2 - \theta)}
\\&\ &\ \ \ =\ \ \ \frac{(\mu_1 - \mu_2)^2 p (1 - p)}{\mu_1 \mu_2 (\mu_1 \mu_2 + \theta)(\mu_1 \mu_2 - \theta) \theta}(\mu_1 \mu_2 - \theta) + \frac{1}{\theta}
\\&\ &\ \ \ =\ \ \ \frac{(\mu_1 - \mu_2)^2 p (1 - p) + \mu_1 \mu_2 (\mu_1 \mu_2 + \theta)}{\mu_1 \mu_2 (\mu_1 \mu_2 + \theta) \theta}
\\&\ &\ \ \ =\ \ \ \frac{\mu_1 (\mu_2 - 1)\mu_2(1 - \mu_1) + \mu_1 \mu_2 (\mu_1 \mu_2 + \theta)}{\mu_1 \mu_2 (\mu_1 \mu_2 + \theta) \theta}
\\&\ &\ \ \ =\ \ \ \frac{\theta + \mu_1 + \mu_2 - 1}{\theta (\theta + \mu_1 \mu_2)},
\end{eqnarray*}
completing the proof.  $\qed$
\endproof

\subsection{Proof of Lemma\ \ref{isposlower1}.}
\proof{Proof of Lemma\ \ref{isposlower1}:}
Let $t_1 \stackrel{\Delta}{=} \frac{1}{4 \theta} \min\big( \frac{|B|}{(\mu_1 + \mu_2)(|C_3| + |C_4|)} , 1 \big)$.  Since $C_1 > 4 \frac{|B|}{\theta}$, and $\exp(- \theta t) > 1 - \theta t$, it follows that for all $t \in [0,t_1]$, 
\begin{eqnarray*}
\exp(-\theta t) C_1 &\geq& \big(1 - \theta t) \frac{4 |B|}{\theta}
\\&\geq& \big(1 - \theta \times \frac{1}{4 \theta} \big) \times \frac{4 |B|}{\theta}\ \ \ \geq\ \ \ 3 \frac{|B|}{\theta};
\end{eqnarray*}
and
\begin{eqnarray*}
(\mu_1 + \mu_2) \big( |C_3| + |C_4| \big) t \exp(- \theta t) &\leq& (\mu_1 + \mu_2) \big( |C_3| + |C_4| \big) \times \frac{|B|}{4 \theta (\mu_1 + \mu_2)(|C_3| + |C_4|)}
\\&<& \frac{|B|}{\theta};
\end{eqnarray*}
from which we conclude that for all $t \in [0,t_1]$,
\begin{equation}\label{getonit1}
\BUF(t) > \frac{|B|}{\theta}.
\end{equation}
It is also easily verified from the definition of $\BUF$ and the fact that $C_1 > 0$ and $t \exp(- \theta t) \leq \frac{1}{\theta}$ for all $t \geq 0$ that
\begin{equation}\label{getonit1b}
\inf_{t \geq 0} \BUF(t) \geq - \frac{(\mu_1 + \mu_2) \big( |C_3| + |C_4| \big) + |B|}{\theta}.
\end{equation}
Recall that $\kappa^1, \kappa^2$ are two independent standard Brownian motions.  It follows from the well-known representation of a non-stationary OU process as a time-changed Brownian motion \cite{Ross.09} that we may construct $\OU_{0,a}, \OU_{0,b}, \kappa^1, \kappa^2, \Pi$ on a common probability space s.t. w.p.1, for all $t \geq 0$,
$$\OU_{0,a}(t) =  \big( \frac{p (1 - p)}{\mu_1 \mu_2} \big)^{\frac{1}{2}} \exp(- \mu_1 \mu_2 t) \kappa^1\big( \exp(2 \mu_1 \mu_2 t) - 1 \big);
$$
$$
\OU_{0,b}(t) =  \theta^{-\frac{1}{2}} \exp(- \theta t) \kappa^2\big( \exp(2 \theta t) - 1 \big);
$$
and $\Pi(t)$ equals
\begin{eqnarray}
\ &\ &\ (\mu_1 - \mu_2) \big( \frac{p (1 - p)}{\mu_1 \mu_2} \big)^{\frac{1}{2}} \exp\big( - (\mu_1 \mu_2 + \theta) t \big) \int_0^t \exp(\theta y) \kappa^1\big( \exp(2 \mu_1 \mu_2 y) - 1 \big) dy \label{getonit2}
\\&\ &\ \ \ \ +\ \theta^{-\frac{1}{2}} \exp(- \theta t) \kappa^2\big( \exp(2 \theta t) - 1 \big). \nonumber
\end{eqnarray}
W.l.o.g, suppose $\mu_1 > \mu_2$ (if not an identical argument works but with the roles of $\mu_1,\mu_2$ interchanged).  It follows from (\ref{getonit2}) that for $t \in [0,t_1]$, 
\begin{eqnarray*}
\Pi(t) &\geq& \bigg( \inf_{y \in [0,t_1]} \kappa^1\big( \exp(2 \mu_1 \mu_2 y) - 1 \big) \bigg) (\mu_1 - \mu_2) \big( \frac{p (1 - p)}{\mu_1 \mu_2} \big)^{\frac{1}{2}} \exp\big( - (\mu_1 \mu_2 + \theta) t \big) \int_0^t \exp(\theta y) dy 
\\&\ &\ \ \ +\ \ \ \theta^{-\frac{1}{2}} \exp(- \theta t) \bigg( \inf_{y \in [0,t_1]} \kappa^2\big( \exp(2 \theta y) - 1 \big) \bigg),
\end{eqnarray*}
itself at least 
$${\mathcal H}_1 \stackrel{\Delta}{=} \bigg(1 + (\mu_1 - \mu_2) \big( \frac{p (1 - p)}{\mu_1 \mu_2} \big)^{\frac{1}{2}} t_1 \bigg) \bigg( \inf_{y \in [0,t_1]} \kappa^1\big( \exp(2 \mu_1 \mu_2 y) - 1 \big) \bigg) + \big(1 + \theta^{-\frac{1}{2}} \big) \bigg( \inf_{y \in [0,t_1]} \kappa^2\big( \exp(2 \theta y) - 1 \big) \bigg).
$$
Let us also define
$${\mathcal H}_2 \stackrel{\Delta}{=} \min\bigg( \inf_{y \in [t_1,T]} \kappa^1\big( \exp(2 \mu_1 \mu_2 y) - 1 \big), \inf_{y \in [t_1,T]} \kappa^2\big( \exp(2 \theta y) - 1 \big) \bigg).$$
Here we recall that the basic properties of continuous Gaussian processes ensure that $\inf_{y \in [0,t_1]} \kappa^2\big( \exp(2 \theta y) - 1 \big) < 0, \inf_{y \in [0,t_1]} \kappa^1\big( \exp(2 \mu_1 \mu_2 y) - 1 \big) < 0$.  Note that for $\epsilon, M > 0$, the event $\lbrace {\mathcal H}_1 > - \epsilon, {\mathcal H}_2 > M \rbrace$ implies the event 
$\lbrace \inf_{t \in [0,t_1]} \Pi(t) > - \epsilon \rbrace$, as well as the event (for all $t \in [t_1,T]$)
\begin{eqnarray*}
\int_0^t \exp(\theta y) \kappa^1\big( \exp(2 \mu_1 \mu_2 y) - 1 \big) dy &\geq& - \int_0^{t_1} \exp(\theta y) \epsilon dy + \int_{t_1}^t \exp(\theta y) \kappa^1\big( \exp(2 \mu_1 \mu_2 y) - 1 \big) dy
\\&\geq& - t_1 \exp(\theta t_1) \epsilon + (t - t_1) M\ \ \ \geq\ \ \ - t_1 \exp(\theta t_1) \epsilon;
\end{eqnarray*}
and also the event (for all $t \in [t_1,T]$)
$$
\theta^{-\frac{1}{2}} \exp(- \theta t) \kappa^2\big( \exp(2 \theta t) - 1 \big) \geq \theta^{-\frac{1}{2}} \exp(- \theta T) M.
$$
Summarizing the above, and combining with (\ref{getonit1}) - (\ref{getonit2}), we conclude that for $\epsilon, M > 0$, the event $\lbrace {\mathcal H}_1 > - \epsilon, {\mathcal H}_2 > M \rbrace$ implies the event
\begin{equation}\label{event1lower}
\bigg\lbrace \inf_{t \in [0,t_1]} \big( \Pi(t) + \BUF(t) \big) \geq - \epsilon + \frac{|B|}{\theta} \bigg\rbrace,
\end{equation}
 as well as the event
\begin{eqnarray}
\ &\ &\ \Bigg\lbrace \inf_{t \in [t_1,T]} \big( \Pi(t) + \BUF(t) \big) \geq \theta^{-\frac{1}{2}} \exp(- \theta T) M \label{event2lower}
\\&\ &\ \ \ \ \ \ \ \ \ \ \ \ - (\mu_1 - \mu_2) \big( \frac{p (1 - p)}{\mu_1 \mu_2} \big)^{\frac{1}{2}} t_1 \exp(\theta t_1) \epsilon - \frac{(\mu_1 + \mu_2) \big( |C_3| + |C_4| \big) + |B|}{\theta}\bigg\rbrace. \nonumber
\end{eqnarray}
Setting $\epsilon = \frac{|B|}{2 \theta}$ and
$$M = 2 \theta^{\frac{1}{2}} \exp(\theta T) \bigg( (\mu_1 - \mu_2) \big( \frac{p (1 - p)}{\mu_1 \mu_2} \big)^{\frac{1}{2}} t_1 \exp(\theta t_1) \epsilon + \frac{(\mu_1 + \mu_2) \big( |C_3| + |C_4| \big) + |B|}{\theta} \bigg),$$
it follows that
$$\E\big[ I\bigg( \inf_{t \in [0,T]} \big( \Pi(t) + \BUF(t) \big) > 0 \bigg) \big] \geq P\big( \lbrace {\mathcal H}_1 > - \epsilon, {\mathcal H}_2 > M \rbrace \big).$$
That for this choice of $M$ and $\epsilon$ (and in fact any strictly positive choice for $M,\epsilon$) it holds that $P\big( \lbrace {\mathcal H}_1 > - \epsilon, {\mathcal H}_2 > M \rbrace \big) > 0$ follows from well-known properties of Brownian motion and associated standard arguments (see e.g. \cite{Borodin.12}), and we omit the details.  Combining the above completes the proof.  $\qed$
\endproof

\subsection{Proof of non-insensitivity.}
To prove that the true large deviations exponent $-\frac{\theta (\theta + \mu_1 \mu_2)}{2 (\theta + \mu_1 + \mu_2 - 1)}$ is not insensitive, it suffices to exhibit two hyper-exponential distributions, with respective parameters $p^a,\mu^a_1,\mu^a_2,p^b,\mu^b_1,\mu^b_2$, s.t. these two hyper-exponential distributions have the same first two moments (with first moment equal to one), yet there exists strictly positive $\theta < \min(\mu^a_1,\mu^a_2,\mu^b_1,\mu^b_2)$ s.t. 
$\frac{\theta (\theta + \mu^a_1 \mu^a_2)}{2 (\theta + \mu^a_1 + \mu^a_2 - 1)} \neq 
\frac{\theta (\theta + \mu^b_1 \mu^b_2)}{2 (\theta + \mu^b_1 + \mu^b_2 - 1)}$.  The existence of such a pair of hyper-exponential distributions is explicitly given in the following lemma.
\begin{lemma}
Let $S_a$ be a hyper-exponentially distributed r.v. with parameters $\mu^a_1 = \frac{1}{31} (15 - 2 \times \sqrt{2}), \mu^a_2 = 2 \times \sqrt{2}, p^a = \frac{1}{257} (121 - 32 \times \sqrt{2})$; and $S_b$ be a hyper-exponentially distributed r.v. with parameters $\mu^b_1 = \frac{1}{7} (3 - \sqrt{2}), \mu^b_2 = \sqrt{2}, p^b = \frac{1}{17} (7 - 4 \times \sqrt{2})$.  Let $\theta_0 = \frac{1}{5}$.  Then $\mu^a_1, \mu^a_2, \mu^b_1, \mu^b_2 \in (0,\infty); \theta_0 \in \bigg(0, \min\big( \mu^a_1, \mu^a_2, \mu^b_1, \mu^b_2 \big) \bigg)$; $\E[S_a] = \E[S_b] = 1$; $\E[S^2_a] = \E[S^2_b] = 4$; but $\frac{\theta_0 (\theta_0 + \mu^a_1 \mu^a_2)}{2(\theta_0 + \mu^a_1 + \mu^a_2 - 1)} \neq \frac{\theta_0 (\theta_0 + \mu^b_1 \mu^b_2)}{2(\theta_0 + \mu^b_1 + \mu^b_2 - 1)}$.
\end{lemma}
\proof{Proof:} 
That $\mu^a_1, \mu^a_2, \mu^b_1, \mu^b_2 \in (0,\infty)$, and $p^a, p^b \in (0,1)$ is easily verified.  That $\mu^a_1 > \theta_0$ follows from the fact that $\frac{1}{31} \times (15 - 2 \times \sqrt{2}) > \frac{1}{36} \times (15 - 4) = \frac{1}{4} > \frac{1}{5}$.  That $\mu^b_1 > \theta_0$ follows from the fact that $\frac{1}{7} \times (3 - \sqrt{2}) > \frac{1}{7} \times (3 - 1.5) = \frac{3}{14} > \frac{1}{5}$.  That $\theta_0 < \min(\mu^a_2, \mu^b_2)$ is trivial.  Furthermore, we claim that $\E[S_a] = \E[S_b] = 1$ and $\E[S^2_a] = \E[S^2_b] = 4$.  Indeed, a straightforward calculation yields the following:
\begin{eqnarray*}
\E[S_a] &=& \frac{ \frac{1}{257} (121 - 32 \times \sqrt{2}) }{\frac{1}{31} (15 - 2 \times \sqrt{2})} + \frac{1 - \frac{1}{257} (121 - 32 \times \sqrt{2})}{2 \times \sqrt{2}}
\\&=& \frac{2 \times \sqrt{2} \times \frac{1}{257} (121 - 32 \times \sqrt{2}) + \frac{1}{31} (15 - 2 \times \sqrt{2}) \times \bigg(1 - \frac{1}{257} (121 - 32 \times \sqrt{2})\bigg)}{\frac{1}{31} (15 - 2 \times \sqrt{2}) \times 2 \times \sqrt{2}}
\\&=& \frac{62 \times \sqrt{2} \times (121 - 32 \times \sqrt{2}) + (15 - 2 \times \sqrt{2}) \times \bigg(136 + 32 \times \sqrt{2})\bigg)}{257 \times (15 - 2 \times \sqrt{2}) \times 2 \times \sqrt{2}}
\\&=& \frac{ \big( - 62 \times 32 \times 2 + 15 \times 136 - 32 \times 2 \times 2 \big) + \sqrt{2} \times \big( 62 \times 121 + 15 \times 32 - 2 \times 136 \big) }{- 257 \times 8 + 514 \times 15 \times \sqrt{2}}
\\&=& \frac{ -2056 + 7710 \times \sqrt{2}}{ - 2056 + 7710 \times \sqrt{2} }\ \ \ =\ \ \ 1.
\end{eqnarray*}
\begin{eqnarray*}
\E[S^2_a] &=& 2 \times \frac{ \frac{1}{257} (121 - 32 \times \sqrt{2}) }{\big(\frac{1}{31} (15 - 2 \times \sqrt{2})\big)^2} + 2 \times \frac{1 - \frac{1}{257} (121 - 32 \times \sqrt{2})}{\big(2 \times \sqrt{2}\big)^2}
\\&=& \frac{2 \times \frac{1}{257} (121 - 32 \times \sqrt{2}) \times \big(2 \times \sqrt{2}\big)^2 
+ 2 \times \big(\frac{1}{31} (15 - 2 \times \sqrt{2})\big)^2 \times \big( 1 - \frac{1}{257} (121 - 32 \times \sqrt{2}) \big)}{\big(\frac{1}{31} (15 - 2 \times \sqrt{2})\big)^2 \times \big(2 \times \sqrt{2}\big)^2}
\\&=& \frac{16 \times 31^2 \times (121 - 32 \times \sqrt{2}) 
+ 2 \times (15 - 2 \times \sqrt{2})^2 \times \big( 257 - (121 - 32 \times \sqrt{2}) \big)}{257 \times 8 \times (15 - 2 \times \sqrt{2})^2}
\\&=& \frac{16 \times 31^2 \times (121 - 32 \times \sqrt{2}) 
+ (466 - 120 \times \sqrt{2}) \times \big( 136 + 32 \times \sqrt{2} \big)}{479048 - 123360 \times \sqrt{2}}
\\&=& \frac{
\bigg( 16 \times 31^2 \times 121 + 466 \times 136 - 120 \times 32 \times 2 \bigg) + \bigg(
-16 \times 31^2 \times 32 + 32 \times 466 - 120 \times 136\bigg) \times \sqrt{2}}{479048 - 123360 \times \sqrt{2}}
\\&=& \frac{1916192 - 493440 \times \sqrt{2}}{479048 - 123360 \times \sqrt{2}}\ \ \ =\ \ \ 4;
\end{eqnarray*}
\begin{eqnarray*}
\E[S_b] &=& \frac{ \frac{1}{17} (7 - 4 \times \sqrt{2}) }{\frac{1}{7} (3 - \sqrt{2})} + \frac{1 - \frac{1}{17} (7 - 4 \times \sqrt{2})}{\sqrt{2}}
\\&=& \frac{ \frac{1}{17} \times (7 - 4 \sqrt{2}) \times \sqrt{2} + \big( 1 - \frac{1}{17} (7 - 4 \times \sqrt{2}) \big) \times \frac{1}{7} \times (3 - \sqrt{2}) }{\frac{1}{7} \times (3 - \sqrt{2}) \times \sqrt{2}} 
\\&=& \frac{ \big(7 \times 7 \times \sqrt{2} - 7 \times 4 \times 2 \big) + (10 + 4 \sqrt{2}) \times (3 - \sqrt{2})}{17 \times 3 \times \sqrt{2} - 17 \times 2}
\\&=& \frac{ \big( -7 \times 4 \times 2 + 10 \times 3 - 8 \big) + \big(7 \times 7 - 10 + 12 \big) \sqrt{2} }{-34 + 51 \sqrt{2}}\ \ \ =\ \ \ 1;
\end{eqnarray*}
\begin{eqnarray*}
\E[S^2_b] &=& 2 \times \frac{ \frac{1}{17} (7 - 4 \times \sqrt{2}) }{\big(\frac{1}{7} (3 - \sqrt{2})\big)^2} + 2 \times \frac{1 - \frac{1}{17} (7 - 4 \times \sqrt{2})}{\big(\sqrt{2}\big)^2}
\\&=& \frac{ 2 \times \frac{1}{17} \times (7 - 4 \sqrt{2}) + \big( 1 - \frac{1}{17}(7 - 4 \sqrt{2}) \big)\big(\frac{1}{7} (3 - \sqrt{2})\big)^2 }{\big(\frac{1}{7} (3 - \sqrt{2})\big)^2}
\\&=& \frac{ \big( 98 \times 7 - 98 \times 4 \times \sqrt{2} \big) + (10 + 4\sqrt{2}) \times (11 - 6 \sqrt{2})}{17 \times (11 - 6 \sqrt{2})}
\\&=& \frac{ \big( 98 \times 7 + 10 \times 11 - 6 \times 4 \times 2 \big) + \big( -98 \times 4 - 6 \times 10 + 4 \times 11 \big) \sqrt{2} }{17 \times 11 - 17 \times 6 \times \sqrt{2}}\ \ \ =\ \ \ 4.
\end{eqnarray*}
Next, we show that although $S_a$ and $S_b$ have the same first two moments, they lead to different large deviations exponents.  Indeed, note that it suffices to demonstrate that
$$
(\theta_0 + \mu^a_1 \mu^a_2) \times (\theta_0 + \mu^b_1 + \mu^b_2 - 1) \neq 
(\theta_0 + \mu^b_1 \mu^b_2) \times (\theta_0 + \mu^a_1 + \mu^a_2 - 1).$$
Then a straightforward calculation yields
\begin{eqnarray*}
\ &\ &\ (\theta_0 + \mu^a_1 \mu^a_2) \times (\theta_0 + \mu^b_1 + \mu^b_2 - 1) - (\theta_0 + \mu^b_1 \mu^b_2) \times (\theta_0 + \mu^a_1 + \mu^a_2 - 1)
\\&\ &\ \ \ =\ \ (\frac{1}{5} + \frac{30}{31} \sqrt{2} - \frac{8}{31}) \times (\frac{1}{5} + \frac{3}{7} + \frac{6}{7} \sqrt{2} - 1) 
\\&\ &\ \ \ \ \ \ \ \ \ \ -\ \ (\frac{1}{5} + \frac{3}{7} \sqrt{2} - \frac{2}{7}) \times \big(\frac{1}{5} + \frac{15}{31} + (2 - \frac{2}{31}) \sqrt{2} - 1\big) 
\\&\ &\ \ \ = (5 \times 7 \times 31)^{-1} \times \bigg( \big( -63 + 1050 \sqrt{2} \big) \times \big( -403 + 930 \sqrt{2} \big) - \big( -93 + 465 \sqrt{2} \big) \times \big( -343 + 2100 \sqrt{2} \big) \bigg)
\\&\ &\ \ \ = - (5 \times 7 \times 31)^{-1} \times (6510 + 126945 \sqrt{2})\ \ \ \neq 0, 
\end{eqnarray*}
completing the proof.  $\qed$
\subsection{Proof that $- \frac{\theta( \theta + \mu_1 \mu_2 ) }{2( \theta + \mu_1 + \mu_2 - 1)} < - \frac{\mu_1 \mu_2 \theta}{2 (\mu_1 + \mu_2 - 1)}.$}
Since $\E[S] = 1$, and hence $\mu_1 + \mu_2 - 1 > 0$, to prove that $- \frac{\theta( \theta + \mu_1 \mu_2 ) }{2( \theta + \mu_1 + \mu_2 - 1)} < - \frac{\mu_1 \mu_2 \theta}{2 (\mu_1 + \mu_2 - 1)}$, it suffices to demonstrate that 
$$( \theta + \mu_1 \mu_2 ) (\mu_1 + \mu_2 - 1) - ( \theta + \mu_1 + \mu_2 - 1) \mu_1 \mu_2 > 0,$$
equivalently that $\theta (\mu_1 + \mu_2 - 1) - \theta \mu_1 \mu_2 > 0$.  It thus suffices to demonstrate that $\mu_1 + \mu_2 - 1 > \mu_1 \mu_2$, equivalently that $\mu_1 (1 - \mu_2) > 1 - \mu_2$.  First, suppose $\mu_2 < 1$.  Then the fact that $\E[S] = 1$ implies that $\mu_1 > 1$, and hence by the implied non-negativity of $1 - \mu_2$ we conclude that $\mu_1 (1 - \mu_2) > 1 - \mu_2$, completing the proof.  Alternatively, suppose that $\mu_2 > 1$, in which case $\mu_1 < 1$.  It then follows from a nearly identical argument that $\mu_1 (\mu_2 - 1) < \mu_2 - 1$, completing the proof in this case as well.  $\qed$  

\bibliographystyle{plainnat}

\end{document}